\documentclass[smallextended]{svjour3}
\usepackage[margin=1in]{geometry}
\usepackage{siunitx}
\usepackage{amsfonts,amsmath,amssymb}
\usepackage{graphicx}
\usepackage[hidelinks]{hyperref}
\usepackage{url}
\usepackage{float}

\usepackage{todonotes}

\renewcommand{\b}{\mathbf}
\renewcommand{\L}{\mathcal{L}}

\usepackage{subcaption}

\renewcommand{\L}{\mathcal{L}}
\renewcommand{\b}{\boldsymbol}

\usepackage{color,soul}
\newcommand{\diff}[1]{#1}

\title{Discretization of non-uniform rational B-spline (NURBS) models for meshless isogeometric analysis}

\author{Urban Duh \and Varun Shankar \and Gregor Kosec}

\institute{Urban Duh 
	\at Faculty of Mathematics and Physics,
	University of Ljubljana, Jadranska 19, 1000 Ljubljana, Slovenia
	\\\email{urban.duh@fmf.uni-lj.si}
	\and
	Varun Shankar
	\at  School of Computing, University of Utah, Salt Lake City, UT USA 84112 
	\\\email{shankar@cs.utah.edu}
	\and
	Gregor Kosec
	\at ``Jožef Stefan'' Institute, Department E6, Parallel
	and Distributed Systems Laboratory, Jamova cesta 39, 1000 Ljubljana,
	Slovenia
	\\\email{gregor.kosec@ijs.si}
}

\hypersetup{
	pdftitle={Discretization of non-uniform rational B-spline (NURBS) models for meshless isogeometric analysis},
	pdfauthor={U. Duh, V. Shankar and G. Kosec}
}

\providecommand{\keywords}[1]{\textbf{Key words:} #1}

\begin{document}
	
\maketitle

\begin{abstract}
We present an algorithm for fast generation of quasi-uniform and 
\diff{variable-spacing} nodes on domains whose boundaries are represented as
computer-aided design (CAD) models, more specifically non-uniform rational
B-splines (NURBS). This new algorithm enables the solution of partial
differential equations (PDEs) within the volumes enclosed by these CAD models
using (collocation-based) meshless numerical discretizations. Our hierarchical
algorithm first generates quasi-uniform node sets directly on the NURBS surfaces
representing the domain boundary, then uses the NURBS representation in
conjunction with the surface nodes to generate nodes within the volume enclosed
by the NURBS surface. We provide evidence for the quality of these node sets by
analyzing them in terms of local regularity and separation distances. Finally,
we demonstrate that these node sets are well-suited (both in terms of accuracy
and numerical stability) for meshless radial basis function generated finite differences
(RBF-FD) discretizations of the Poisson, Navier-Cauchy, and heat equations. Our
algorithm constitutes an important step in bridging the field of node
generation for meshless discretizations with isogeometric
analysis.
\end{abstract}

\keywords{meshless, CAD, RBF-FD, advancing front algorithms, NURBS}
	
\section{Introduction} 
\label{sec:intro}
A key element of any numerical method for solving partial differential
equations (PDE) is discretization of the domain. In traditional numerical
methods such as the finite element method (FEM), this discretization is
typically performed by partitioning the domain into a mesh, \emph{i.e.},
a finite number of elements that entirely cover it. Despite substantial
developments in the field of mesh generation, the process of meshing
often remains the most time consuming part of the whole solution procedure
while the mesh quality limits the accuracy and stability of the numerical
solution~\cite{liu_introduction_2005}. In contrast, meshless methods
for PDEs work directly on point clouds; in this context, points are
typically referred to as ``nodes''. In particular, meshless methods based
on radial basis function generated finite difference (RBF-FD) formulas
allow for high-order accurate numerical solutions of PDEs on complicated
time-varying domains~\cite{shankar2021efficient,jancic_monomial_2021}
and even manifolds~\cite{petras2018rbf}. The generation of suitable
nodes is an area of ongoing research, with much work in recent
years~\cite{slak2019generation,duh2021fast,fornberg2015fast,van2021fast,shankar2018robust}. In this work, we focus primarily on the generation of nodes
suitable for RBF-FD discretizations, although our node generation approach is
fully independent of the numerical method used. Node sets may be generated
in several different ways. For instance, one could simply generate a mesh using
an existing tool and discard the connectivity information~\cite{liu2002mesh}.
However, such an approach is obviously computationally expensive, not easily
generalized to higher dimensions, and in some scenarios even fails to generate
node distributions of sufficient quality~\cite{shankar2018robust}. Another
possible approach is to use randomly-generated nodes~\cite{liu2002mesh,milewski2021higher}; this approach has been used (with some modifications)
in areas such as compressive sensing~\cite{adcock2021improved} and function
approximation in high dimensions~\cite{narayan2012stochastic}. Other
approaches include iterative optimization~\cite{hardin2004discretizing,liu2010node,kosec_local_2018}, sphere-packing~\cite{li2000point}, QR
factorization~\cite{suchde2022point}, and repulsion~\cite{fornberg2015fast,van2021fast}. It is generally accepted that quasi-uniformly-spaced node sets
improve the stability of meshless methods~\cite{wendland2004scattered}. In this
context, methods based on Poisson disk sampling are particularly appealing as
they produce quasi-uniformly spaced nodes, scale to arbitrary dimension, are
computationally efficient, and can be fully automated~\cite{shankar2018robust,slak2019generation,duh2021fast}. A related consideration is the quality
of the domain discretization. In the context of meshes, for instance, it is
common to characterize mesh quality using element aspect ratios or determinants
of Jacobians~\cite{gokhale2008practical,zala2018curvilinear}. Analogously,
the node generation literature commonly characterizes node quality in
terms of two measures: the minimal spacing between any pair of nodes (the
separation distance), and the maximal empty space without nodes (the fill
distance). Once again, in this context, Poisson disk sampling via advancing
front methods constitutes the state of the art~\cite{slak2019generation,duh2021fast}. More specifically, the DIVG algorithm~\cite{slak2019generation}
allows for variable \diff{spacing} Poisson disk sampling on complicated domains
in arbitrary dimension, while its generalization (sDIVG)~\cite{duh2021fast} allows for
sampling of arbitrary-dimensional parametric surfaces. DIVG has since
been parallelized~\cite{depolli2022parallel}, distributed as a standalone
node generator~\cite{standalone}, and is also an important component
of the open-source meshless project Medusa~\cite{slak2021medusa}. Despite these rapid advances in node generation for meshless methods
(and in meshless methods themselves), the generation of node sets
on domains whose boundaries are specified by computer-aided design
(CAD) models is still in its infancy. Consequently, the application of
meshless methods in CAD supplied geometries is rare and limited either
to smooth geometries~\cite{mirfatah2019solution} or to the use of surface
meshes~\cite{drumm2008finite,gerace2017model,jacquemin2022smart}. In
contrast, mesh generation and the use of FEM in CAD geometries is a mature and
well-understood field~\cite{gokhale2008practical,cottrell2009isogeometric}.
In our experience, current node generation approaches on CAD geometries
violate quasi-uniformity near the boundaries and are insufficiently robust
or automated for practical use in engineering applications. 

In this work,
we extend the sDIVG method to the generation of variable \diff{spacing} node sets on
parametric CAD surfaces specified by non-uniform rational B-splines (NURBS). We
then utilize the \diff{variable-spacing} node sets generated by sDIVG in conjunction
with the DIVG method to generate node sets in the volume enclosed by the
NURBS surface. Our new framework is automated, computationally efficient,
scalable to higher dimensions, and generates node sets that retain 
quasi-uniformity all the way up to the boundary. This framework also inherits the
quality guarantees of DIVG and sDIVG, and is consequently well-suited for
stable RBF-FD discretizations of PDEs on complicated domain geometries. 

The remainder of the paper is organized as follows. The NURBS-DIVG
algorithm is presented in Section~\ref{sec:algorithm} along with analysis
of specific components of the algorithm. The quality of generated nodes is
discussed in Section~\ref{sec:node_qulity}. Its application to the RBF-FD
solution of PDEs is shown in Section~\ref{sec:PDE}. The paper concludes in
Section~\ref{sec:conc}.

\section{The NURBS-DIVG Algorithm}
\label{sec:algorithm}
CAD surfaces are typically described as a union of multiple, non-overlapping, parametric patches (curves in 2D, surfaces in 3D), positioned so that the transitions between them are either smooth or satisfying some geometric conditions. A popular choice for representing each patch is a NURBS~\cite{piegl_nurbs_2012}, which is the focus of our work. Here, we present a NURBS-DIVG algorithm that has three primary components:
\begin{enumerate}
\item First, we extend the sDIVG algorithm~\cite{duh2021fast} (Section \ref{sec:sdivg}) for sampling parametric surfaces to sampling individual NURBS patches and also the union of multiple NURBS patches (Section \ref{sec:nurbs-patches}). 
\item Next, we deploy the DIVG algorithm~\cite{slak2019generation,depolli2022parallel} in the interior of the domain using the sDIVG generated samples as seed nodes (Section \ref{sec:divg}).
\item To ensure that DIVG generates the correct node sets in the interior of the domain whose boundary consists of multiple parametric NURBS patches, we augment sDIVG with a supersampling parameter (Section \ref{sec:inout}).
\end{enumerate}
\diff{In the following subsections}, we first briefly present the DIVG algorithm. We then describe the sDIVG algorithm, which generalizes DIVG to parametric surfaces, focusing on sampling a single NURBS surface. We then describe how the sDIVG algorithm is generalized to surfaces consisting of multiple NURBS patches, each of which have their own boundary curves. Finally, we describe the inside check utilized by our algorithm needed to generate nodes within the NURBS patches, and the complications therein.

\subsection{The DIVG algorithm}
\label{sec:divg}
We now describe the DIVG algorithm for generating node sets within an arbitrary domain. As mentioned previously, this algorithm forms the foundation of the sDIVG and NURBS-DIVG algorithms. 

DIVG is an iterative algorithm that begins with a given set of nodes called
``seed nodes''; in our case, these will later be provided by the sDIVG part of the NURBS-DIVG. The
seed nodes are placed in an \emph{expansion queue}. In each iteration $i$ of
the DIVG algorithm, a single node $\b{p}_i$ is dequeued and ``expanded''. Here,
``expansion'' means that a set $C_i$ of $n$ candidates for new nodes is
uniformly generated on a sphere centered at the node $\b{p}_i$, with some
radius $r_i$ and a random rotation. Here, $r_i$ stands for target nodal
\diff{spacing} and can be thought of as derived from a \diff{spacing} function
$h$, so that $r_i = h(\b{p}_i)$~\cite{slak2019generation,slak_adaptive_2019}.
Of course, the set $C_i$ may contain candidates that lie outside the domain
boundary or are too close to an existing node. Such candidates are rejected.
The candidates that are not rejected are simply added to the domain and to the
expansion queue; this is illustrated in Figure~\ref{fig:fill_scheme}. The
iteration continues until the queue is empty. A full description of \diff{the
DIVG algorithm} can be found in~\cite{slak2019generation}. Its parallel variant
is described in~\cite{depolli2022parallel}.

\subsection{The sDIVG algorithm}
\label{sec:sdivg}
The sDIVG algorithm is a generalization of \diff{the DIVG algorithm} to parametric surfaces. Unlike
DIVG (which fills volumes with node sets), sDIVG instead places nodes on a
target parametric surface in such a way that the spacing between nodes on the
surface follows a supplied \diff{spacing} function. While other algorithms typically
achieve this through direct Cartesian sampling and
elimination~\cite{yuksel2015,shankar2018robust}, the sDIVG algorithm samples
the parametric domain corresponding to the surface with an
appropriately-transformed version of the supplied \diff{spacing} function. More
concretely, given a domain $\Omega \subset \mathbb{R}^d$, sDIVG iteratively
samples its boundary $\partial \Omega \subset \mathbb{R}^{d}$ by sampling a
parametrization $\Lambda$ of its boundary instead. The advantage of this
approach over direct Cartesian sampling is obtained from the fact that $\Lambda
\subset \mathbb{R}^{d-1}$ (or $\mathbb{S}^{d-1}$) is a lower-dimensional
representation of $\partial \Omega$, leading to an increase in efficiency.

\begin{figure}[t]
	\centering
	\includegraphics[height=0.37\linewidth]{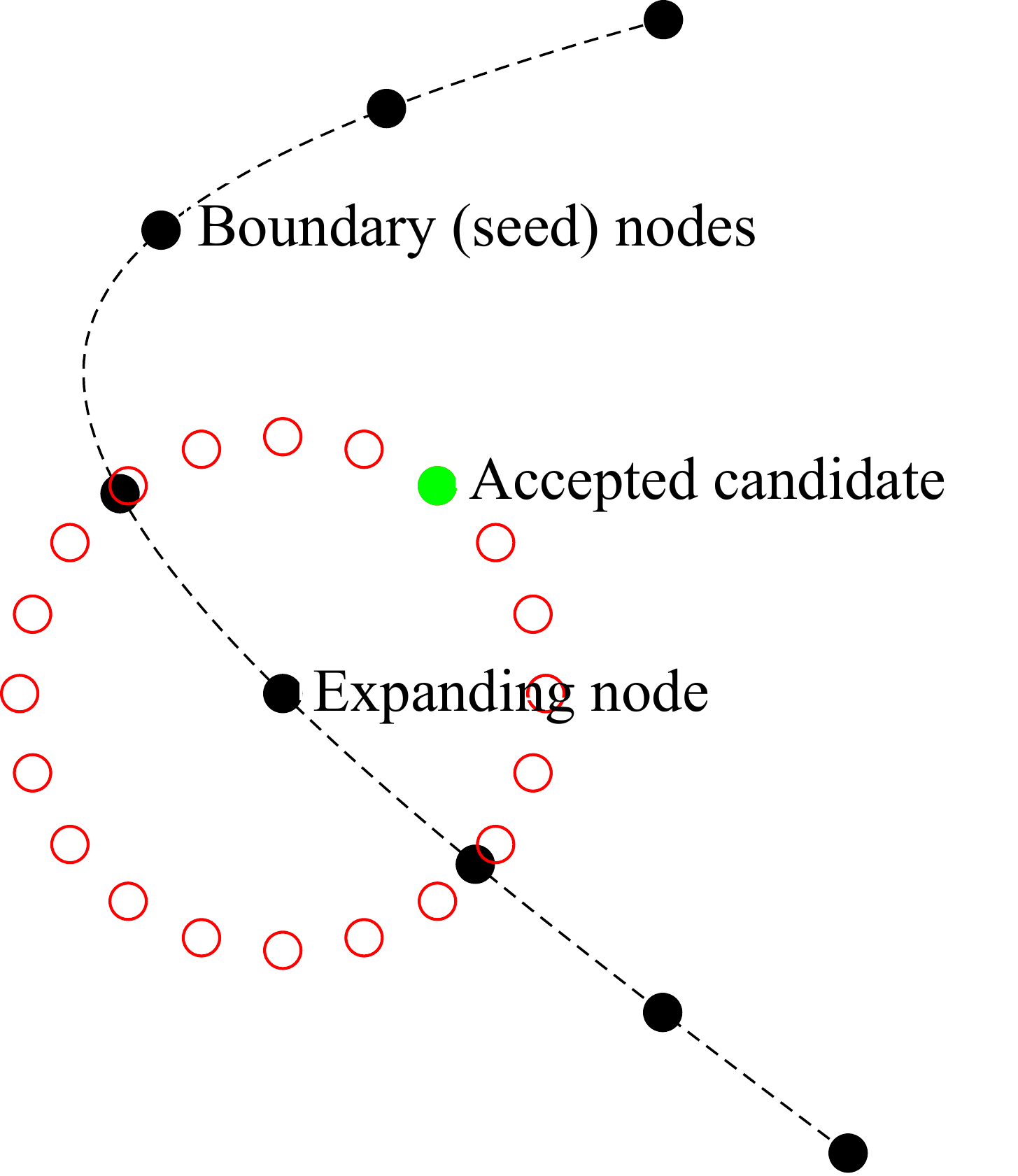}
	\includegraphics[height=0.32\linewidth]{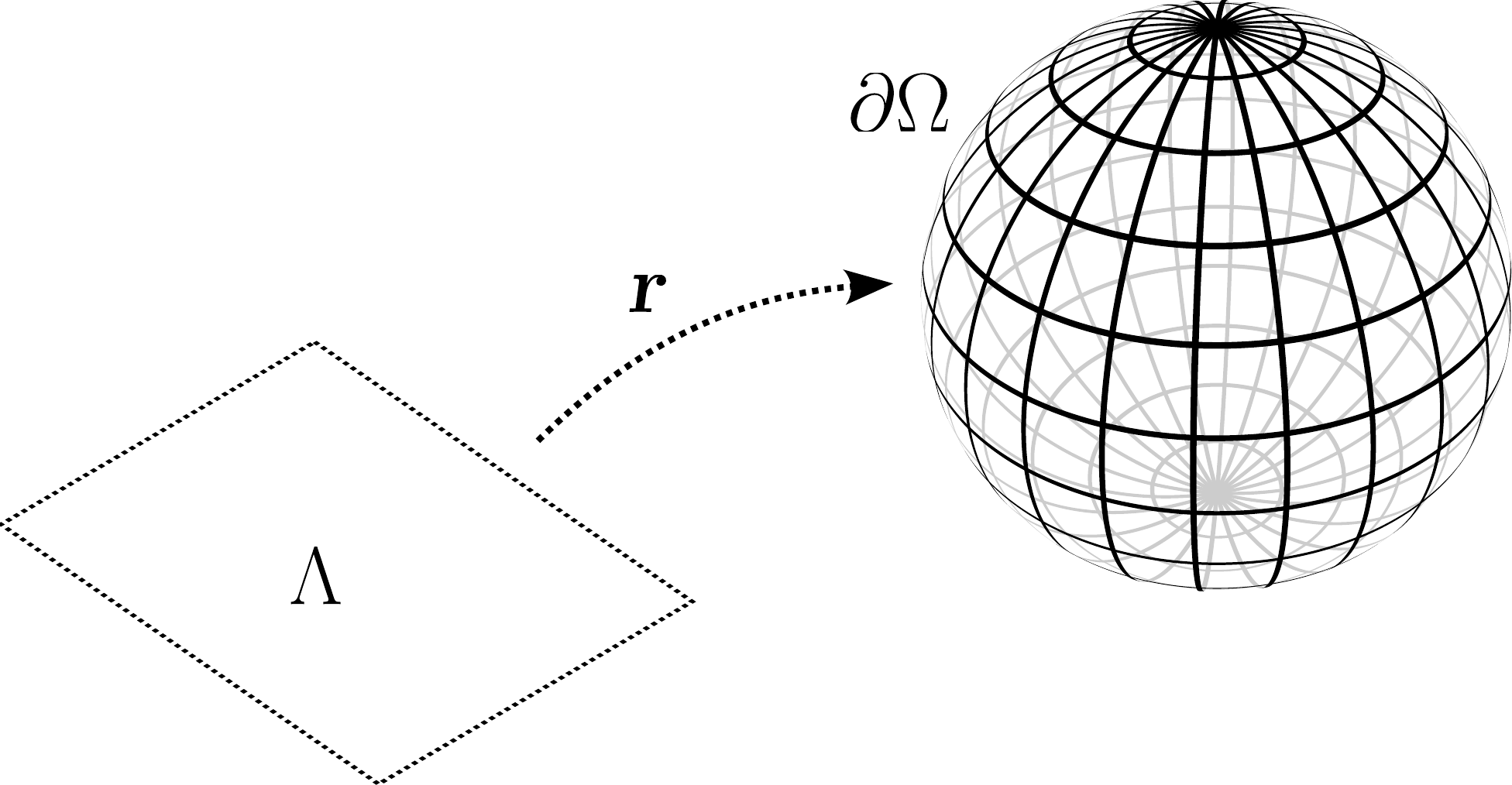}
	\caption{The DIVG expansion scheme (left) and the sDIVG mapping scheme (right).}
	\label{fig:fill_scheme}
\end{figure}

We now briefly describe the \diff{spacing} function transformation utilized by sDIVG
\diff{to generate a candidate set for expansion analogous to the one in DIVG}.
We first define a parametrization $\mathbf{r} : \Lambda \to \partial \Omega$,
\emph{i.e.}, a map from the parametric domain $\Lambda \subset
\mathbb{R}^{d-1}$ to the manifold $\partial\Omega \subset \mathbb{R}^d$; the
Jacobian of this function is denoted by $\nabla \mathbf{r}$. As in the DIVG
algorithm, let $h$ denote the desired \diff{spacing} function. Now, given a node
$\boldsymbol{\lambda}_i \in \Lambda$, we wish to generate a set of $n$ 
candidates for expanding $\boldsymbol{\lambda}_i$, which we write as
\begin{align}\label{eq:candset-sdivg}
	C_i = \{\boldsymbol \eta_{i, j} \in \Lambda;\ j = 1, \ldots, n\}.
\end{align}
It is important to note that the candidates $\boldsymbol \eta_{i, j}$ all lie
in the parametric domain. Our goal is to determine how far from $\boldsymbol
\lambda_i$ must each candidate lie. From the definition of $\mathbf{r}$ and
$h$, \diff{the target spacing between the candidate $\boldsymbol \eta_{i, j}$
and the node being expanded $\boldsymbol \lambda_i$ is}
\begin{align}\label{eq:r_to_h}
\|\mathbf r(\boldsymbol \eta_{i, j}) - \mathbf r(\boldsymbol 	\lambda_i)\| = h(\mathbf r(\boldsymbol \lambda_i)), 
\end{align}
for all $j = 1,\ldots,n$. The candidates $\boldsymbol \eta_{i, j}$ can be thought of as lying on some manifold around $\boldsymbol \lambda_i$. This allows us to rewrite $\boldsymbol \eta_{i, j}$ as
\begin{equation}\label{eq:eta_to_lambda}
\boldsymbol \eta_{i, j} = \boldsymbol \lambda_i + \alpha_{i, j} 		\vec s_{i, j},
\end{equation}
for some constant $\alpha_{i, j} > 0$ and unit vector $\vec s_{i, j}$. Here, \diff{we must appropriately
choose the unit vectors $\vec s_{i, j}$ (more on that later)} and
$\alpha_{i,j}$ must be determined by an appropriate transformation of
$h(\mathbf r(\boldsymbol \lambda_i))$, \emph{i.e.} the parametric distances \diff{$\alpha_{i, j}$ between
the candidate and the node being expanded} must be
obtained by a transformation of the \diff{spacing} function $h$ specified on
$\partial\Omega$. We may now use Eq.~\eqref{eq:eta_to_lambda} to Taylor expand
$\mathbf{r}(\boldsymbol \eta_{i, j})$ as
\begin{equation}\label{eq:taylor1}
\mathbf r(\boldsymbol \eta_{i, j}) = \mathbf r(\boldsymbol \lambda_i + \alpha_{i, j} \vec s_{i, j}) \approx \mathbf r(\boldsymbol \lambda_i) + \alpha_{i, j} \nabla \mathbf r(\boldsymbol \lambda_i) \vec s_{i, j}.
\end{equation}
We can now use the \diff{Taylor expansion in Eq.~{\eqref{eq:taylor1}} to approximate the actual
spacing between $\boldsymbol \lambda_i$ and $\boldsymbol \eta_{i, j}$ in Eq.~{\eqref{eq:r_to_h}}} to obtain the following expression for $h(\mathbf r(\boldsymbol \lambda_i))$ in terms of $\alpha_{i,j}$:
\begin{equation}
h(\mathbf r(\boldsymbol \lambda_i)) \approx \|\mathbf r(\boldsymbol 	\lambda_i) + \alpha_{i, j} \nabla \mathbf r(\boldsymbol \lambda_i) \vec s_{i, j} - \mathbf r(\boldsymbol \lambda_i)\| = \alpha_{i, j} \|\nabla \mathbf r(\boldsymbol \lambda_i) \vec s_{i, j}\|.
\end{equation}
This in turn allows us to express $\alpha_{i,j}$ as
\begin{equation}\label{eq:alpha_to_h}
\alpha_{i, j} = \frac{h(\mathbf r(\boldsymbol \lambda_i))}{\|\nabla \mathbf r(\boldsymbol \lambda_i) \vec s_{i, j}\|}.
\end{equation}
\diff{It is important to note here that for such $\alpha_{i, j}$ the target
  spacing defined in Eq.~{\eqref{eq:r_to_h}} holds only approximately, \emph{i.e.}, to the
  first order in the Taylor series expanded in $\alpha_{i, j}$. This is not an
  issue in practice, since in order to solve PDEs, we typically require node
spacings that are small compared to the curvature of the domain boundary
$\partial\Omega$. Higher-order approximations can also be computed if needed.}
We can now use Eq.~\eqref{eq:alpha_to_h} within Eq.~\eqref{eq:candset-sdivg} to obtain an explicit expression for the candidate set $C_i$ purely in terms of the \diff{spacing} function $h$ and the parametrization $\mathbf{r}$. Thus, we have
\begin{equation}
C_i = \left\{\boldsymbol \lambda_i + \frac{h(\mathbf r(\boldsymbol \lambda_i))}{\|\nabla \mathbf r(\boldsymbol \lambda_i) \vec s_{i, j}\|} \vec s_{i, j}; \vec s_{i, j} \in S_i \right\},
\end{equation}
where $S_i$ is set of $n$ random unit vectors on a unit ball. All other steps are identical to the DIVG algorithm, albeit in the parametric domain $\Lambda$. The final set of points on $\partial\Omega$ is then obtained by evaluating the function $\mathbf{r}$ at these parametric samples; a schematic of this is shown in Figure \ref{fig:fill_scheme} (right). A full description of the sDIVG algorithm \diff{and
an analysis of its potential weakness} can be found in~\cite{duh2021fast}.

\subsection{NURBS-DIVG}
\label{sec:nurbs-divg}
In principle, sDIVG can be used with any map $\mathbf{r} : \Lambda \to \partial \Omega$. NURBS-DIVG, however, is a specialization of sDIVG to surfaces comprised of NURBS patches (collections of non-overlapping and abutting NURBS). We first describe the use of NURBS for generating the surface representation $\mathbf{r}$, then discuss the generalization to a surface containing multiple patches.

\subsubsection{An overview of NURBS surfaces}
\label{sec:nurbs}
To define a NURBS representation, it is useful to first define the corresponding B-spline basis in one-dimension. Given a sequence of nondecreasing real numbers $T = \{t_0, t_1, \dots, t_k\}$ called the \emph{knot vector}, the degree-$p$ B-spline basis functions $N_{i, p}(u)$ are defined recursively as~\cite{piegl_nurbs_2012}
\begin{align}
    N_{i, 0}(u) &= \begin{cases}
        1; & t_i \leq u < t_{i + 1} \\
        0; & \text{otherwise}
    \end{cases} \\
        N_{i, p}(u) &= \frac{u - t_i}{t_{i + p} - t_i} N_{i, p - 1}(u) +
        \frac{t_{i + p + 1} - u}{t_{i + p + 1} - t_{i + 1}} N_{i + 1, p - 1}(u)
\end{align}
We can now use these basis function to define a \emph{NURBS curve} in $\mathbb{R}^d$. \diff{Given a knot vector of the form $T = \{\overbrace{a, \dots, a}^{p + 1 \
\text{times}}, t_{p + 2}, \dots, t_{k - p - 1}, \overbrace{b, \dots, b}^{p + 1 \
\text{times}}\}$}, $n$ control points $\mathbf{p}_i \in \mathbb{R}^d$ and $n$
weights $w_i \in \mathbb{R}$, the degree-$p$ NURBS curve is defined as~\cite{piegl_nurbs_2012}
\begin{align}
    \mathbf{s}(u) = \frac{\sum_{i = 0}^{n - 1} N_{i, p}(u) w_i \mathbf{p}_i}{\sum_{i =
    0}^{n - 1} N_{i, p}(u) w_i}, \qquad \text{for} \ a \leq u \leq b.
\end{align}
In practice, it is convenient to evaluate $\mathbf{s}(u) \subset \mathbb{R}^d$
as a  B-spline curve in $\mathbb{R}^{d+1}$, and then project that curve down to
$\mathbb{R}^d$. For more details, see~\cite{piegl_nurbs_2012}. We evaluate the
B-spline curve in a numerically stable and efficient fashion using the de Boor algorithm~\cite{de_boor_1971}, which is itself a generalization of the well-known de Casteljau algorithm for Bezier curves~\cite{farin_cagd_2000}. It is also important to note that derivatives of the NURBS curves with respect to $u$ are needed for computing surface quantities (tangents and normals, for instance, or curvature). Fortunately, it is well-known that these parametric derivatives are also NURBS curves and can therefore also be evaluated using the de Boor algorithm~\cite{piegl_nurbs_2012}. We adopt this approach in this work.

Finally, the map $\mathbf r$ which represents a surface of co-dimension one in $\mathbb{R}^d$ can be represented as a NURBS surface. \diff{This surface is obtained in a fairly standard fashion as a tensor-product in knot space, followed by evaluation of the product space through the 1D spline maps}. This allows all NURBS surface operations to be computed via the de Boor algorithm applied to each parametric dimension. The sDIVG algorithm can then be used to sample this surface as desired, which in turn allows for node generation in the interior of the domain via DIVG. 

\subsubsection{Sampling surfaces consisting of multiple NURBS patches}
\label{sec:nurbs-patches}
\begin{figure}[!tbp]
	\centering
	\includegraphics[width=0.47\linewidth]{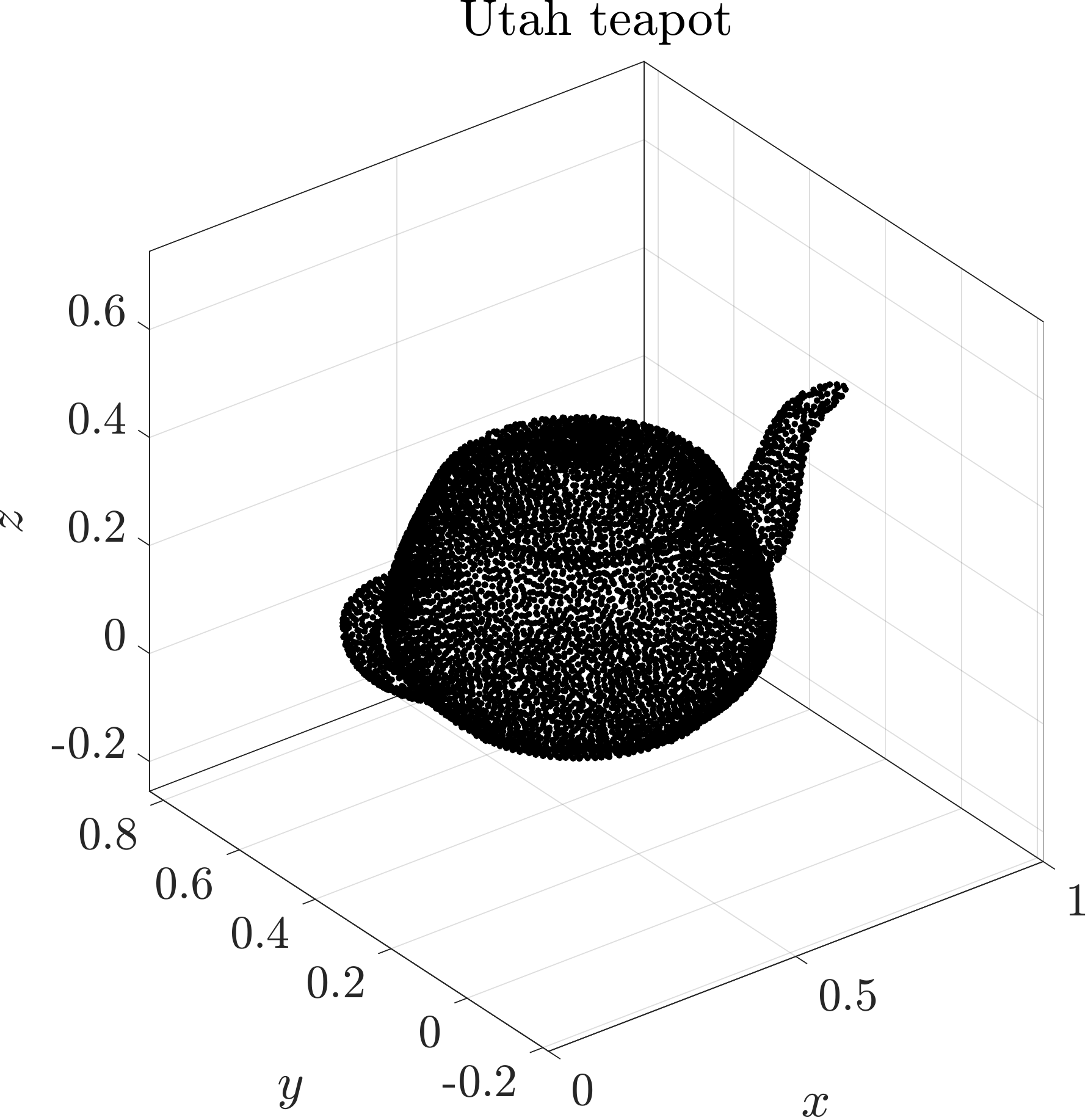}
	\includegraphics[width=0.42\linewidth]{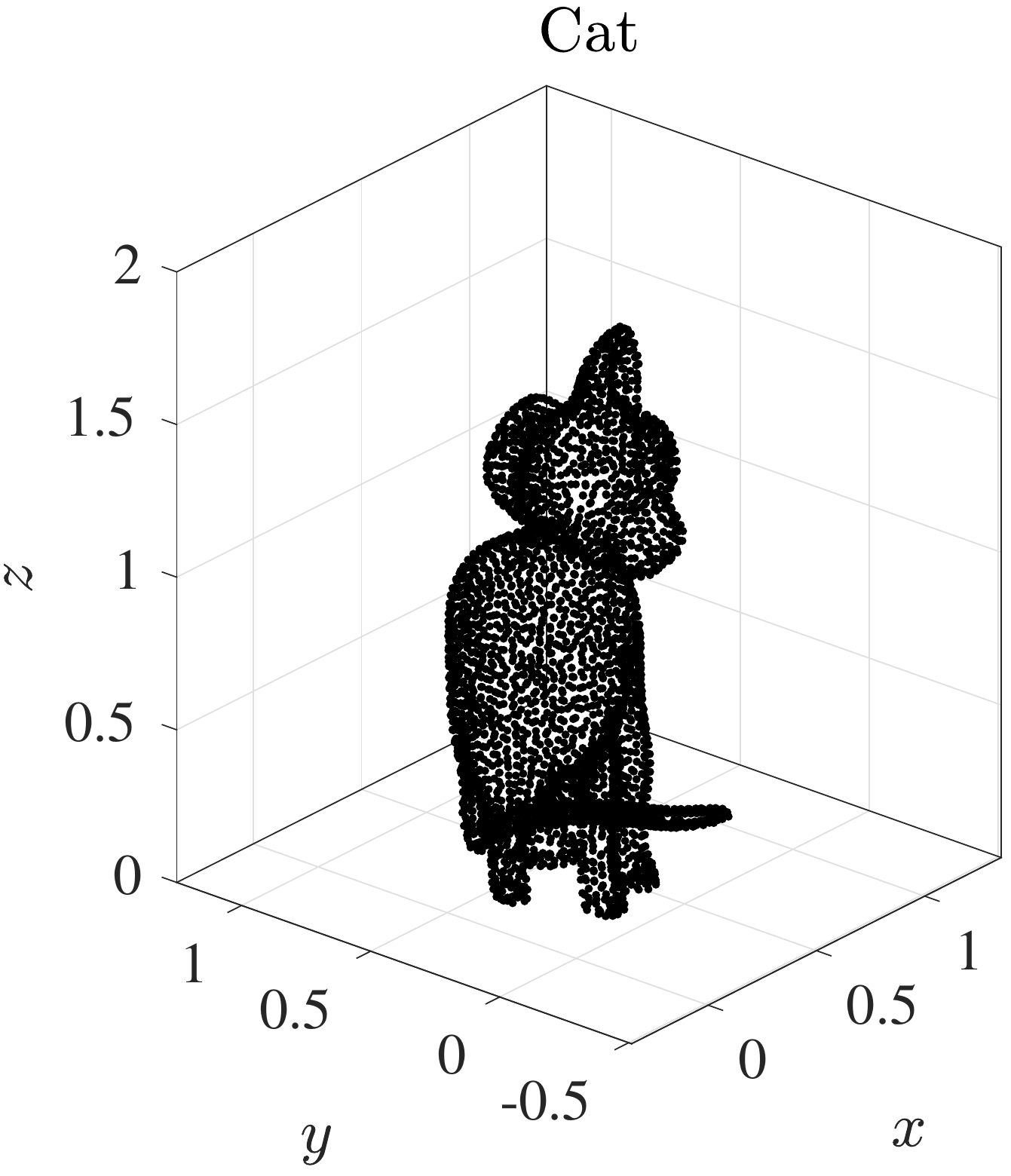}
  \caption{Node sets generated by NURBS-DIVG on the famous CAD Utah Teapot (left) and a CAD model of cat (right) based on~\cite{catmodel}. The Utah Teapot model is made of 32 patches and has 7031 boundary nodes; the cat has 211 patches and 3439 boundary nodes.}
	\label{fig:boundary_example}
\end{figure}

Practical CAD models typically consist of multiple non-overlapping and abutting NURBS surface patches. A NURBS patch meets another NURBS patch at a NURBS curve (the boundary curves of the respective patches). While degenerate situations can easily arise (such as two NURBS patches intersecting at a single point), we \diff{restrict} ourselves in this work with patches that intersect in NURBS curves. Some examples of node sets generated by \diff{the NURBS-DIVG algorithm} on CAD surfaces consisting of multiple NURBS patches are shown in Figure~\ref{fig:boundary_example}. 
\begin{figure}[!tbp]
	\begin{subfigure}[b]{0.5\textwidth}
	\centering
	\includegraphics[width=0.49\linewidth]{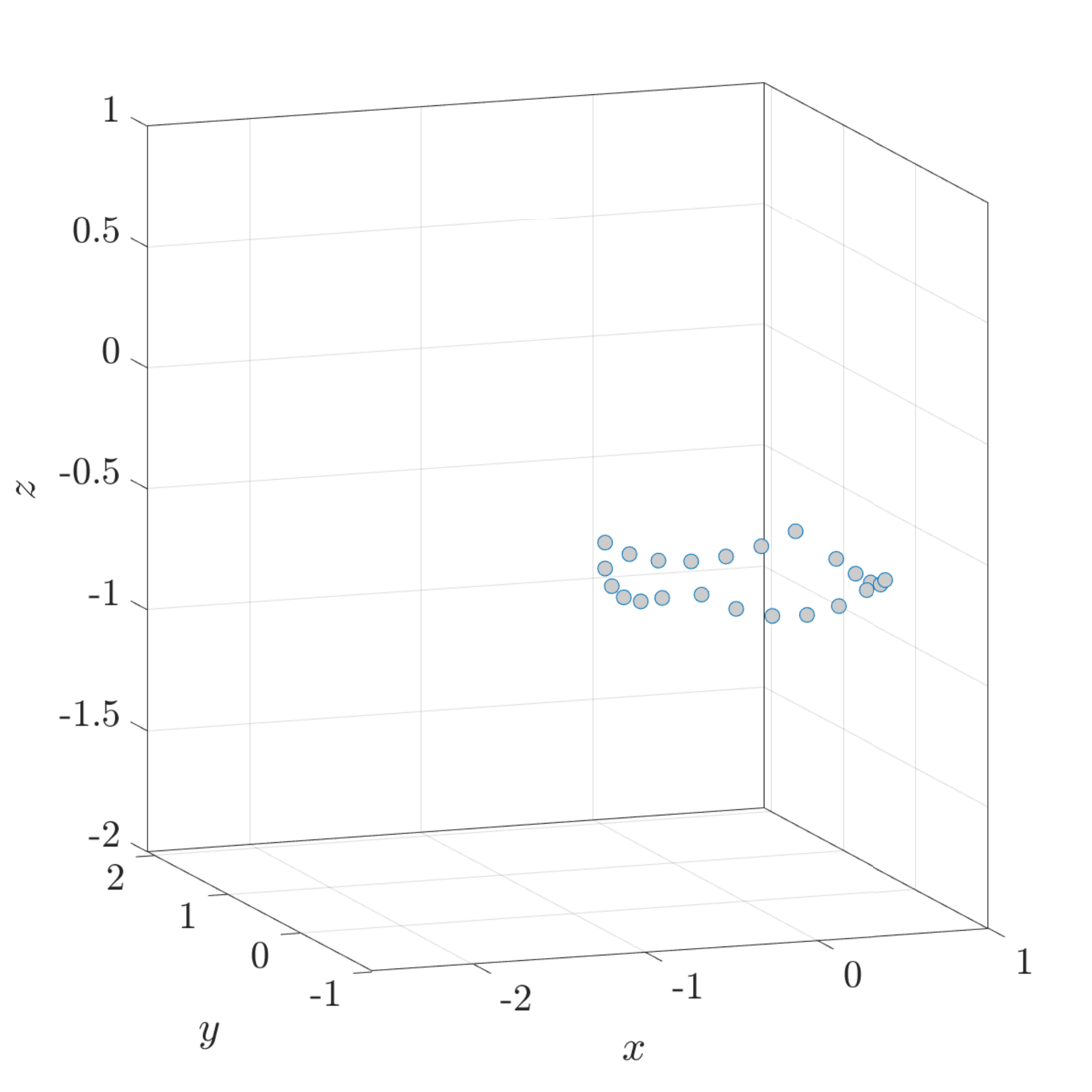}
	\includegraphics[width=0.49\linewidth]{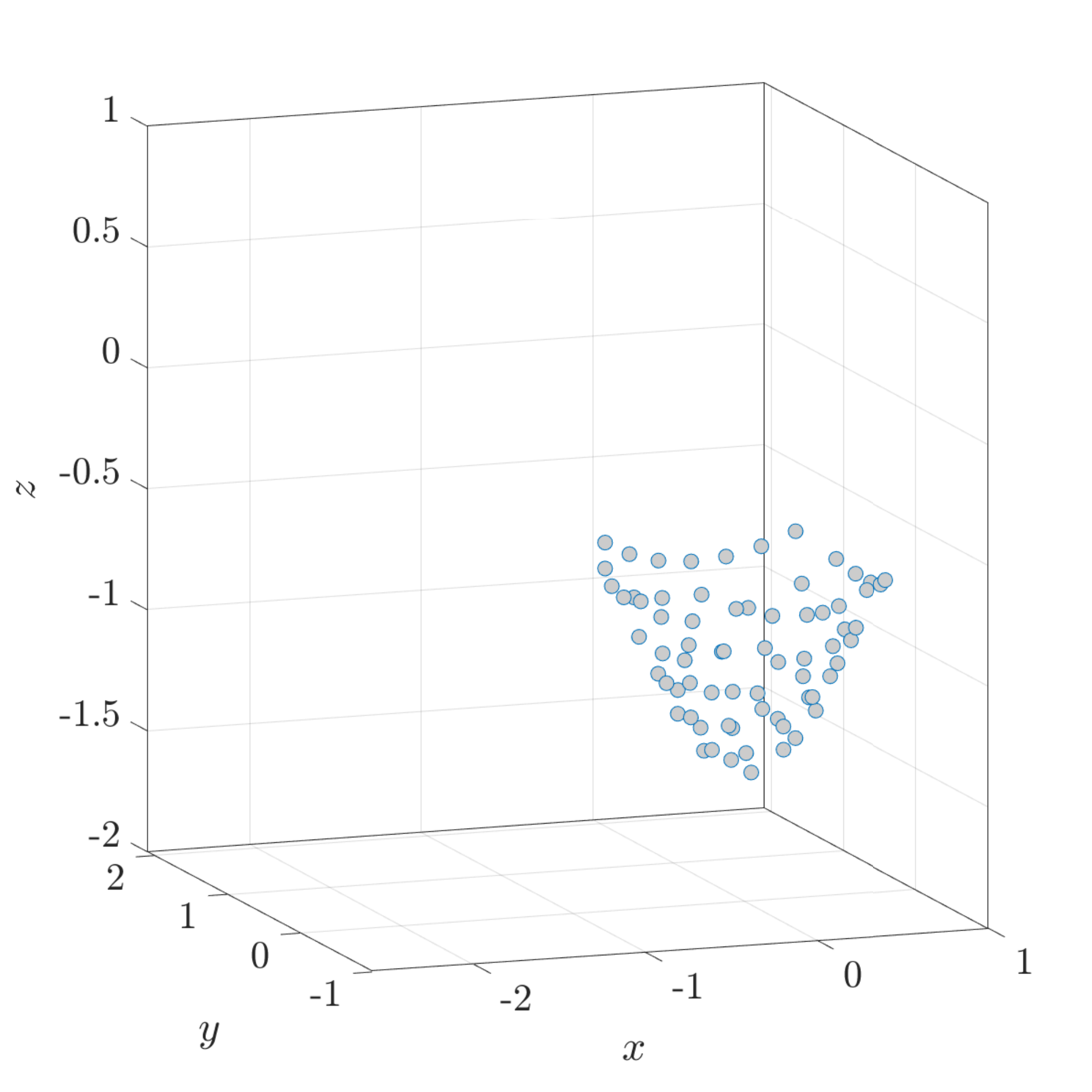}
	\end{subfigure}
	\begin{subfigure}[b]{0.5\textwidth}
	\centering
	\includegraphics[width=0.49\linewidth]{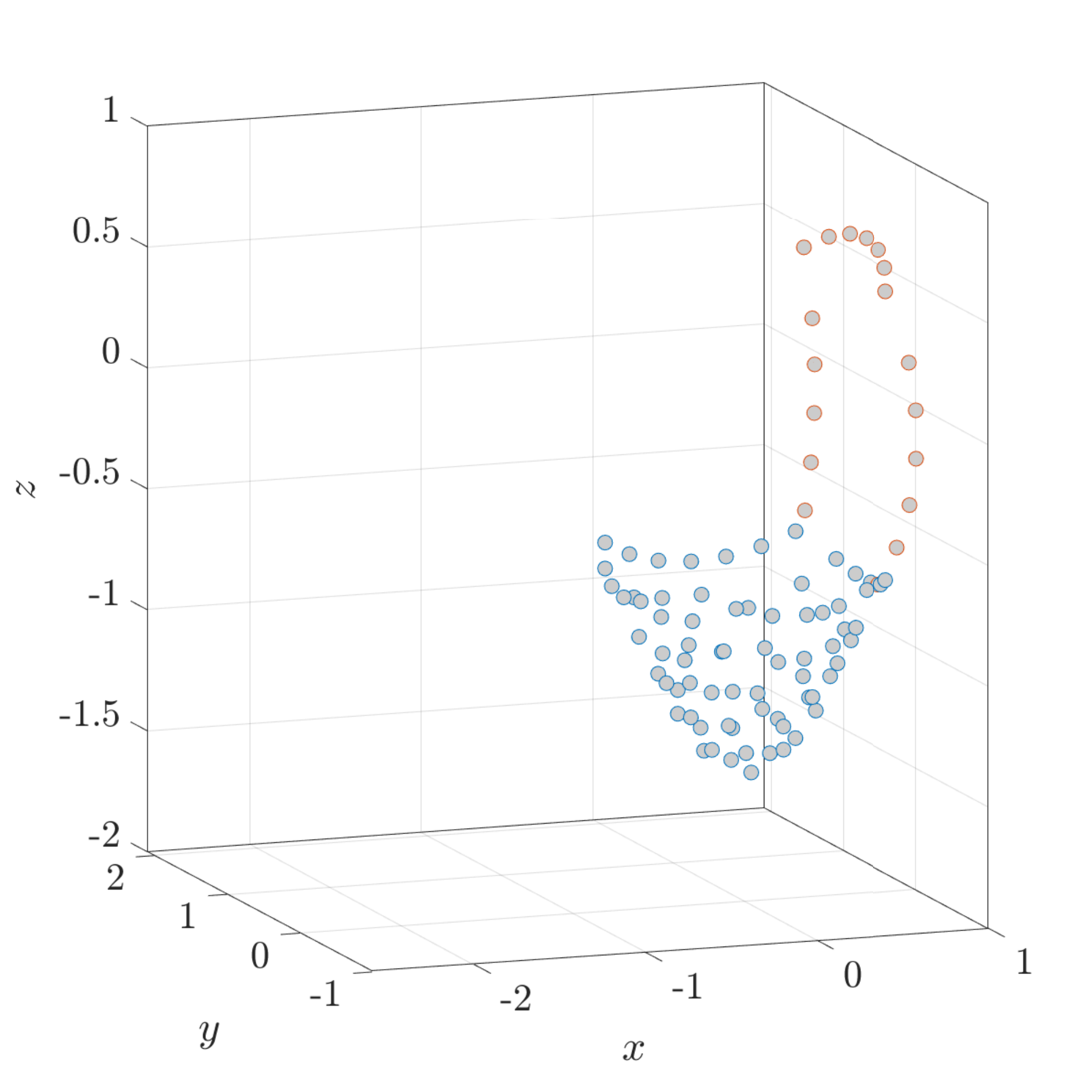}
	\includegraphics[width=0.49\linewidth]{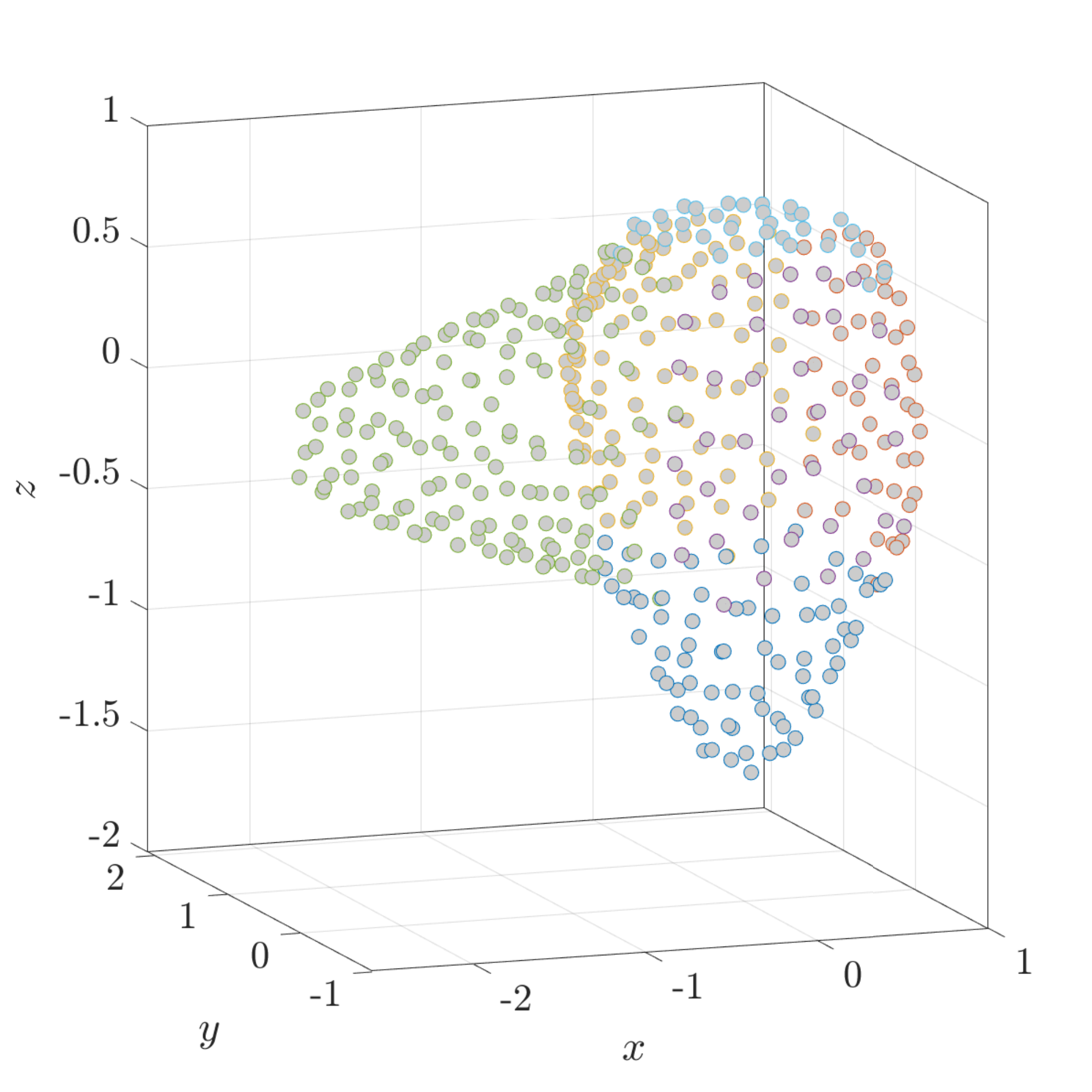}
	\end{subfigure}
	\caption{Illustration of positioning nodes on a deformed sphere made of five NURBS patches. In the first step, the boundary of the first patch is filled (first), followed by filling of that patch interior (second). Once the first patch is processed, the boundary of the second patch is discretized (third); this process is repeated until all patches are fully populated with nodes (fourth).}
\label{fig:patch_exec}
\end{figure}
We now describe the next piece of the NURBS-DIVG algorithm: extending sDIVG to discretize a CAD model that consists of several NURBS patches $\partial \Omega_i$. We proceed as follows: 
\begin{enumerate}
\item We use sDIVG to populate patch boundaries $\partial(\partial \Omega_i)$ with a set of nodes. Recall that these patch boundaries are NURBS curves.
\item We then use these generated nodes as seed nodes within another sDIVG run, this time to fill the NURBS surface patches $\partial \Omega_i$ enclosed by those patch boundaries. 
\end{enumerate}
To populate patch boundaries, the boundary NURBS curve representation obtained from
any of the intersecting patches can be used. But, to use nodes from patch
boundaries as seed nodes in sDIVG for populating surface patches, the corresponding node from the patch's parametric domain $\Lambda$ is required. \diff{However in general, nodes on
  intersecting patch boundaries do not necessarily correspond to the same parametric nodes in all the respective
  parametric domains of intersecting patches. Consequently, the parametric
domains from intersecting patches cannot be joined into one ``global'' parametric domain in a simple and efficient way.} While it is possible to
determine the map $\Omega \to \Lambda$ through a nonlinear solve, we found it \diff{more efficient} to simply populate the patch boundaries twice, once from each of the NURBS
representations obtained from intersecting patches. This produces two sets of seed nodes (one corresponding to each patch), but only the set from one of the representations is used in the final discretization (it does not matter which one, since both node sets are of similar quality). The full process is illustrated in Figure~\ref{fig:patch_exec}.\\

\begin{figure}[!tbph]
	\centering
	\includegraphics[height=0.25\linewidth]{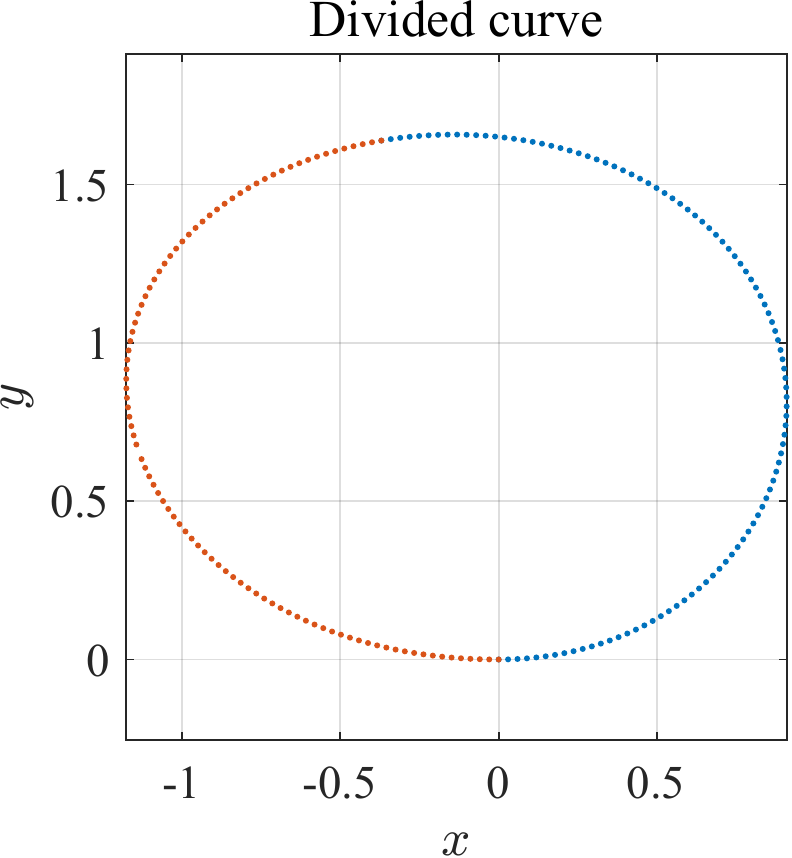}
	\includegraphics[height=0.25\linewidth]{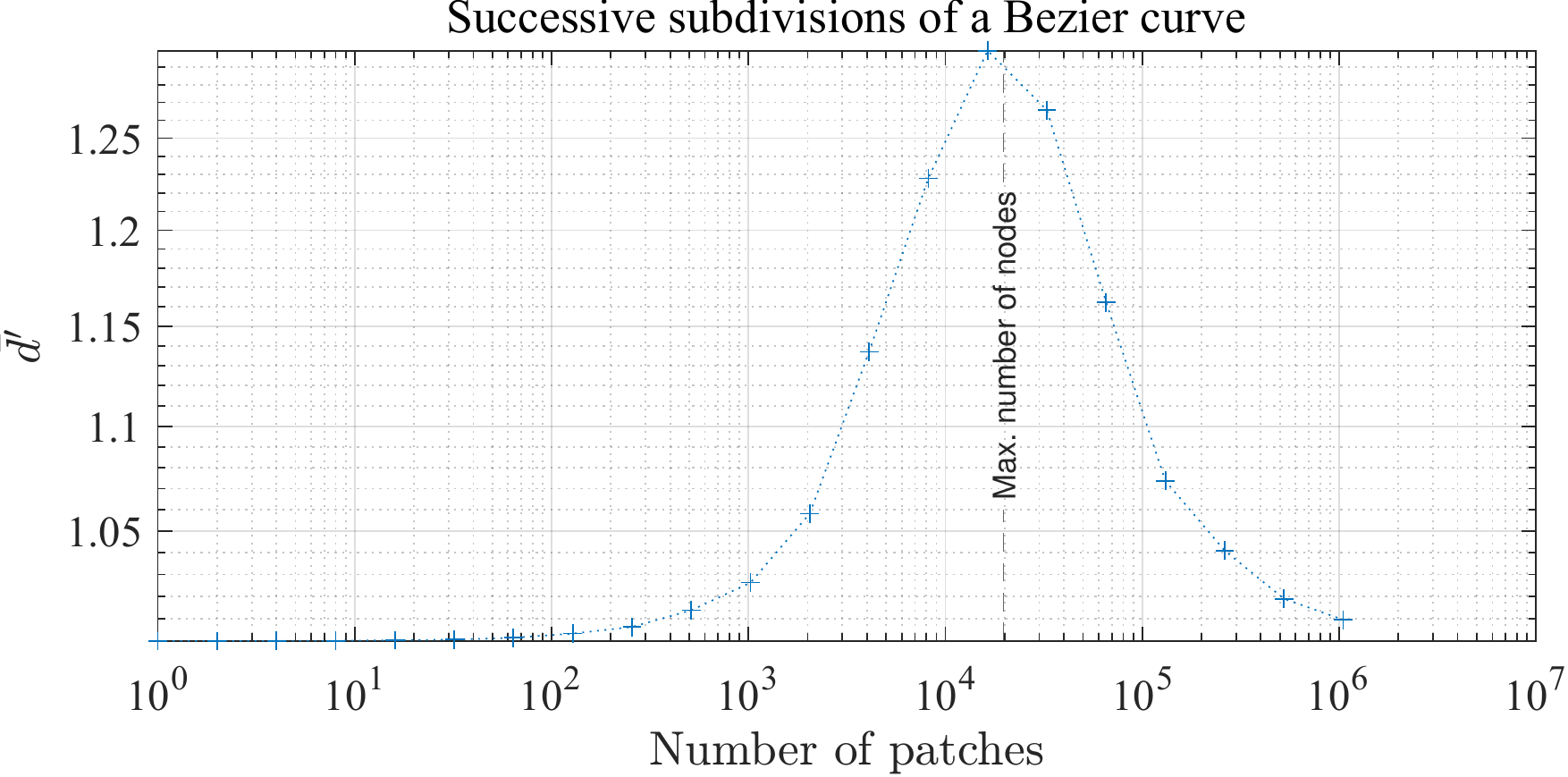}
	
	\includegraphics[height=0.25\linewidth]{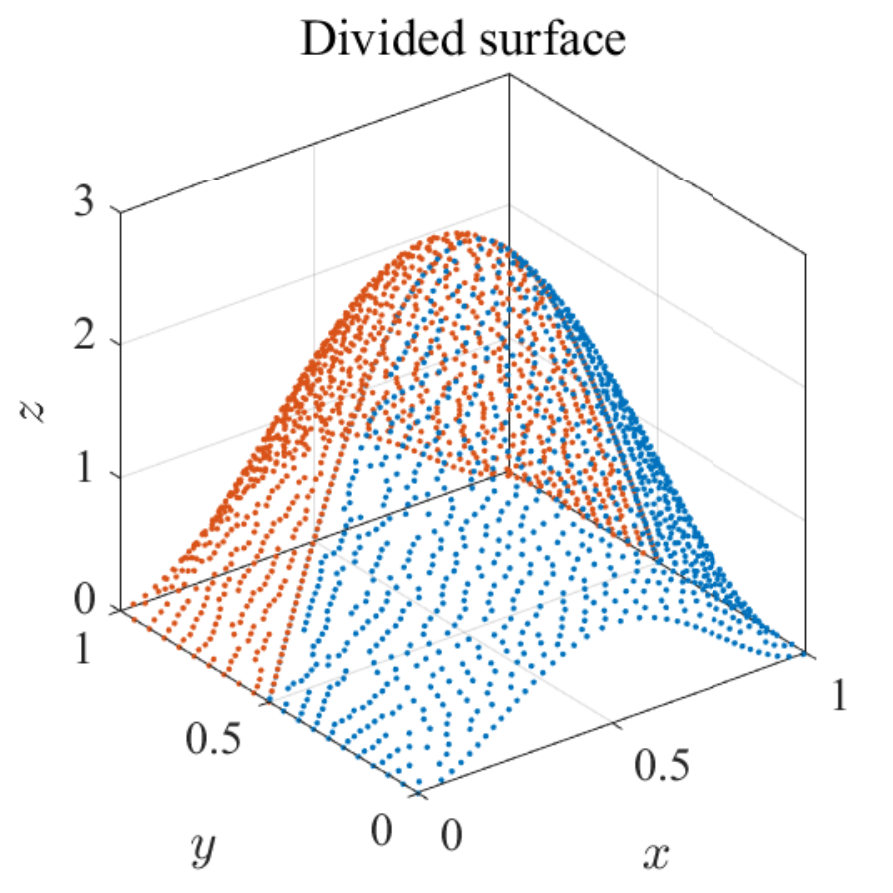}
	\includegraphics[height=0.25\linewidth]{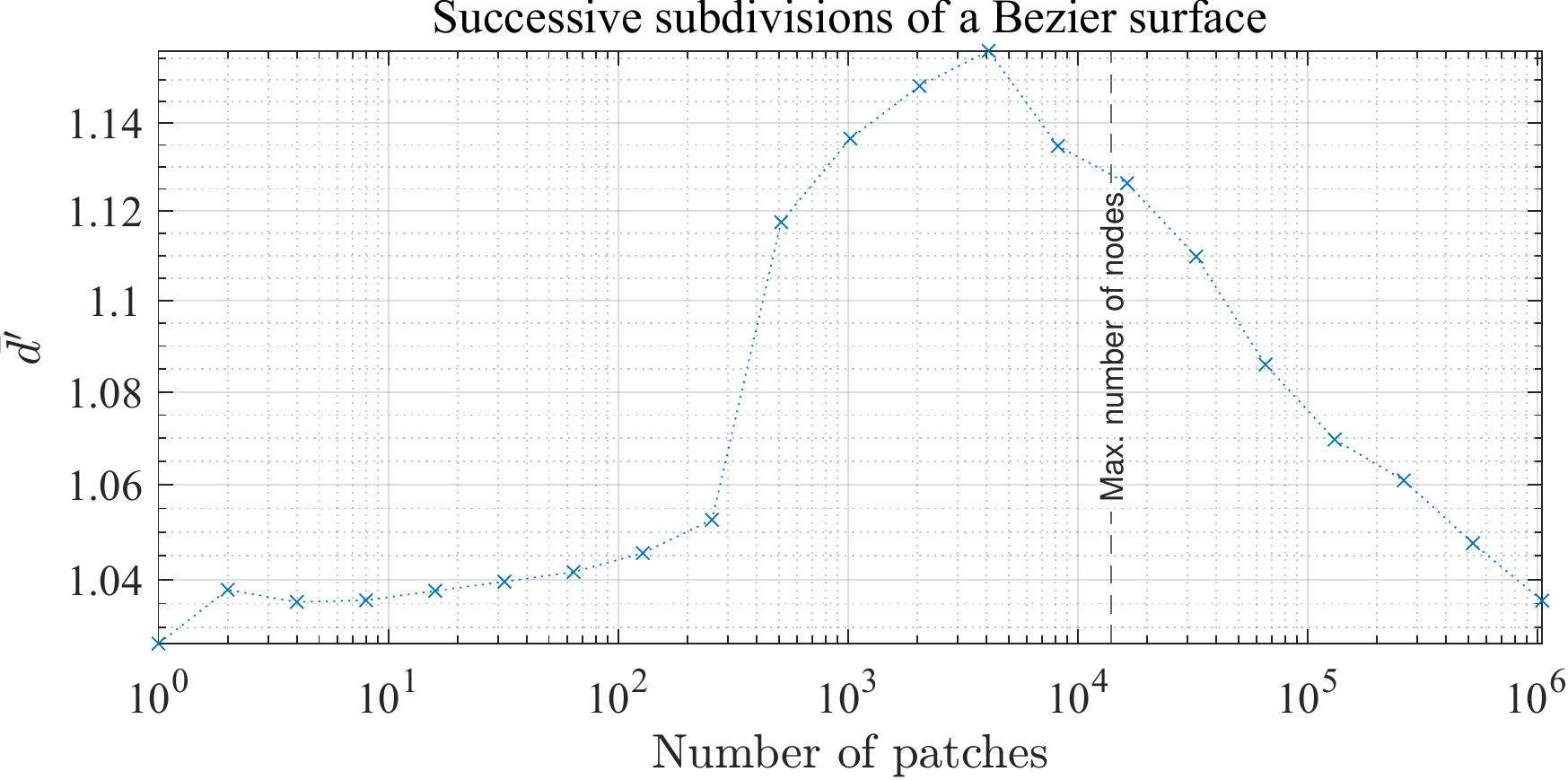}
  \caption{\diff{
		The average normalized distance to $c$ nearest neighbors 
		(see equation Eq.~{\eqref{eq:dist_def}} for the precise definition) averaged over the 
		whole domain, $\overline{d'}$, for a discretization of a successively subdivided Bezier curve in 2D and an analogous Bezier surface in 3D. In 2D $c = 2$ and in 3D $c = 3$ are used. Ideally, one strives for $\overline d' = 1$.
}}
	\label{fig:sub_division}	
\end{figure}

For a given CAD model consisting of a union of NURBS patches and a desired node spacing $h$ (\emph{i.e.} constant spacing function), it is possible that the smallest dimension of the patch becomes comparable to (or even smaller) than $h$. We now analyze the behavior of sDIVG in this regime. To do so, we construct simple models comprising of Bezier surfaces (NURBS with constant weights); to emulate the existence of multiple patches, we simply subdivide the Bezier surfaces to obtain patches. In the 3D  case, each subdivision is performed in a different direction to ensure patches of similar size. The resulting surface, now a union of non-overlapping and abutting NURBS (Bezier) patches, was then discretized with the NURBS-DIVG algorithm using a uniform spacing of $h=10^{-4}$. We then assessed the quality of the resulting node sets on those patches using the normalized local regularity metric $\overline d_i'$ defined in Section~\ref{sec:node_qulity}. The models and this metric are shown in Figure~\ref{fig:sub_division}.  
Figure~\ref{fig:sub_division} shows that NURBS-DIVG works as expected when $h$ is considerably smaller than the patch size. However, once the patch size becomes comparable to the $h$, the NURBS-DIVG algorithm rejects all nodes except those on the boundaries of the patch, as there is not enough space on the patch itself for additional nodes.\footnote{Of course, the figures also show that in the case where $h$ is bigger than the average patch size, NURBS-DIVG works well again, since this scenario is analogous to choosing which patches to place nodes in. However, this scenario is not of practical interest.} In all discussions that follow, we restrict ourselves to the first and most natural regime where $h$ is considerably smaller than the patch size. 

\subsubsection{NURBS-DIVG in the interior of CAD objects}
\label{sec:inout}
As our goal is to generate node sets suitable for meshless numerical analysis \diff{in volumetric domains}, it is vital for the NURBS-DIVG algorithm to be able to generate node sets in the \emph{interior} of volumes whose boundaries are CAD models, in turn defined as a \diff{union of NURBS} patches. While it may appear that the original DIVG algorithm is already well-suited to this task, we encountered a problem of nodes ``escaping'' the domain interior when DIVG was applied naively in the CAD setting.  We now explain this problem and the NURBS-DIVG solution more clearly.

To discretize the interior of CAD objects, we must accurately determine whether
a particular node lies inside or outside the model. The choice of boundary
representation can greatly affect the technique used for such an inside/outside
test. For instance, if the domain boundary is modeled as an implicit surface \diff{(level set)} of
the form $f({\bf x}) = 0$, a node ${\bf x}_k$ is inside if $f({\bf x}_k) < 0$
(up to some tolerance). However, in the case where the domain boundary is
modeled as a parametric surface or a collection of parametric patches (as in
this work), the analogous approach would be to instead solve a nonlinear system of
equations to find the parameter values corresponding to ${\bf x}_k$ and test if
${\bf x}_k$ is inside. A simpler approach used in recent work has been to simply
find the closest point from the boundary discretization ${\bf p}$ (with the
given spacing $h$) to ${\bf x}_k$, and use its unit outward normal to decide if ${\bf x}_k$ is inside the domain. More concretely, if ${\bf n}$ is the unit outward normal vector at ${\bf p}$, ${\bf x}_k$ is inside the domain $\Omega$ when 
\begin{equation}
	{\bf n} \cdot ({\bf x}_k - {\bf p}) < 0.
\end{equation}
\begin{figure}[!htbp]
\centering
	\includegraphics[width=0.40\linewidth]{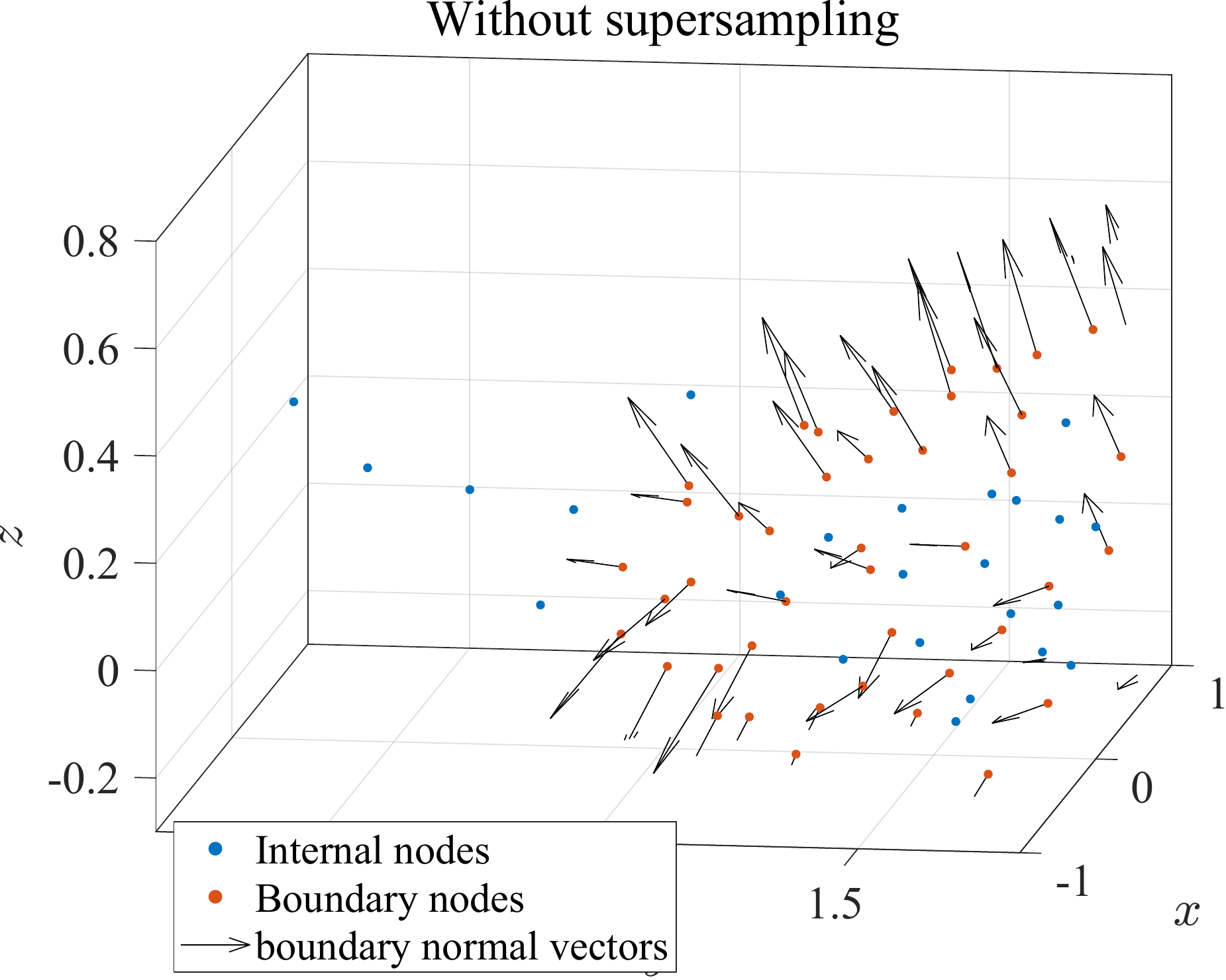}
	\includegraphics[width=0.40\linewidth]{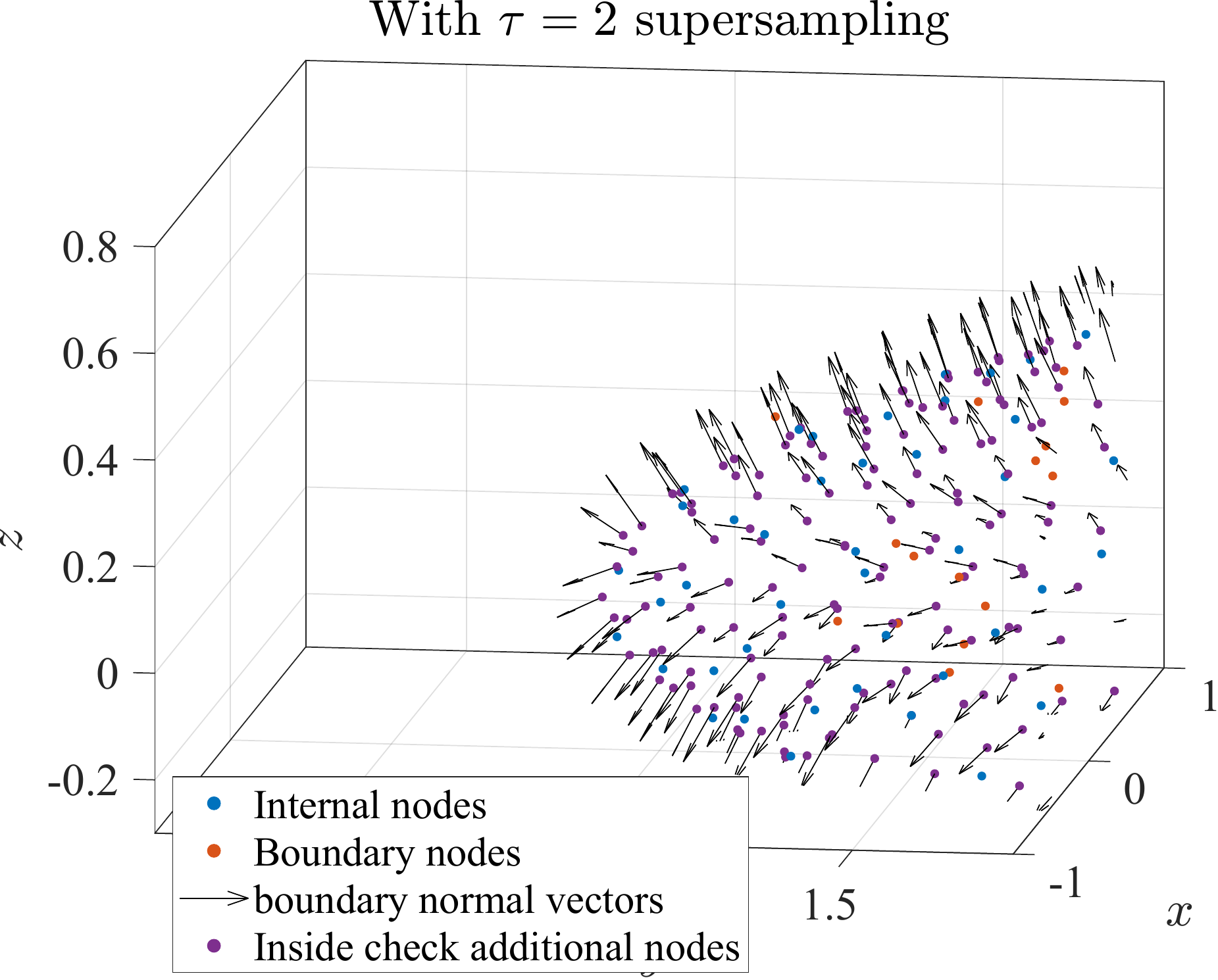}
	\caption{Demonstration of the supersampling approach for the inside/outside test in NURBS-DIVG. The figure on the left shows nodes generated by the naive test used in the DIVG algorithm, with nodes escaping the domain boundary. The figure on the right shows the nodes generated using boundary supersampling in NURBS-DIVG; all non-boundary nodes are enclosed within the volume defined by the boundary.}
	\label{fig:inside_demo}
\end{figure}
This is the approach used by the DIVG algorithm (and many others). However, in our experience, this does not work well for complex geometries with sharp edges and concavities. For an illustration, see Figure \ref{fig:inside_demo} (left); we see nodes marked as ``interior'' nodes that are visually outside the convex hull of the boundary nodes.

An investigation revealed that a relatively coarse sampling of a patch
near its boundary NURBS curves could result in the closest point ${\bf p}$ and
its normal vector $\bf n$ being a bad approximation of the actual closest point and its
normal on the domain $\Omega$, thereby resulting in ${\bf x}_k$ being
erroneously flagged as inside $\Omega$. This problem is especially common on
patch boundaries, where normal vectors do not vary smoothly. NURBS-DIVG uses a simple solution: supersampling. More precisely, we use a secondary set of refined boundary nodes only for the inside check with a reduced spacing $\hat{h}$ given by
\begin{equation}
		\hat h = h / \tau,
\end{equation}
where $\tau > 1$ is a factor that determines the extent of supersampling. Though this potentially requires $\tau$ to be tuned, this solution worked well in our tests with a  minimal additional implementation complexity and computational overhead (see execution profiles for Poisson's equation in section~\ref{sec:PDE}). Figure \ref{fig:inside_demo} (right) shows the effect of setting $\tau = 2$ in the same domain; nodes no longer ``escape'' the boundary. While this approach is particularly useful for boundary represented as a collection of NURBS patches, it is likely to be useful in any setting where the boundary has sharp changes in the derivative of the normal vector (or the node spacing $h$). In fact, intuitively, it seems that the greater the derivative (or the bigger the value of $h$), the greater the value of $\tau$ required to prevent nodes escaping.
\begin{figure}[!htbp]
	\centering
	\includegraphics[width=0.21\linewidth]{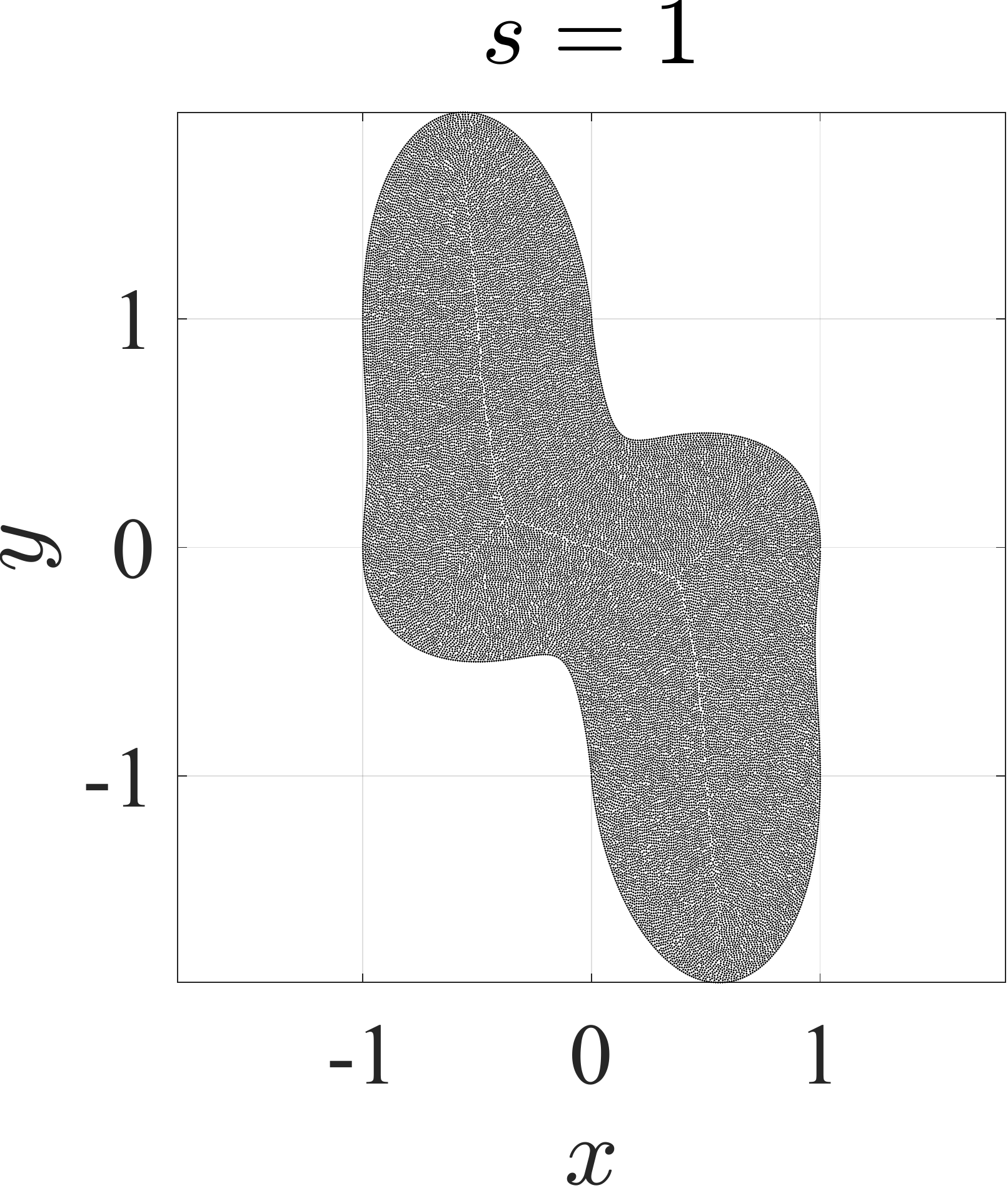}
	\includegraphics[width=0.21\linewidth]{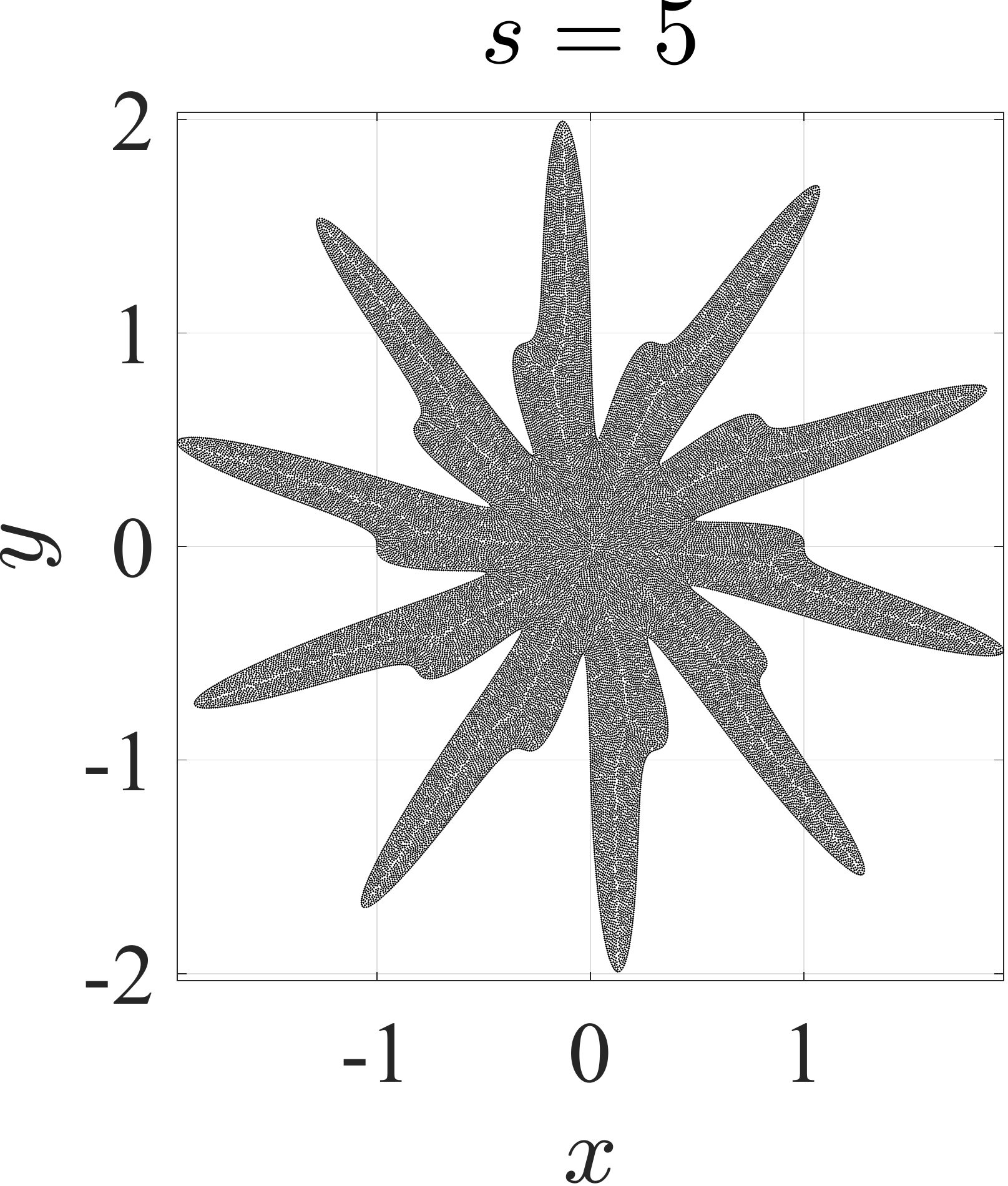}
	\includegraphics[width=0.25\linewidth]{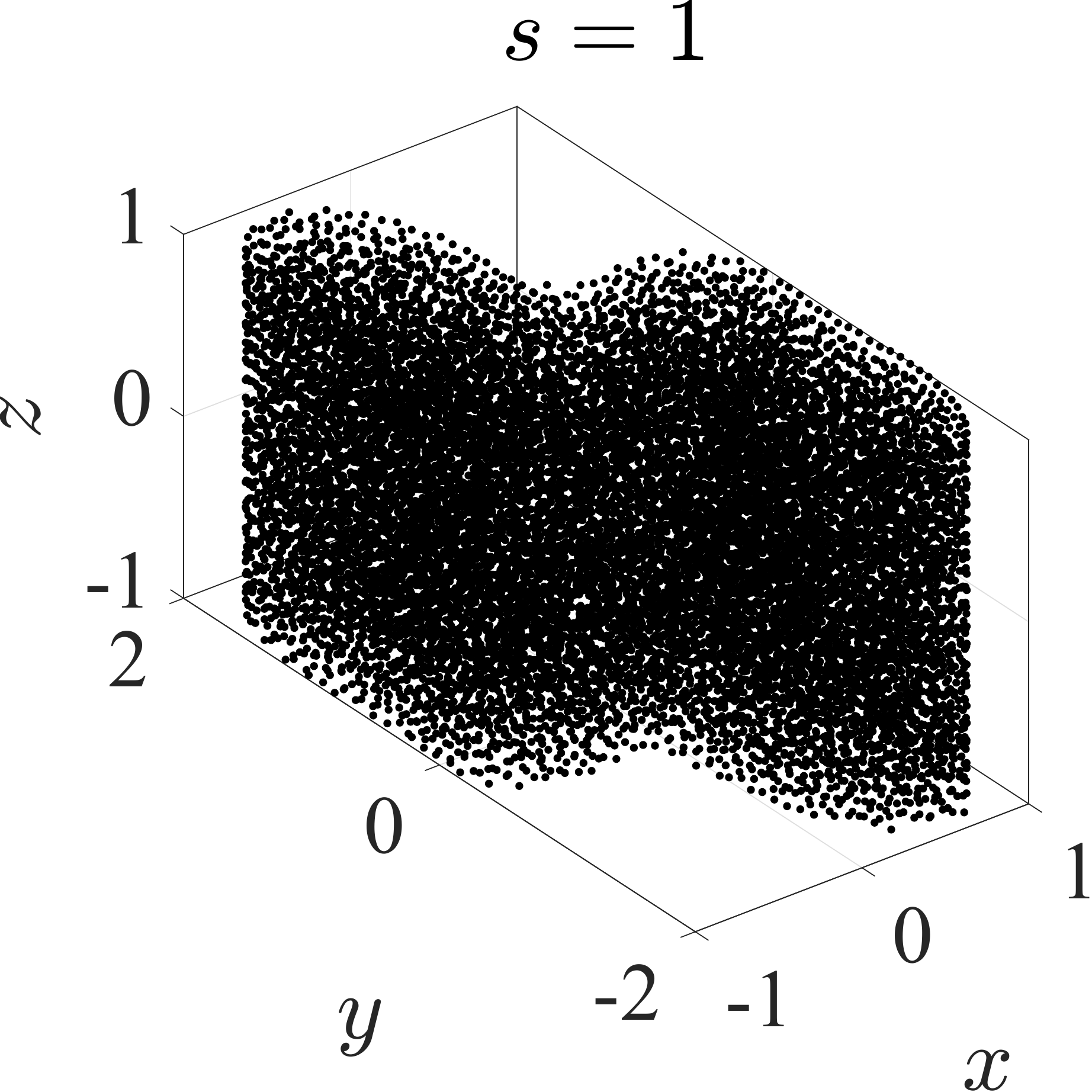}
	\includegraphics[width=0.30\linewidth]{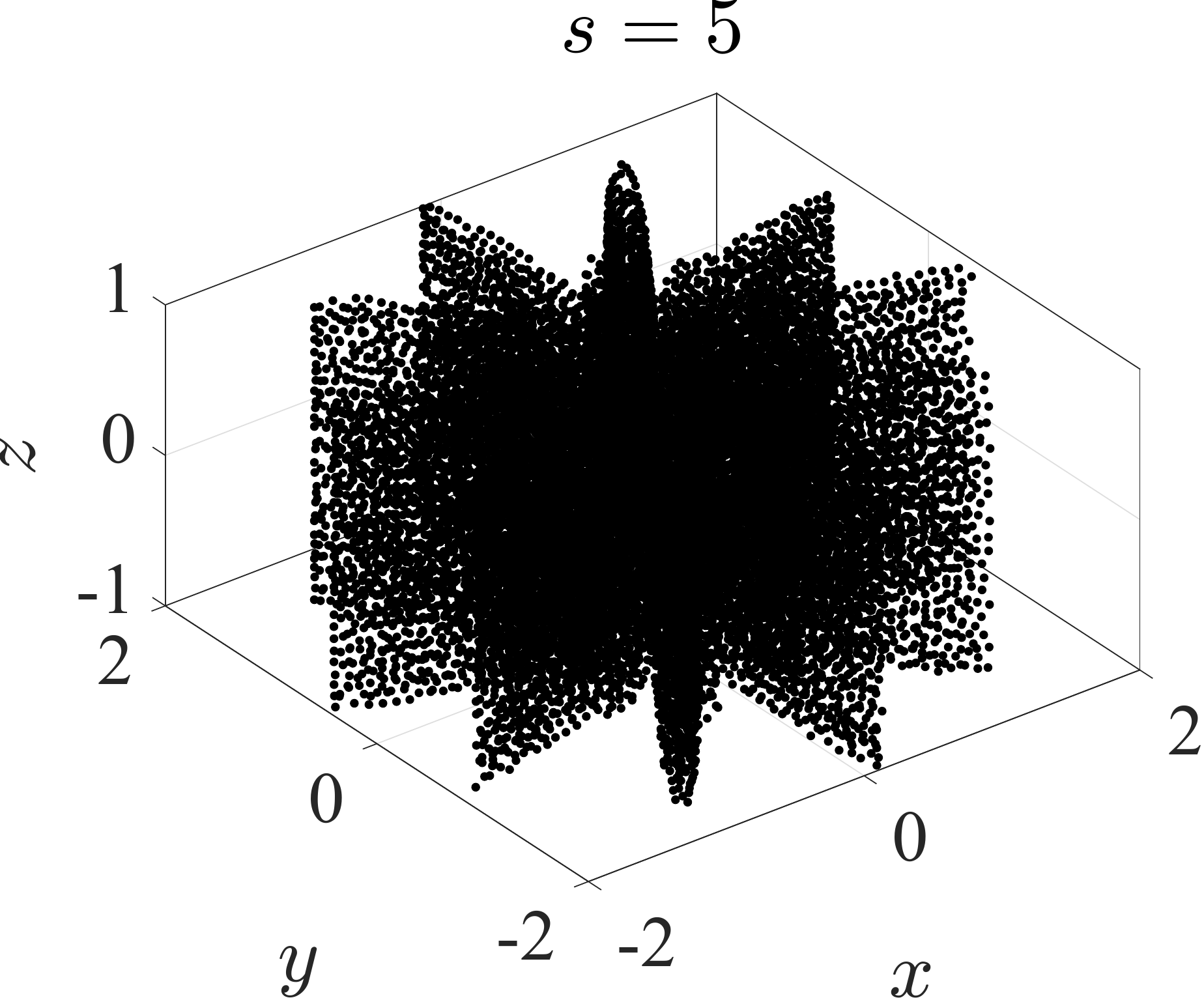}
	\caption{Shapes used to test the supersampling approach in NURBS-DIVG.}
	\label{fig:inside_test_func}
\end{figure}
To confirm this intuition, we run a simple test both in 2D and 3D. In 2D, we define a parametric curve
\begin{equation}
	\label{eg:2dc}
	r(t) = |\cos(st)|^{\sin(2st)}, t \in [0, 2\pi),
\end{equation}
where $s$ is a parameter controlling the complexity of the curve ($s = 1$ gives 2 legs, $s = 1.5$ gives 3 legs, and so forth). These ``legs'' create sharp changes in the derivative of the normal (notice $r(t)$ is not smooth in $t$). In 3D, we simply extrude this curve in the $z$ direction to obtain a surface:
\begin{equation}
	\mathbf r = (r(t) \cos(t), r(t) \sin(t), z), t \in [0, 2\pi), z \in [-1, 1].
\end{equation} 
\begin{figure}[!htbp]
	\centering
	\includegraphics[width=0.49\linewidth]{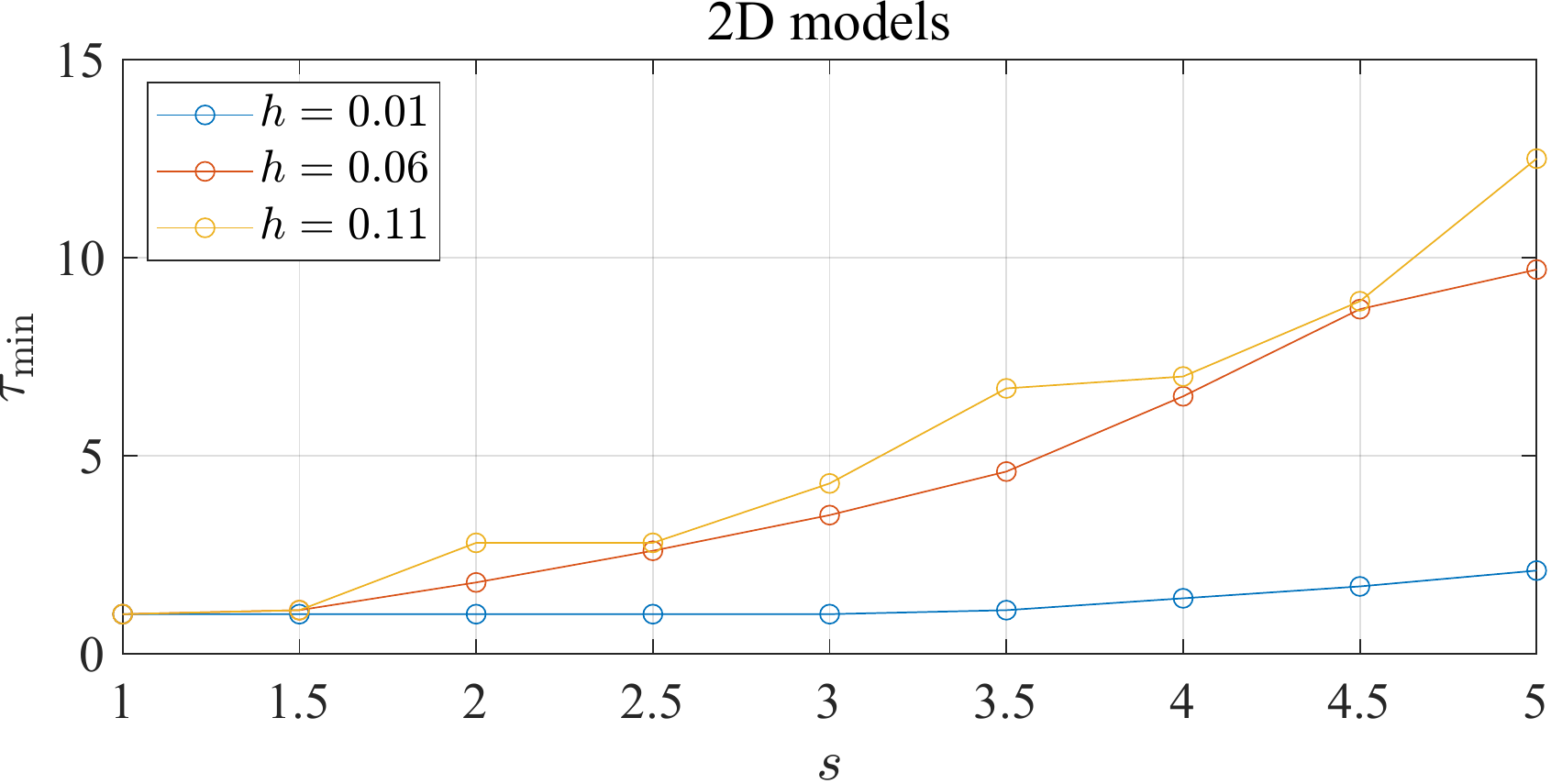}
	\includegraphics[width=0.49\linewidth]{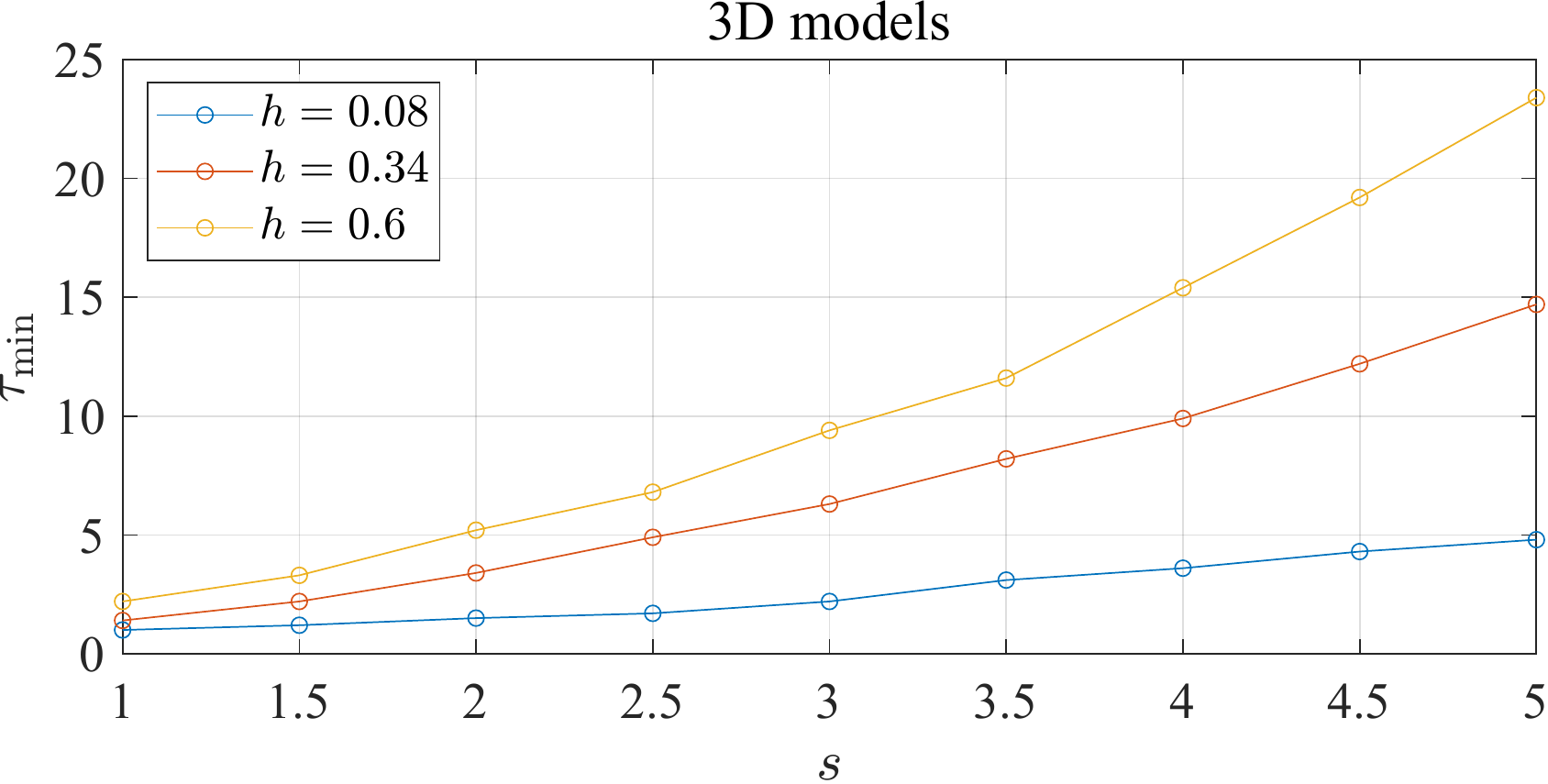}
	\caption{Minimal $\tau$ required to appropriately fill a model.}
	\label{fig:inside_tau}
\end{figure}
Both test domains are depicted in Figure~\ref{fig:inside_test_func}. We then plot the minimum value of $\tau$ required for a successful inside check as a function of $s$, the parameter that controls the number of legs, and the node spacing $h$. The results are shown in Figure \ref{fig:inside_tau}. As expected, increasing the number of legs via $s$ necessitates a greater degree of supersampling ($\tau_{min}$ in the plots) in both 2D and 3D. However, if $h$ is sufficiently small to begin with, smaller values of $\tau$ appear to suffice. In the tests presented in later sections, we selected $h$ to be sufficiently small that $\tau = 2$ sufficed.

\section{Node quality}
\label{sec:node_qulity}
\begin{figure}[h]
	\centering
	\includegraphics[width=0.30\linewidth]{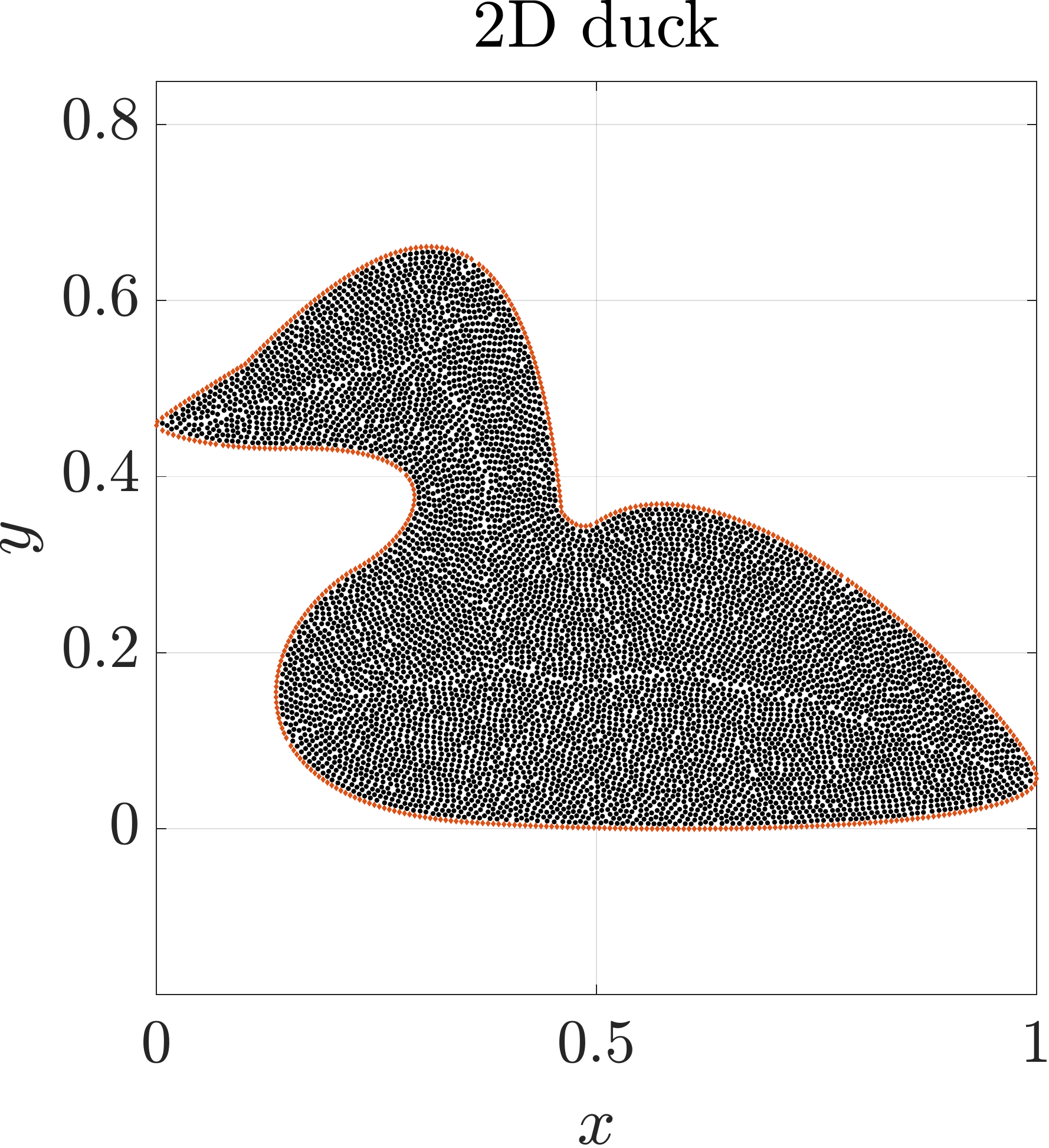}
	\includegraphics[width=0.31\linewidth]{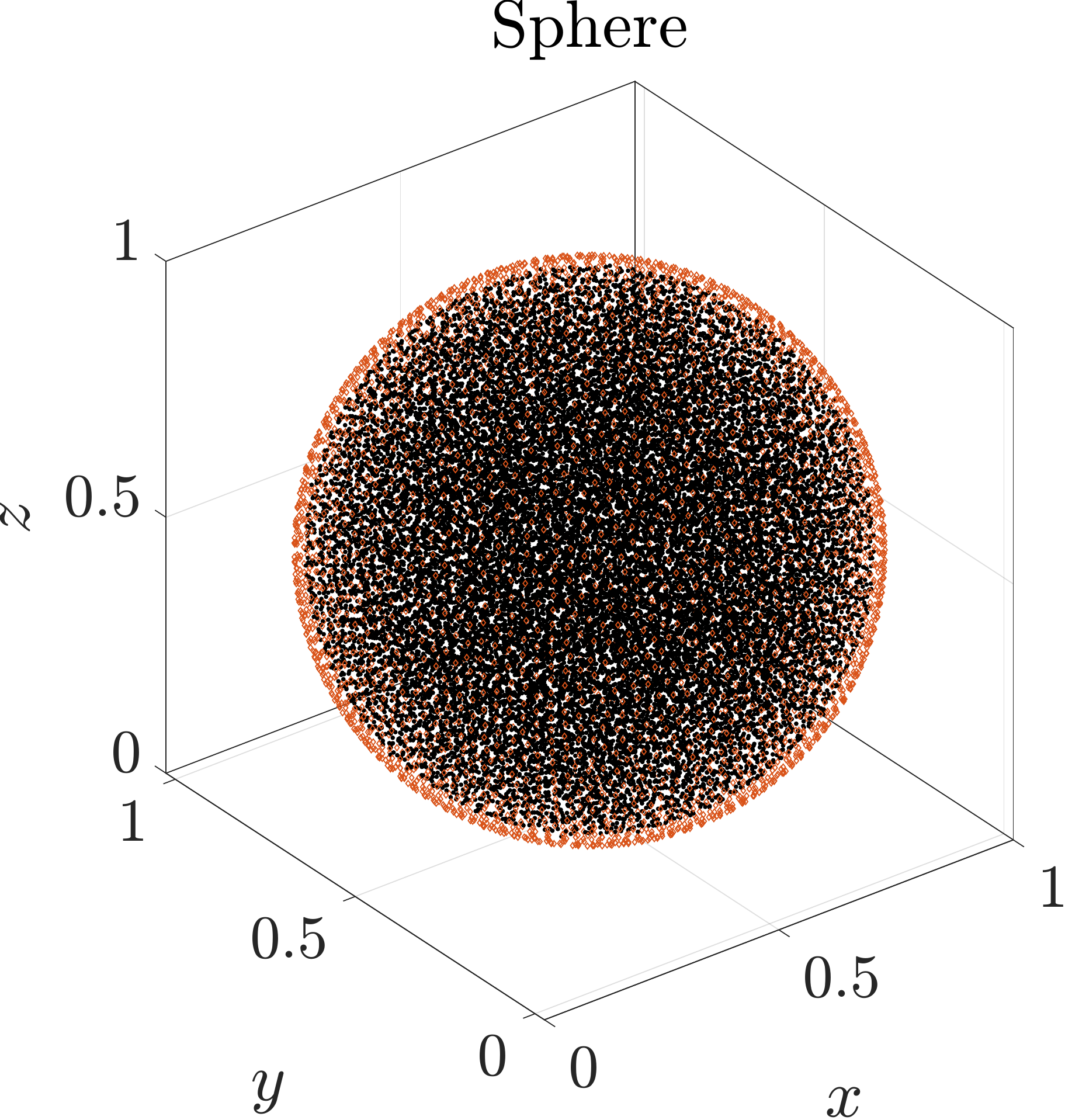}
	\includegraphics[width=0.31\linewidth]{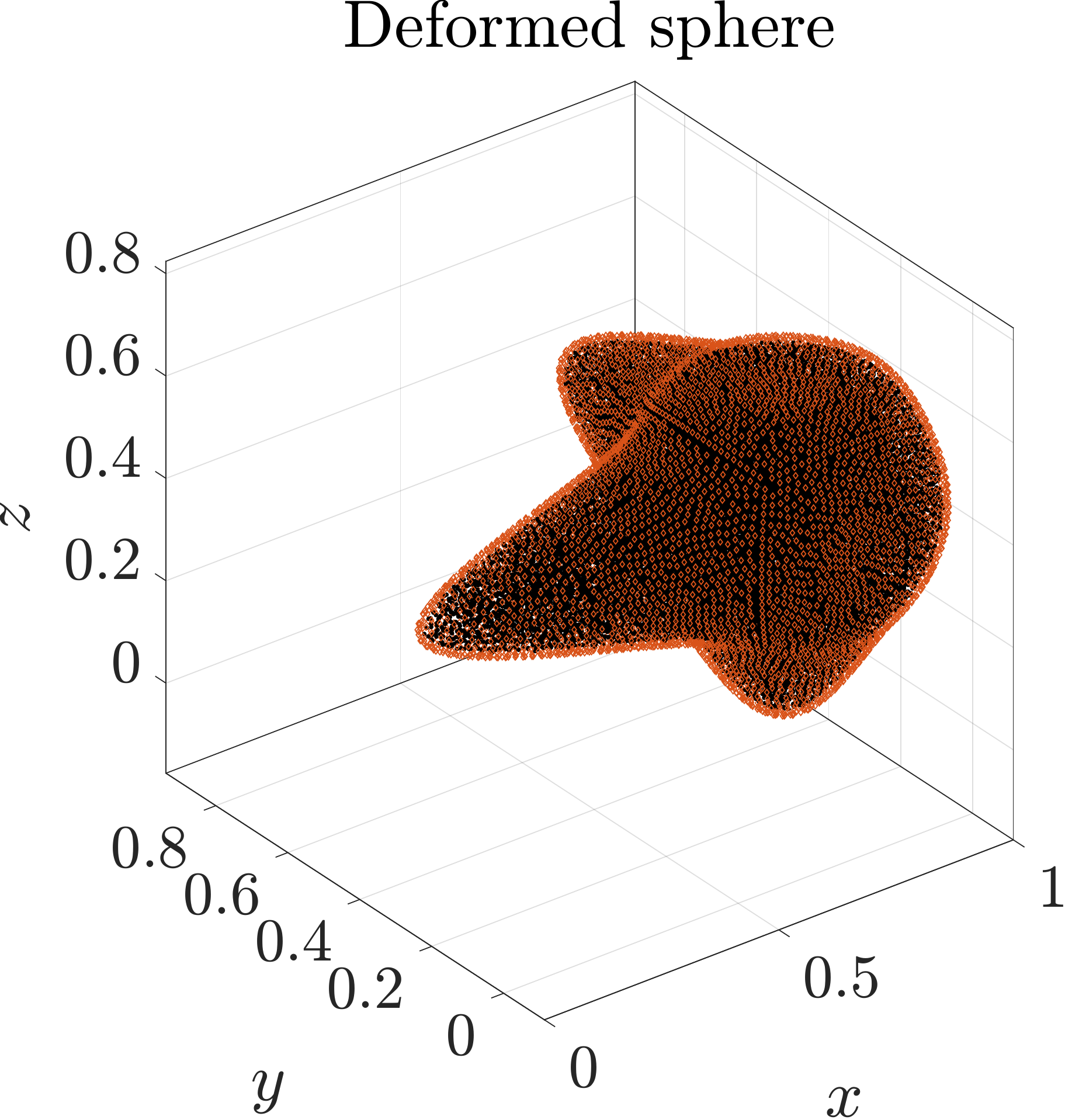}
	\caption{Geometries used in the node quality analysis.}
	\label{fig:node_q_e}
\end{figure}

Although node quality in \diff{the} meshfree context is not as well understood as in mesh
based methods, we can analyze local regularity by examining distance distributions to
nearest neighbors. For each node $\mathbf p_i$ with nearest neighbors $\mathbf p_{i, j}, j = 1, \dots c$ we compute
\begin{align}
  \overline{d_i}    & = \frac{1}{c} \sum_{j = 1}^c ||\mathbf p_i - \mathbf p_{i, j}||, \label{eq:dist_def} \\
	d_i^\mathrm{min} & =  \min_{j=1, \dots, c} ||\mathbf p_i - \mathbf p_{i, j}||,                          \\
	d_i^\mathrm{max} & =  \max_{j=1, \dots, c} ||\mathbf p_i - \mathbf p_{i, j}||.
\end{align}
\diff{In our analysis, the spacing function $h$ is constant over the whole domain, therefore the quantities are normalized as}
\begin{equation}
	\overline{d_i'} = \overline{d_i}/h.
\end{equation}
For the analysis the following models are selected
\begin{itemize}
	\item 2D duck with 8 patches,
	\item sphere with 6 patches, 
	\item deformed sphere with 6 patches,
\end{itemize}
all depicted in Figure~\ref{fig:node_q_e}.

\begin{figure}[h]
	\centering
	\includegraphics[width=0.8\linewidth]{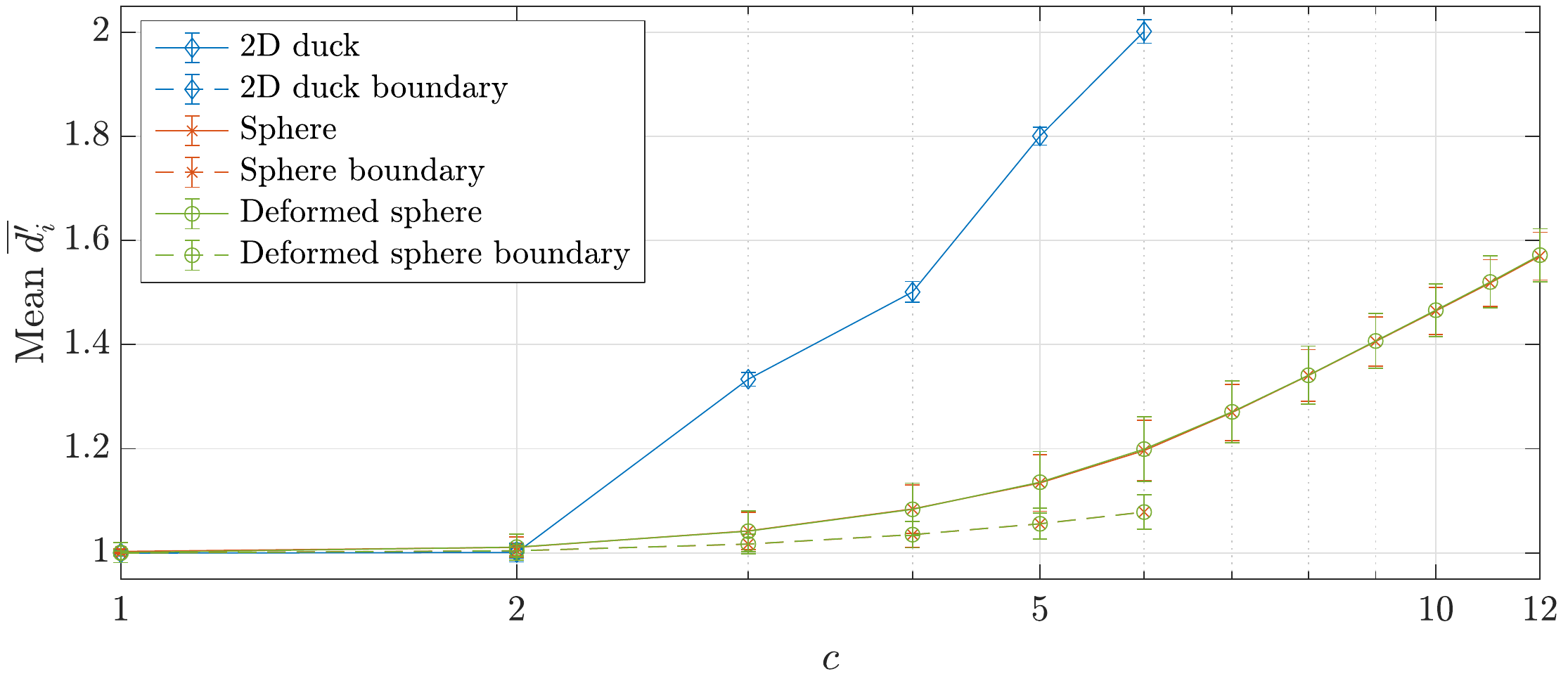}
	\caption{\diff{Mean and standard deviation (depicted as error bars) of local regularity distributions $\overline{d_i'}$
  computed over the whole domain as a function of the number of nearest neighbors $c$ for all three test models.}}
	\label{fig:loc_reg_of_c}
\end{figure}

\begin{figure}[h]
	\centering
	\includegraphics[width=0.32\linewidth]{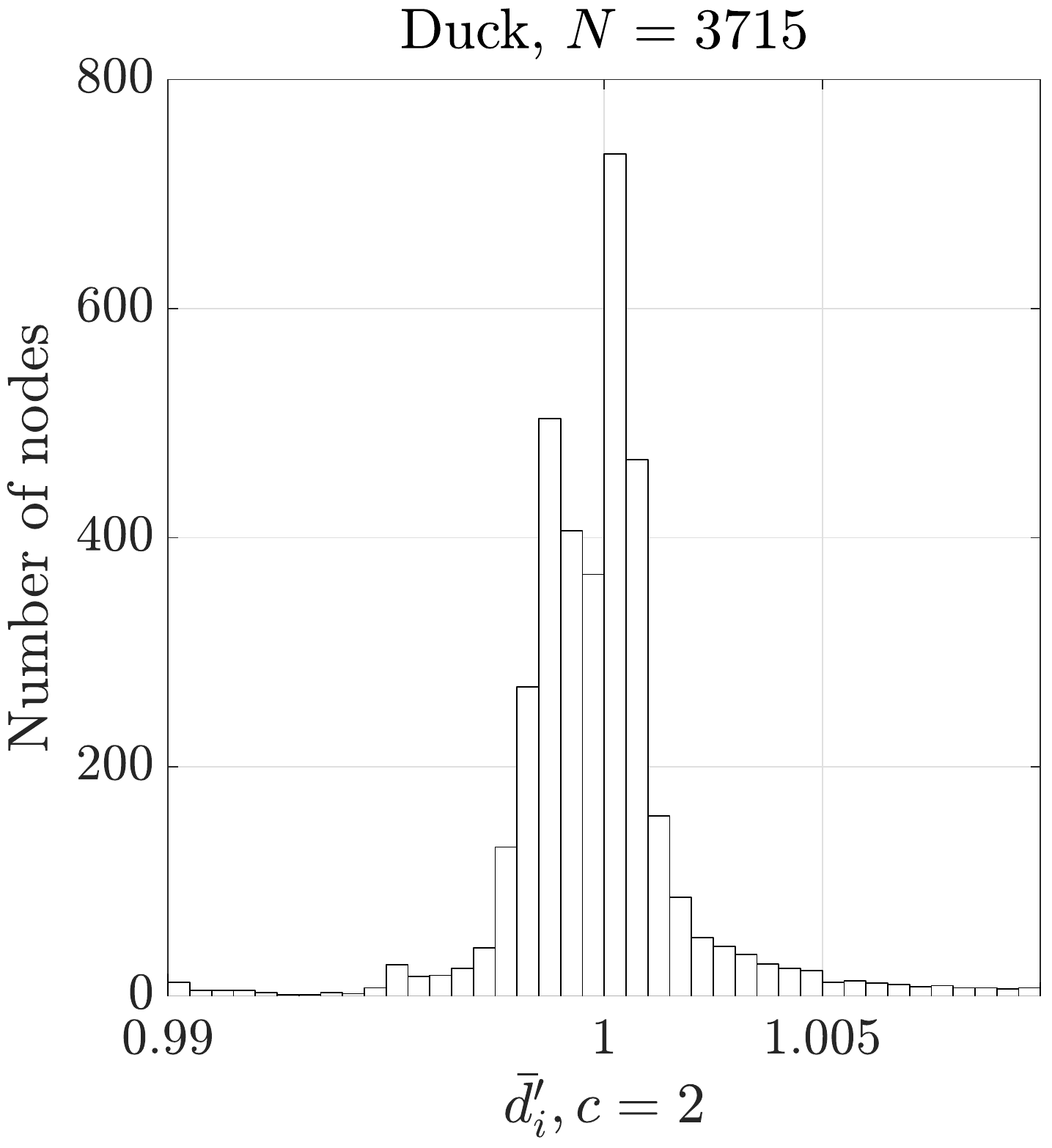}
	\includegraphics[width=0.32\linewidth]{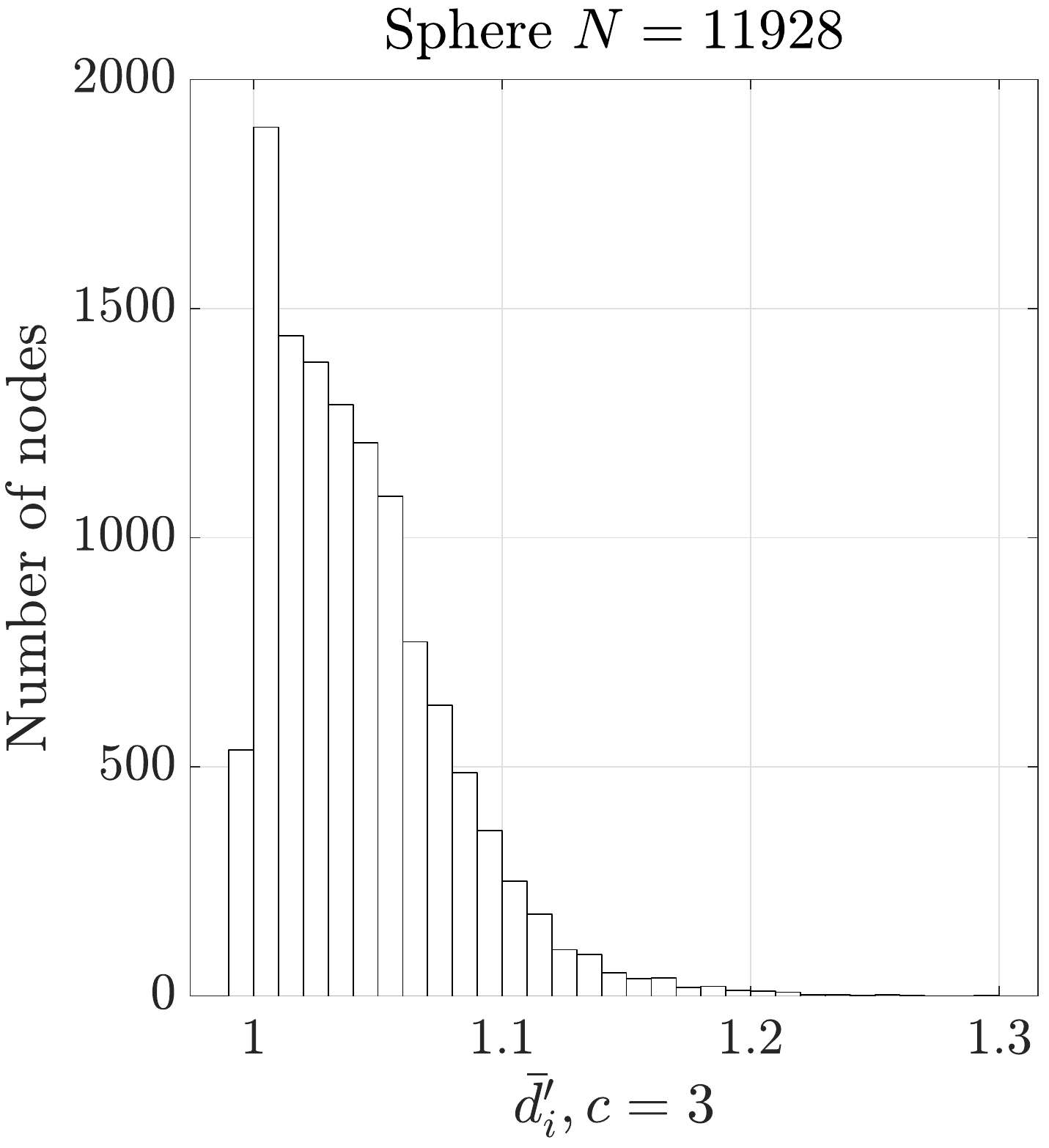}
	\includegraphics[width=0.32\linewidth]{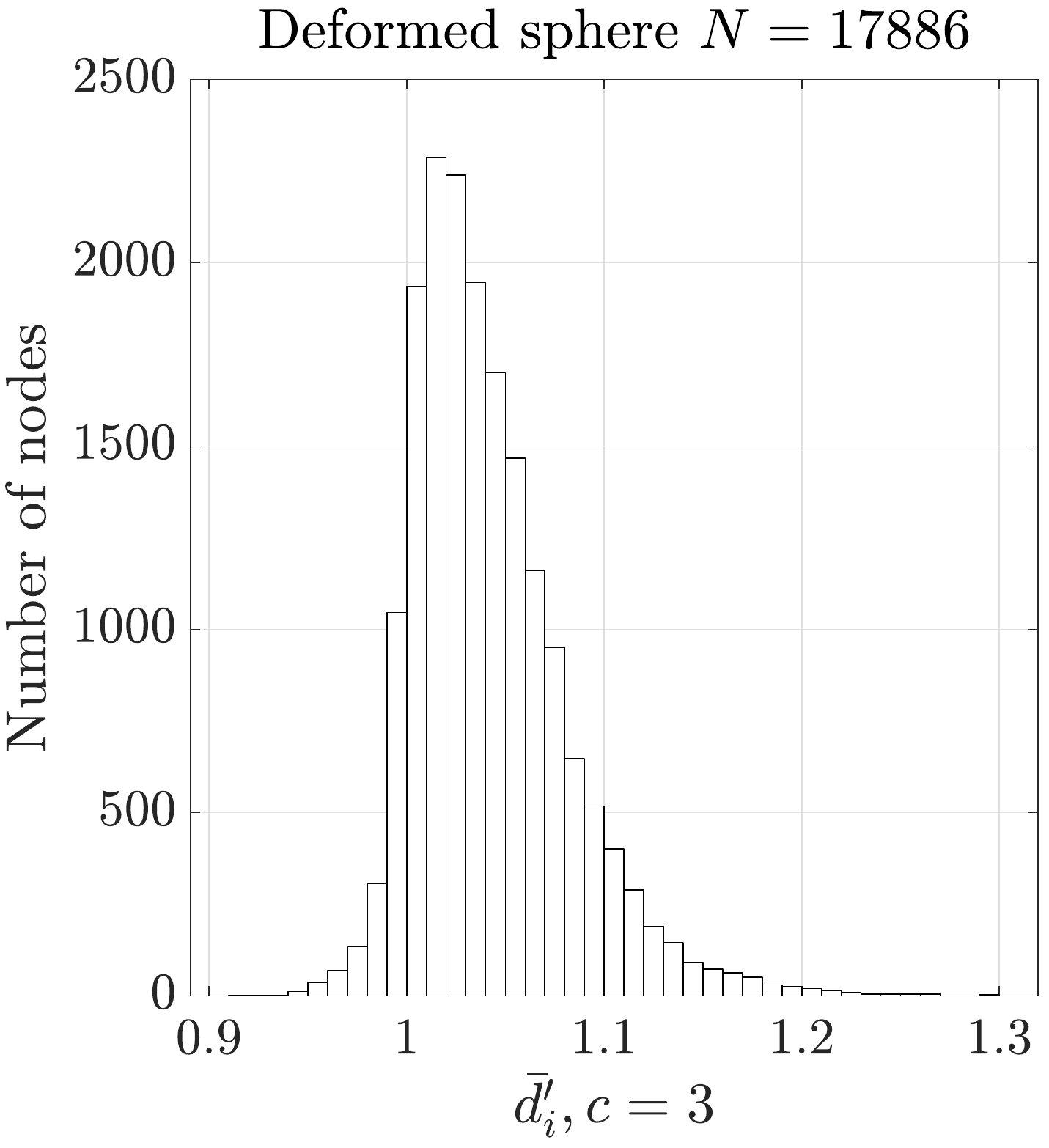}
	\caption{Local regularity distributions for boundary nodes \diff{in the case of a constant spacing function}.}
	\label{fig:loc_reg_1}
\end{figure}

\begin{figure}[h]
	\centering
	\includegraphics[width=0.32\linewidth]{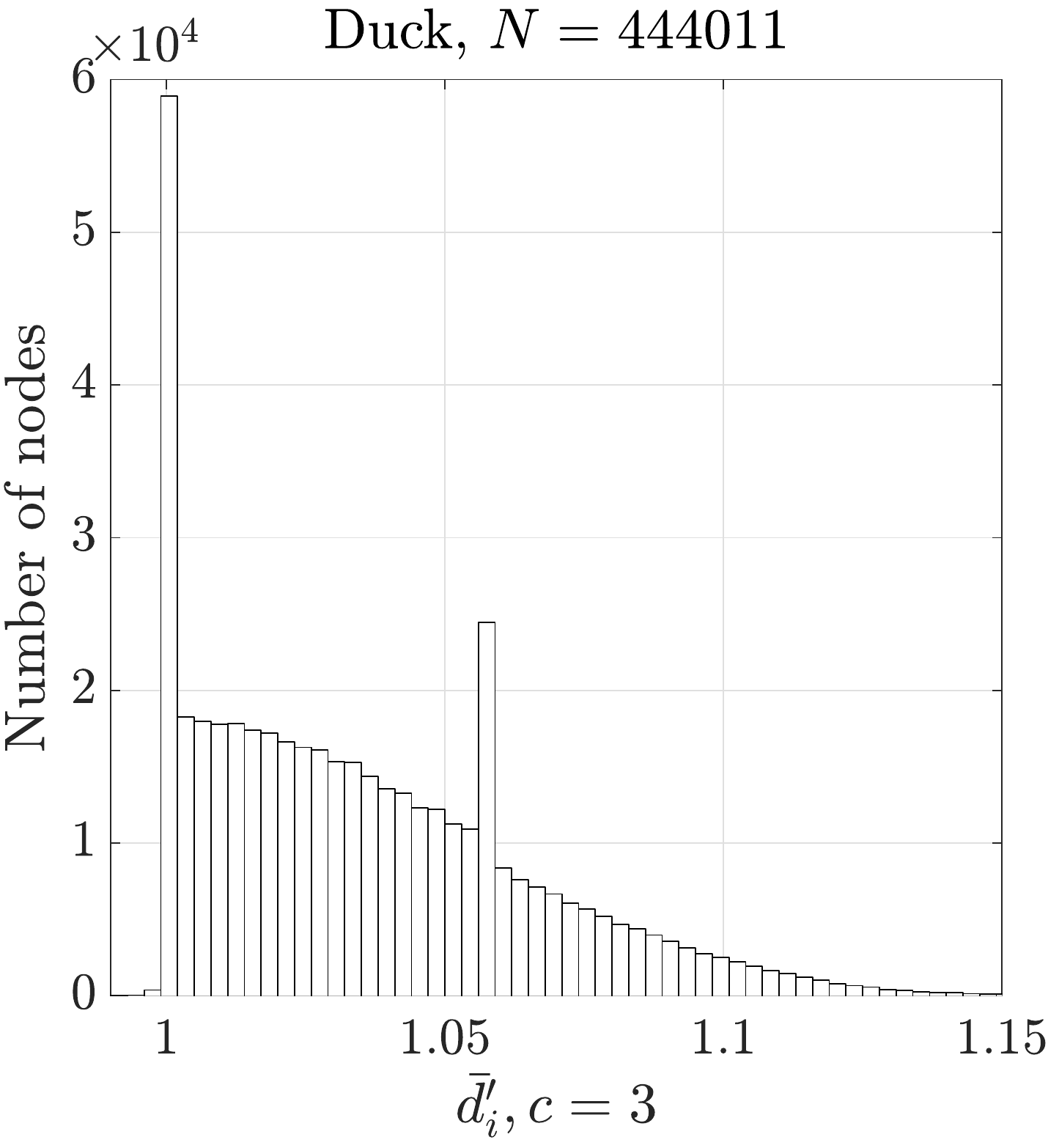}
	\includegraphics[width=0.32\linewidth]{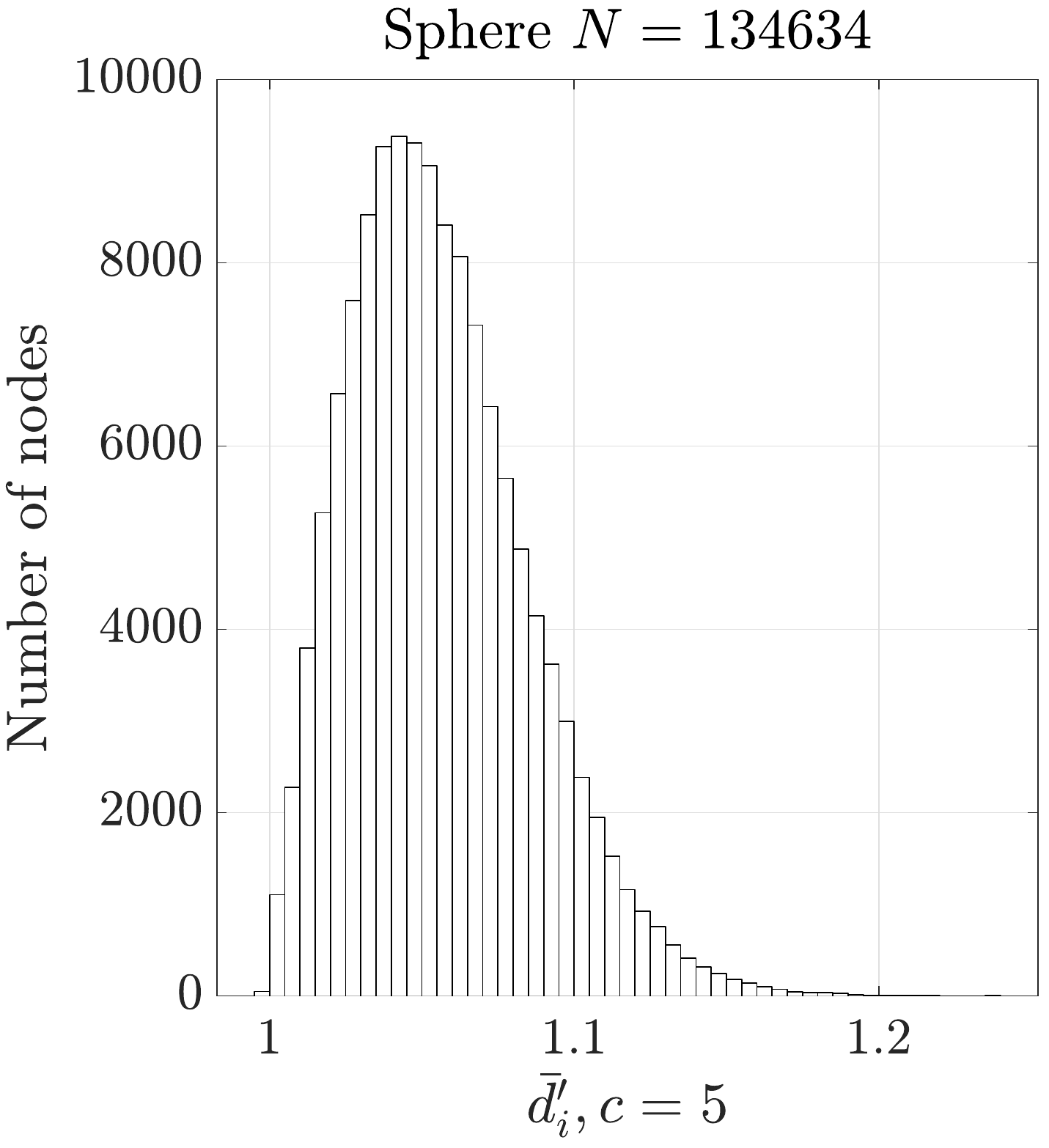}
	\includegraphics[width=0.32\linewidth]{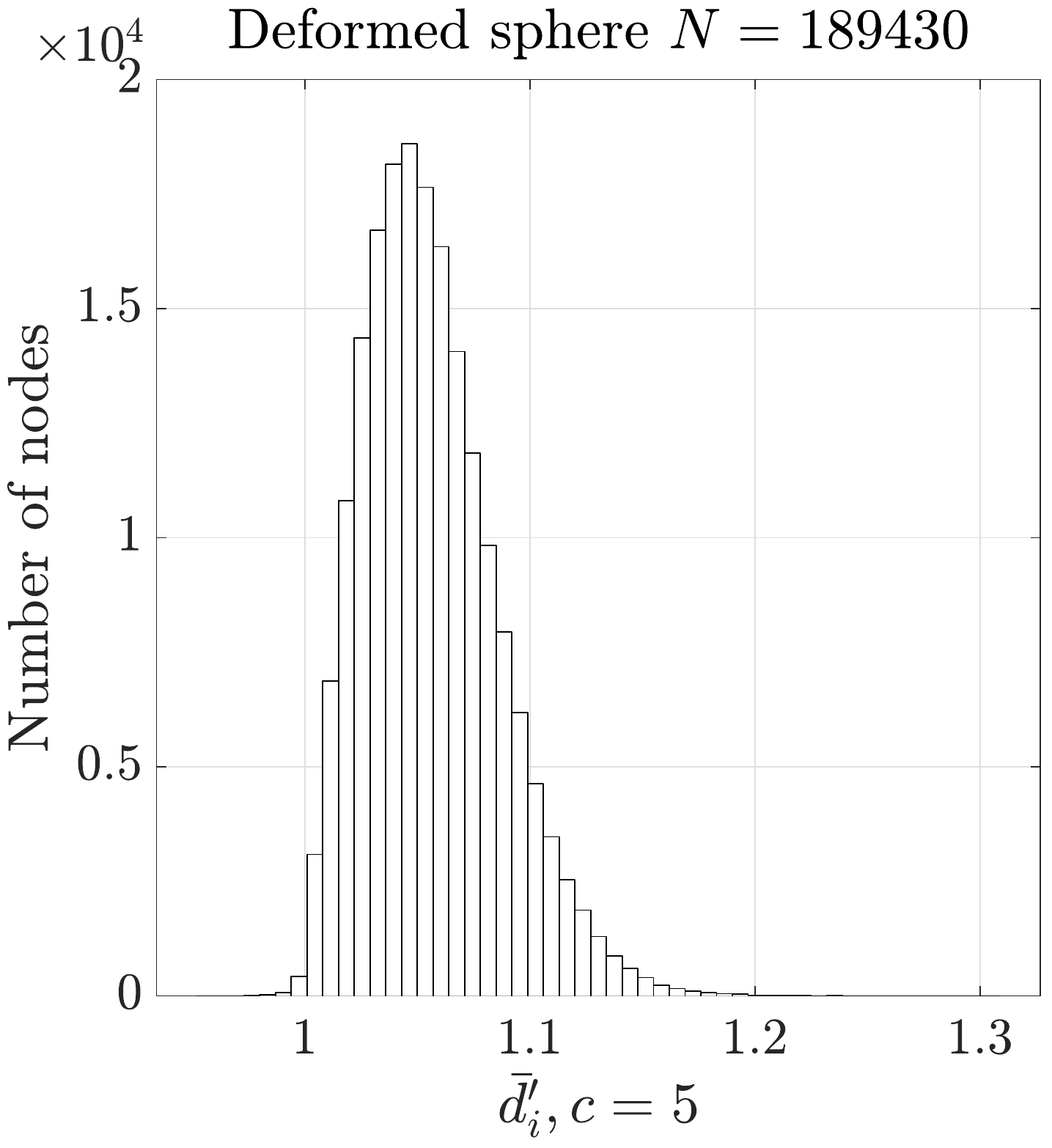}
	\caption{Local regularity distributions for boundary and interior nodes \diff{in the case of a constant spacing function}.}
	\label{fig:loc_reg_2}
\end{figure}

\diff{Since the value of $\overline d_i$ depends on the value of $c$, some reasoning is needed before we analyse our models. Sufficiently far from the
boundary, one would ideally like to consider the value of $c$ equal to the
maximal number of points that can be placed on a unit sphere with mutual
distances greater than or equal to 1 (that is 2 in 1D, 6 in 2D and 12 in
3D~{\cite{conway2010sphere}}). In practice, however, the node distributions
are not close to ideal even without the presence of a
boundary~{\cite{slak2019generation}}, which means that considering just the
ideal $c$ would fail to fairly assess the uniformity of the distribution,
especially in the case of a CAD model (where the boundary often plays an important
role). In Figure~{\ref{fig:loc_reg_of_c}}, the means and standard deviations
of $\overline{d_i'}$ computed over the whole domain (\emph{i.e.} the means and averages of
distributions later shown in Figure~{\ref{fig:loc_reg_1}} and {\ref{fig:loc_reg_2}}) are 
shown as a function of $c$ for each considered model. We see that
both statistical quantities depend on the model, the dimensionality of
the domain, and if we are considering the whole domain or only the boundary.
In general, boundaries of a given model are easier to uniformly discretize
than the interior,
since the boundaries have one dimension less than the interior. This is true
despite the fact that sDIVG uses the first order Taylor expansion to
determine the appropriate spacing, which results in candidates being
generated at spacing only approximately equal to $h$, whereas DIVG makes no
such approximation. Furthermore, a simple
argument considering only the dimensionalities cannot be sufficient for
explaining why the distributions for the 2D duck case are worse than for the 3D
boundaries (which is also a 2D object). Here, we must take into account that
the duck model has more convex vertices, where, even in the ideal case, one
cannot hope to come close to the ideal number of equidistant neighbors for the
case of the empty space. Therefore, the distribution of $\overline{d_i'}$ for a
large number of nearest neighbors $c$ fails to fairly assess the uniformity of
nodes. This effect is also later visible in Figure~{\ref{fig:loc_reg_2}}, where a
spike just after $\overline{d_i'} = 1.05$ (which can be attributed to said
vertices) is visible. The boundary of both 3D objects does not
itself have a boundary $\partial(\partial\Omega)$, which means that higher
values of $c$ give a fairer estimate of node distribution quality there. In
following analyses we therefore used $c = 2$ for 1D objects (\emph{i.e.} domain
boundaries in 2D), $c=3$ for 2D objects (domain boundaries in 3D and domain
interiors in 2D) and $c = 5$ for 3D objects (domain interiors in 3D).}

The distance distributions to nearest neighbors for boundary nodes are
presented in Figure~\ref{fig:loc_reg_1}, and Figure~\ref{fig:loc_reg_2}
shows distributions for all nodes. The quantitative statistics are
presented in Table~\ref{tb:stat_node_q}. It can be seen that the nodes are
quite uniformly distributed as all distributions are condensed near $1$. In
general, \diff{the uniformity of boundary node distribution is on par with the distribution of interior nodes.}

In the 2D duck case, the distribution \diff{of boundary nodes} visually seems much better
than in the 3D cases. This is a consequence of the candidate generation
procedure, which is more optimal when the parametric domains are 1D. Another
feature characteristic for 1D parametric domains is that outliers with distance
to nearest neighbors slightly less than $2h$ are not uncommon. This happens at
nodes where the advancing fronts of the sDIVG algorithm meet and is rarely a
problem in practice.
For these reasons, the standard deviation for duck case shown in
Table~\ref{tb:stat_node_q} is of the same order of magnitude as the 3D cases.
If we remove the 10 most extreme outliers, the standard deviation reduces by an
order of magnitude. See~\cite{duh2021fast} for a deeper analysis and possible
solutions.

In the 3D cases, the distribution of nodes for the deformed sphere is slightly
worse, which can be attributed to a greater complexity of the model.
Nevertheless, the quality of generated nodes is of the same order as the nodes
generated by pure DIVG~\cite{slak2019generation} and sDIVG~\cite{duh2021fast}.


\begin{table}[h]
	\centering
	\caption{Statistics \diff{of local regularity distributions shown in Figure~{\ref{fig:loc_reg_1}} and {\ref{fig:loc_reg_2}}}.}
	\begin{tabular}{c|c|c|c|c}
		                            &                 & $\operatorname{mean}\bar{d}'_i$ & $\operatorname{std}\bar{d}'_i$ &
		$\operatorname{mean}\left(\left(d_i^{\text{max}}\right)' -
		\left(d_i^{\text{min}}\right)'\right)$                                                                                      \\
		\hline \hline
		boundary nodes              & Duck            & $1.0007$                        & $0.017$                        & $0.0022$ \\
		\cline{2-5}
		                            & Sphere          & $1.04$                          & $0.036$                        & $0.10$   \\
		\cline{2-5}
		                            & Deformed sphere & $1.04$                          & $0.039$                        & $0.10$   \\
		\hline
		Boundary and interior nodes & Duck            & $1.036$                         & $0.030$                        & $0.101$  \\
		\cline{2-5}
		                            & Sphere          & $1.055$                         & $0.029$                        & $0.103$  \\
		\cline{2-5}
		                            & Deformed sphere & $1.056$                         & $0.030$                        & $0.140$  \\
	\end{tabular}
	\label{tb:stat_node_q}
\end{table}
Additionally, there are two quantities often considered as node quality
measures, \emph{i.e.}, minimal distance between nodes \diff{(also referred to as separation distance)} and fill distance (\diff{also referred to as the maximal empty sphere radius}) within the domain~\cite{hardin2004discretizing,wendland2004scattered}. The minimal distance is defined for set of nodes $\Xi
	= \{x_1, \dots, x_N\} \subset \Omega$ as
\begin{equation}
	r_{\mathrm{min}, \Xi} = \frac{1}{2} \min_{i \neq j} ||\mathbf x_i - \mathbf x_j||
\end{equation}
and \diff{fill distance} as
\begin{equation}
	r_{\mathrm{max}, \Xi} = \sup_{\mathbf x \in \Omega} \min_i ||\mathbf x - \mathbf x_i||.
\end{equation}
Quantity $r_{\mathrm{min}, \Xi}$ is determined by \diff{finding the nearest} neighbor
for all nodes using a spatial search structure, such as a $k$-d tree. A
$r_{\mathrm{max}, \Xi}$ is estimated numerically by sampling $\Omega$ with
higher node density and searching for the closest node among $\Xi$.

\begin{figure}[h]
	\centering
	\includegraphics[width=0.32\linewidth]{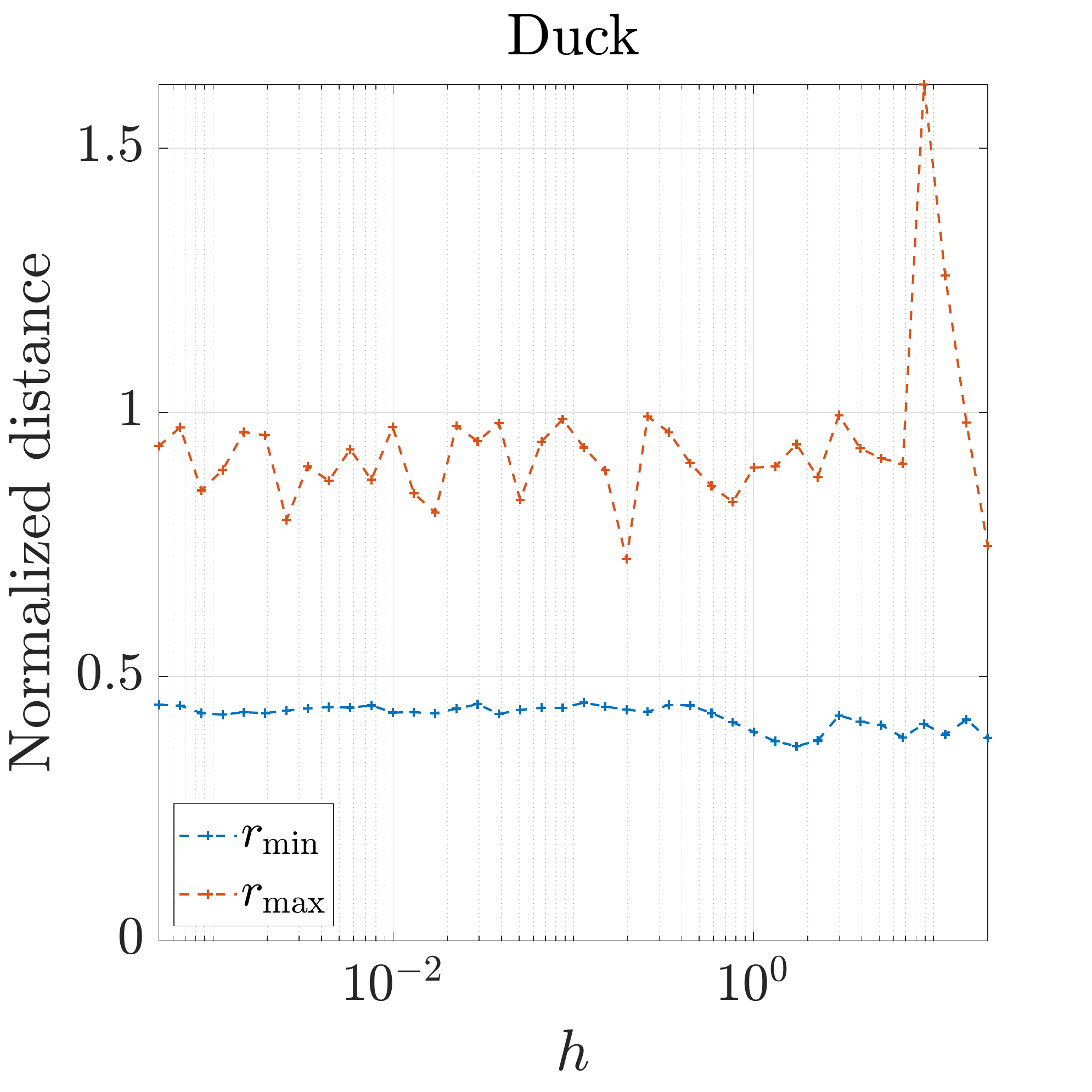}
	\includegraphics[width=0.32\linewidth]{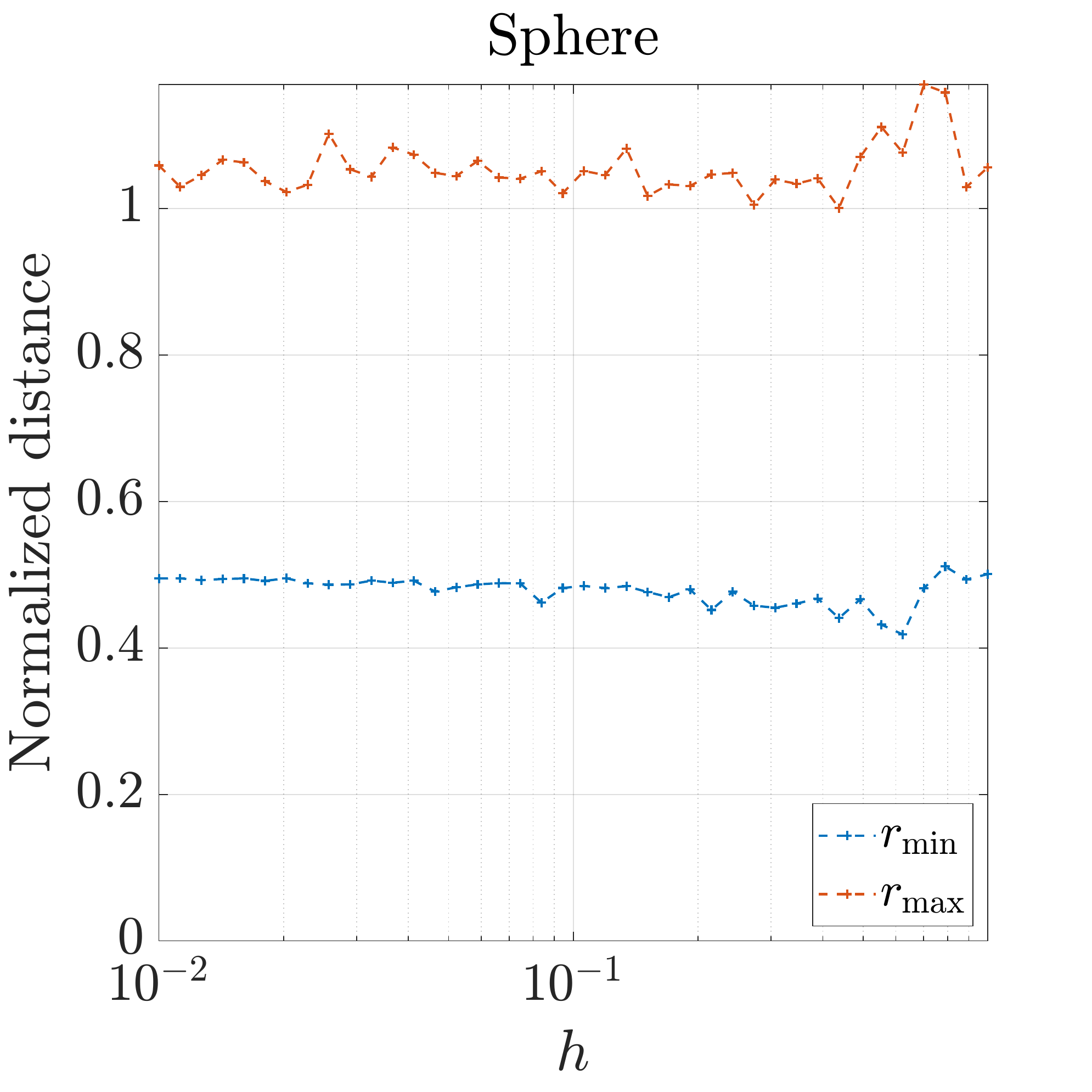}
	\includegraphics[width=0.32\linewidth]{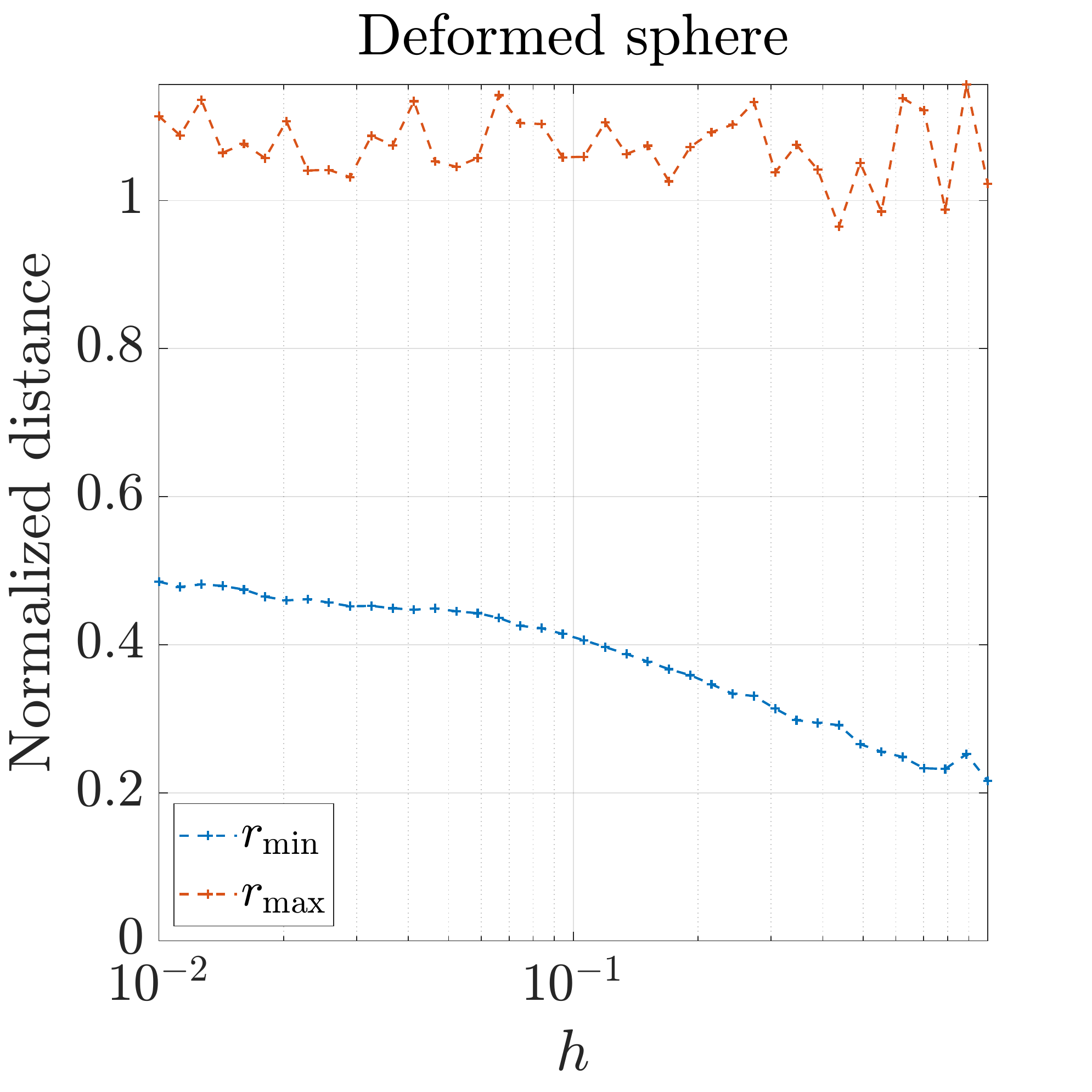}
	\caption{Minimal and \diff{fill distance on the boundary with respect to different constant values of the spacing function $h$}}
	\label{fig:loc_reg_3}
\end{figure}

\begin{figure}[h]
	\centering
	\includegraphics[width=0.32\linewidth]{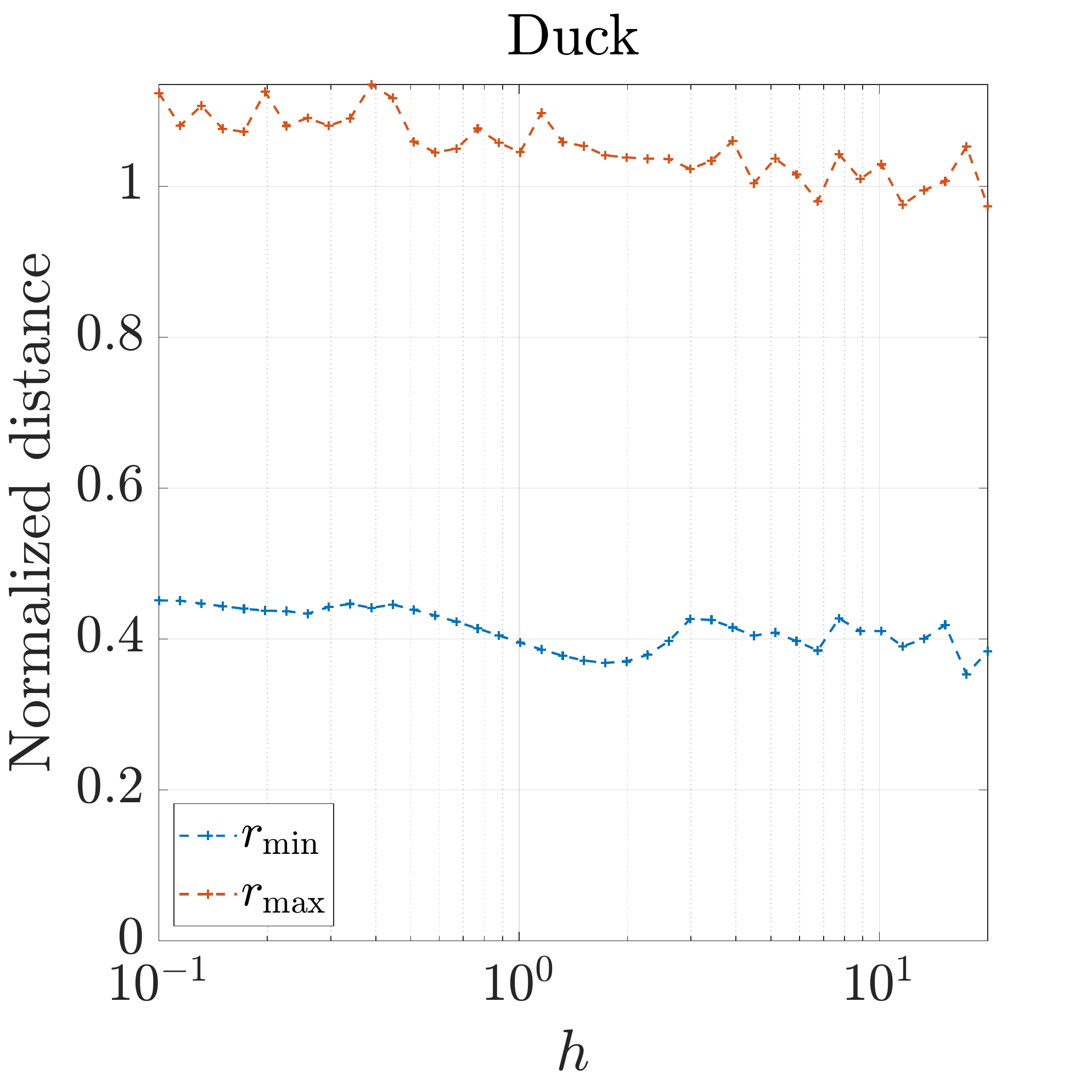}
	\includegraphics[width=0.32\linewidth]{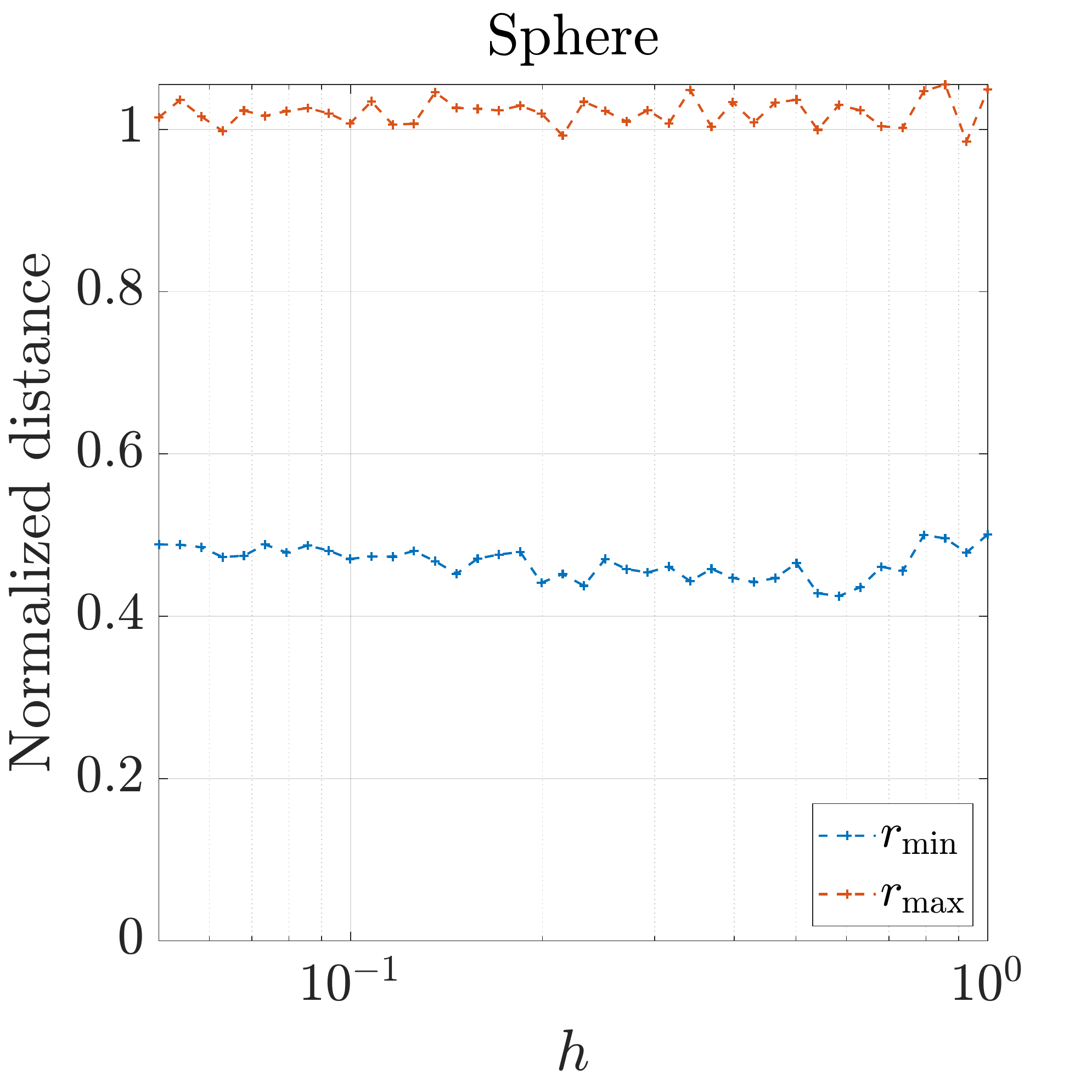}
	\includegraphics[width=0.32\linewidth]{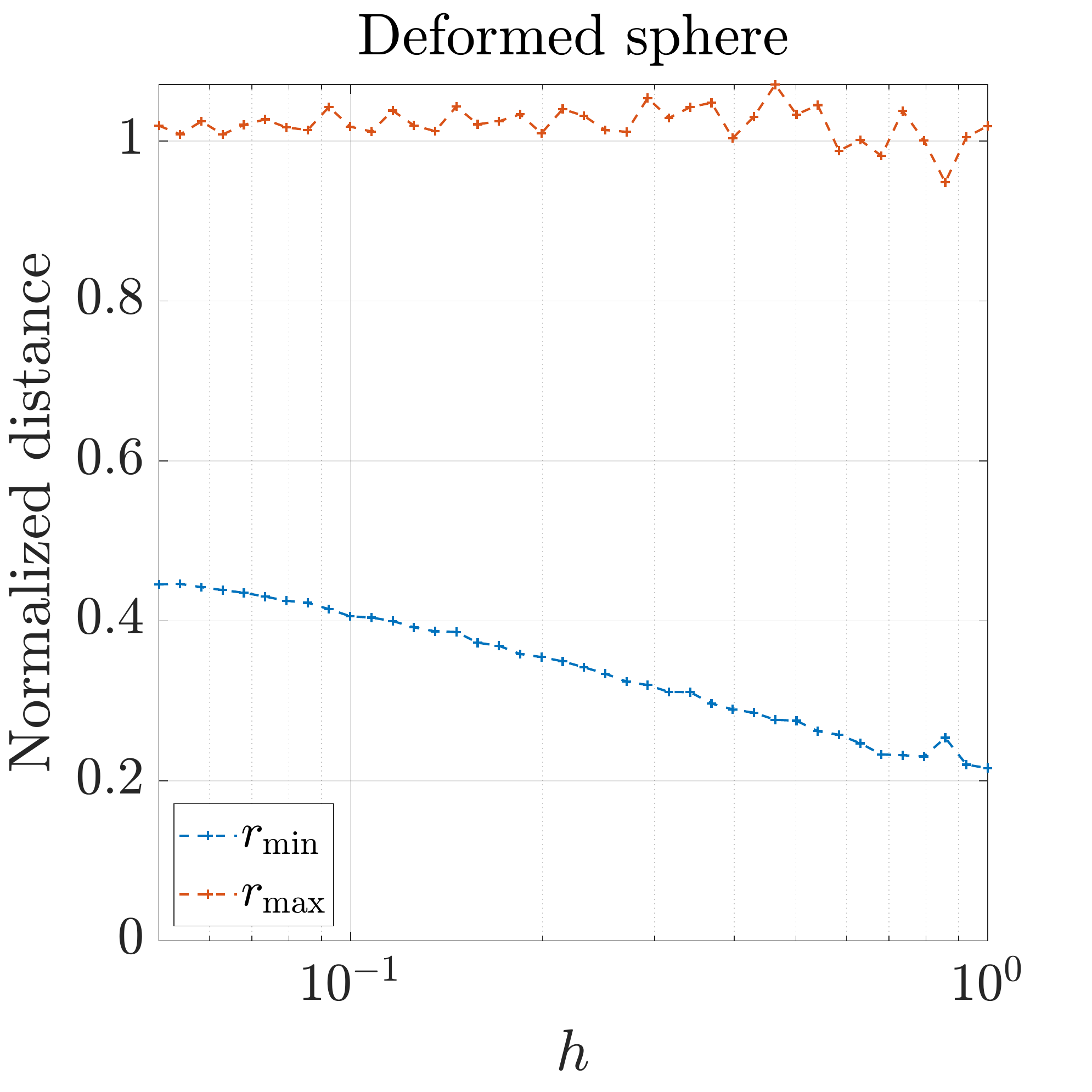}
	\caption{Minimal and \diff{fill distance on the whole domain with respect to different constant values of the spacing function $h$}}

	\label{fig:loc_reg_4}
\end{figure}

\diff{The behaviour} of \diff{the normalized fill distance} and separation distance for
all three cases with respect to target nodal distance $h$ is presented in
Figures~\ref{fig:loc_reg_3} and~\ref{fig:loc_reg_4}. In all cases,
$r_\mathrm{max}$ is relatively stable near an acceptable value of $1$ and
$r_\mathrm{min}$ approaches the optimal value of $0.5$ with decreasing $h$. This
behaviour is consistent with previous results and the analytical bound for
$r_\mathrm{min}$ for sDIVG~\cite{slak2019generation,duh2021fast}.

\section{Solving PDEs on CAD geometry}
\label{sec:PDE}
\diff{In this section, we focus on solving PDEs on domains discretized with scattered nodes $\mathbf x_i$ using the new NURBS-DIVG algorithm. 
In each node $\mathbf x_i$, the partial differential operator $\L$ is approximated using a set of $n$ nearest nodes, commonly referred to as support domain or stencil, as
}
\begin{equation}
	\label{eq:meshless}
	\L u(\mathbf x_i) \approx \sum_{j=1}^n w_j u(\mathbf x_{i,j}),
\end{equation}
\diff{where index $j$ runs over the stencil nodes of a node $\mathbf x_i$, $\mathbf w$ are
weights still to be determined and $u(\mathbf x_{i,j})$ stands for the function $u$
evaluated at the $j$-th stencil node of the node $\mathbf x_i$. The weights are determined by solving a linear
system resulting from enforcing the equality of the equation~{\eqref{eq:meshless}} for the set of approximation basis functions. In our case, the basis consists of polyharmonic
splines (PHS)~{\cite{bayona2017role}} that are centered at the stencil nodes, augmented with polynomials up to
order $m$. Such a setup corresponds to a meshless method
commonly referred to as the Radial basis function-generated finite
differences (RBF-FD)}~\cite{tolstykh2003using,jancic_monomial_2021,sabine_le_borne_potential_2021,bayona2019insight,bayona2017role}.
\diff{For the purposes of this work, we used the RBF-FD implementation
discussed in}~\cite{slak2021medusa,medusa} with augmentation up to order $m \in \{2, 4, 6\}$ on $n = 4 \binom{m + 2}{2}$ closest nodes in 2D and $n = 4 \binom{m + 3}{3}$ in 3D to obtain the mesh-free approximations of
the differential operators involved. 

\subsection{Poisson's equation}
\begin{figure}[h]
	\centering
	\includegraphics[width=0.32\linewidth]{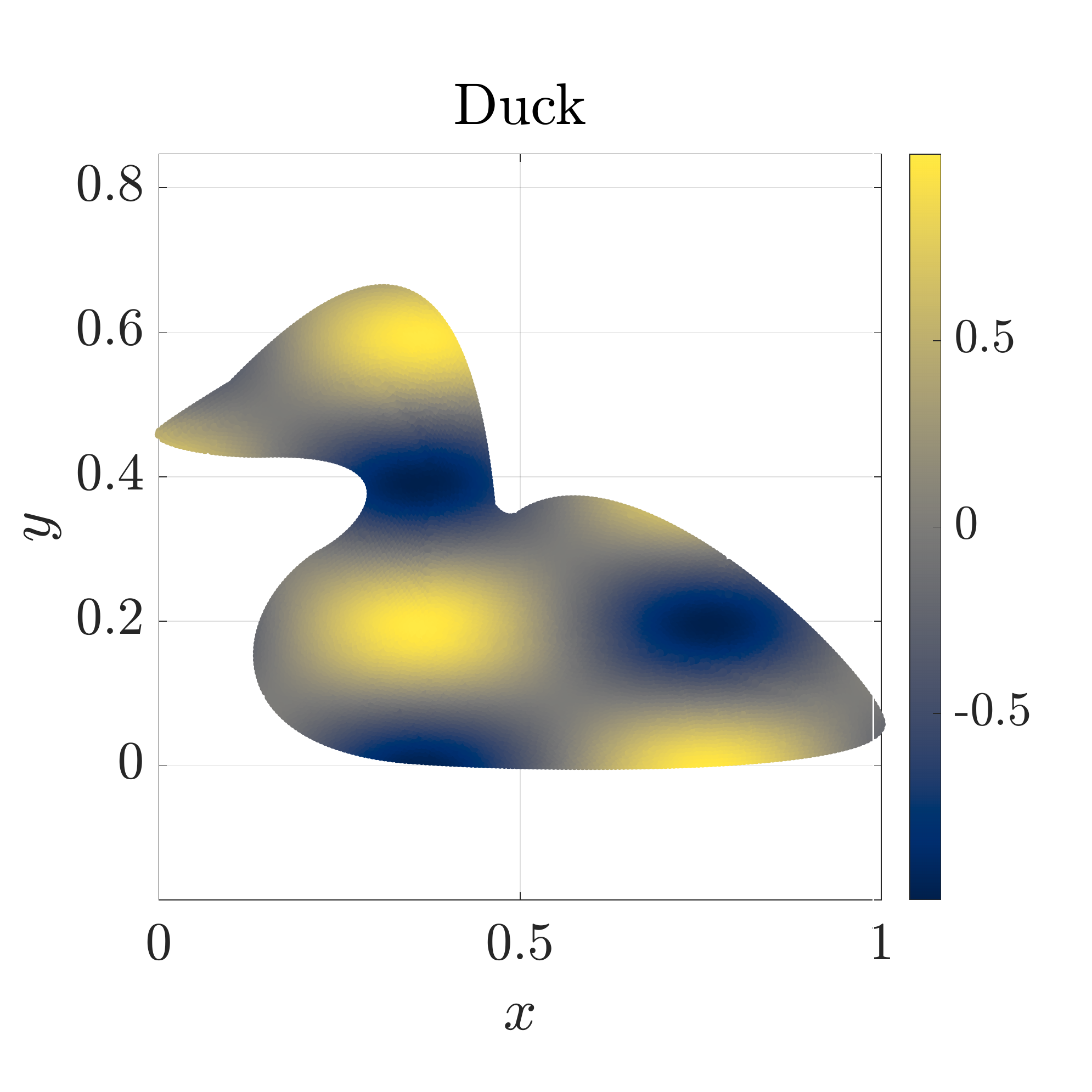}pdeso
	\includegraphics[width=0.32\linewidth]{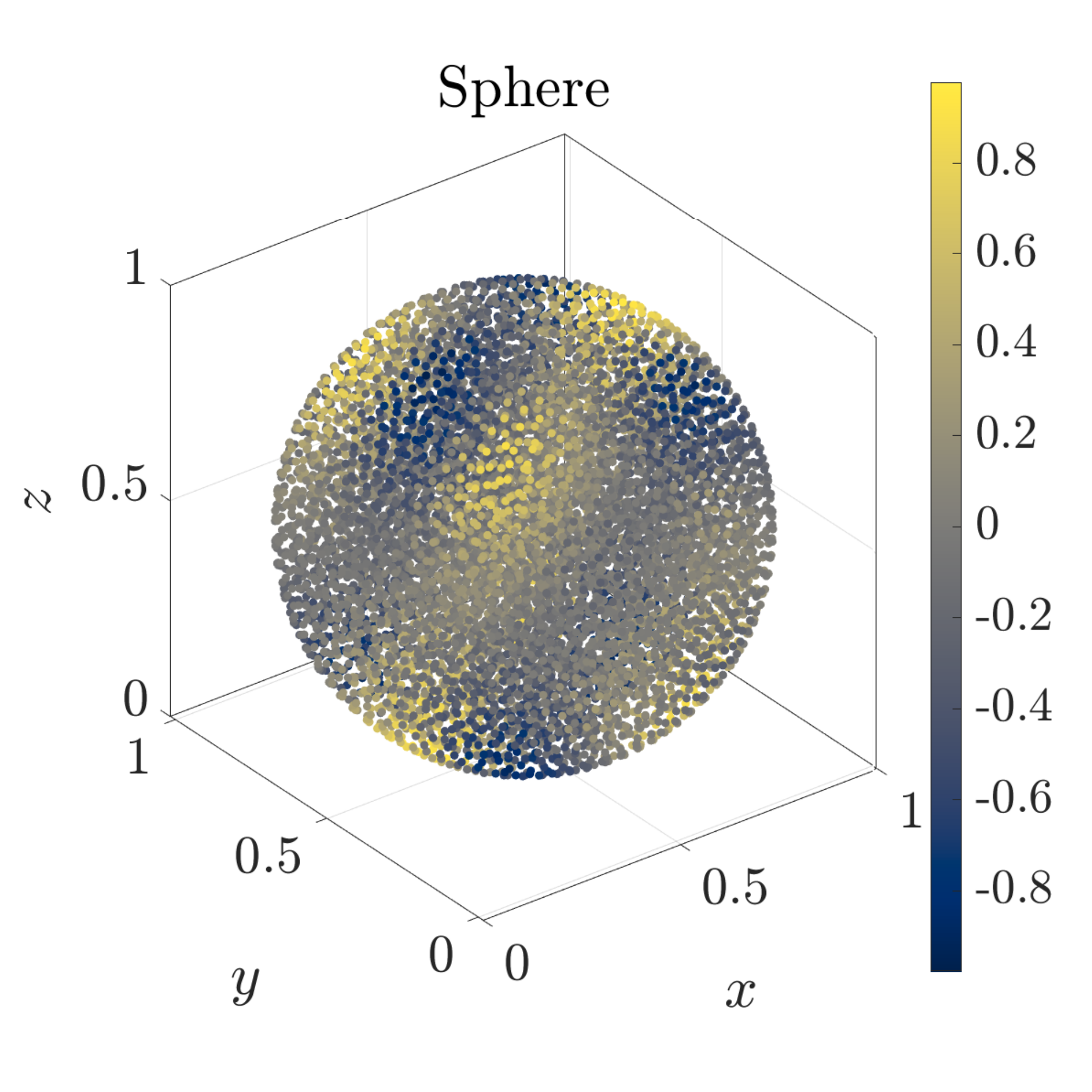}
	\includegraphics[width=0.32\linewidth]{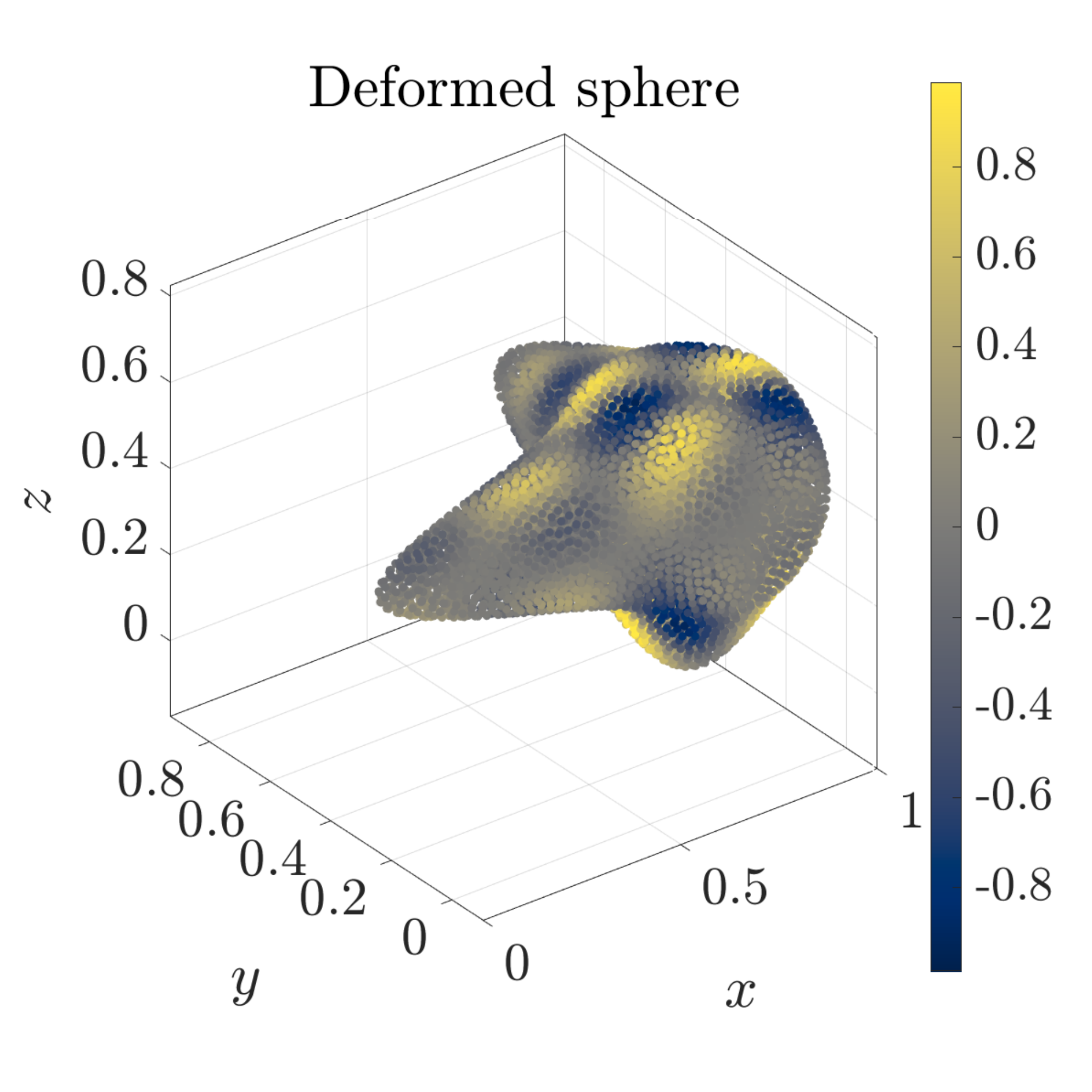}

	\caption{Solution of Poisson's equation on all three test geometries.}
	\label{fig:poiss_sol}
\end{figure}

First, we solve the Poisson's equation 
\begin{equation}
	\nabla^2 u = f
\end{equation}
with a known closed-form solution  
\begin{align}
	u_a(x, y) & = \sin\left(\frac{\pi}{100}x\right) \cos\left(\frac{2\pi}{100} y\right), \quad \text{in 2D},\\
	u_a(x, y, z) & = \sin(\pi x) \cos(2 \pi y) \sin(0.5 \pi z), \quad \text{in 3D,}
\end{align}
using mixed Neumann-Dirichlet boundary conditions
\begin{align}
	u &= u_a, \quad \text{on } \Gamma_d,\\
	\frac{\partial{u}}{\partial\mathbf{n}} &= \frac{\partial{u_a}}{\partial\mathbf{n}}, \quad \text{on } \Gamma_n, 
\end{align}
where $\Gamma_d$ and $\Gamma_n$ stand for Dirichlet and Neumann boundaries. In all cases, the domain boundary is divided into two halves, where we apply a Dirichlet boundary condition to one half and a Neumann boundary condition to another. The numerical solution of the problem is presented in Figure~\ref{fig:poiss_sol}.

\begin{figure}[h]
	\centering
	\includegraphics[width=0.32\linewidth]{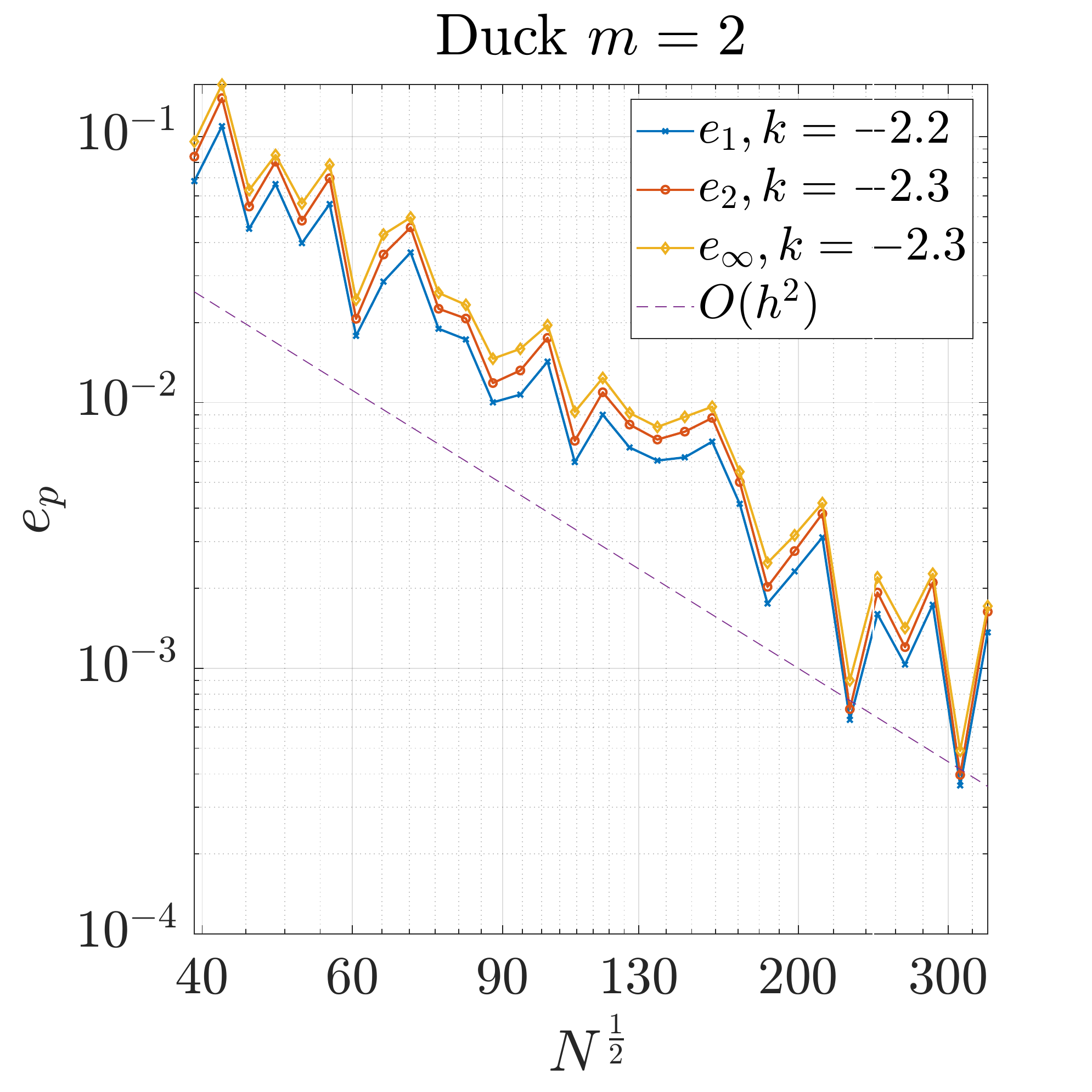}
	\includegraphics[width=0.32\linewidth]{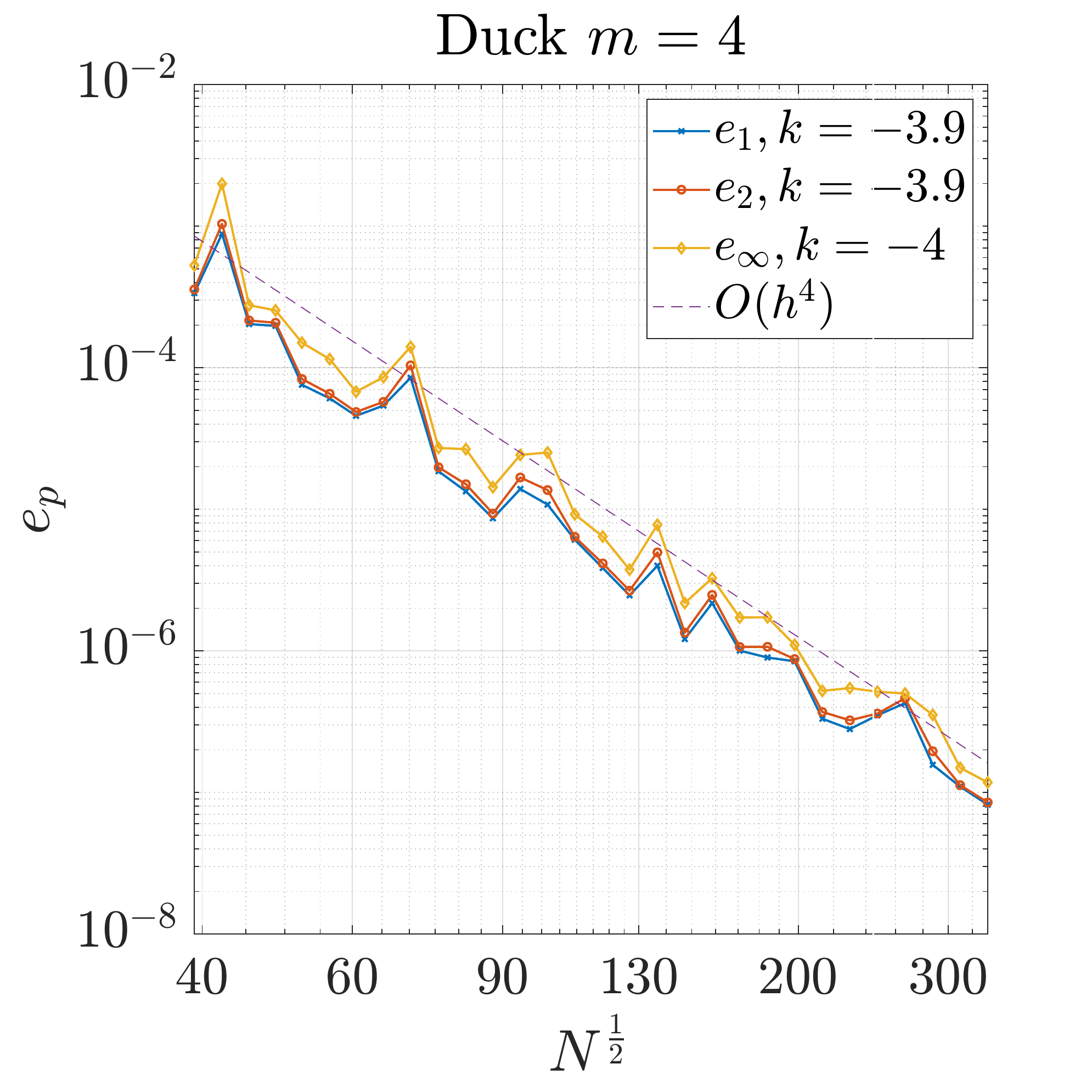}
	\includegraphics[width=0.32\linewidth]{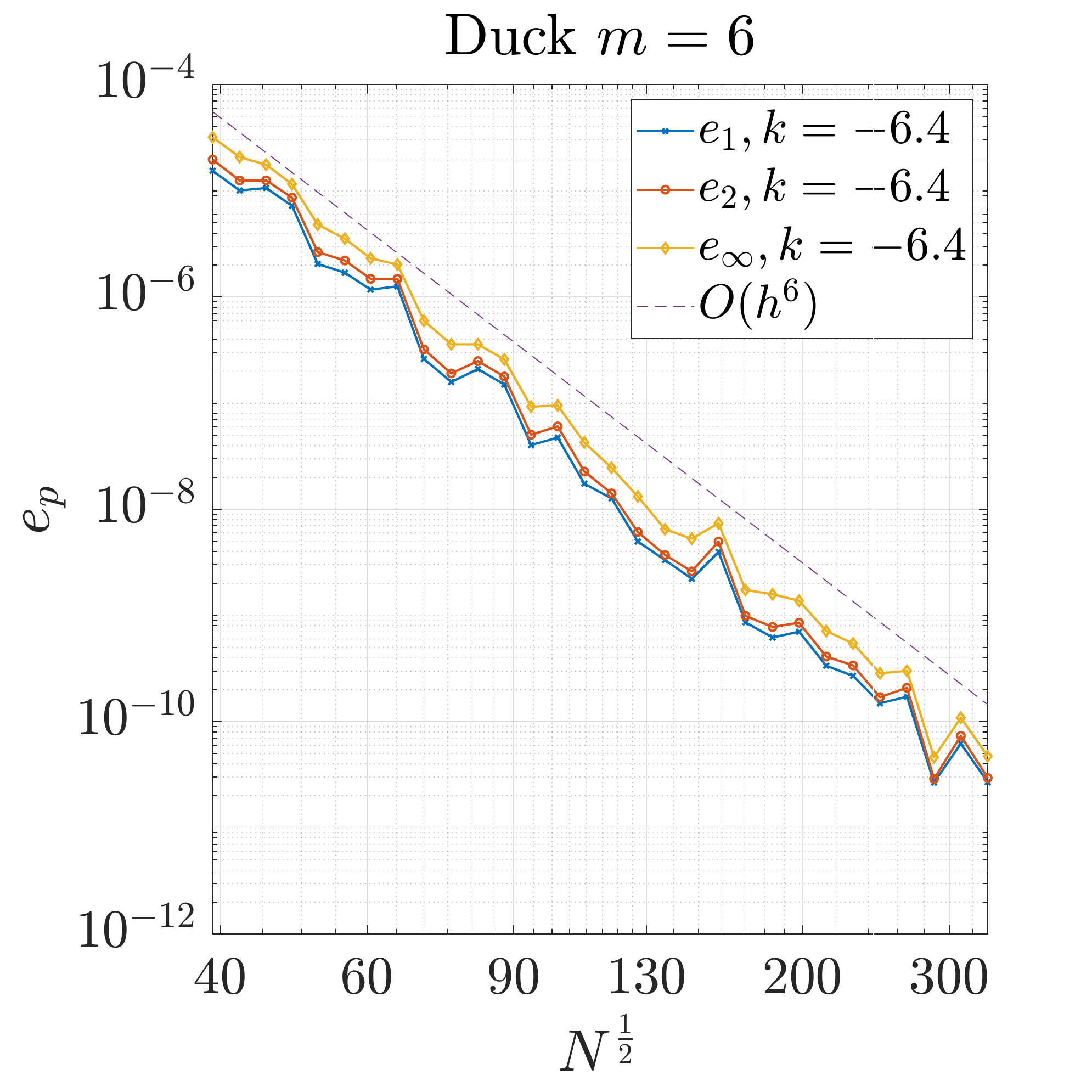}
	
	\includegraphics[width=0.32\linewidth]{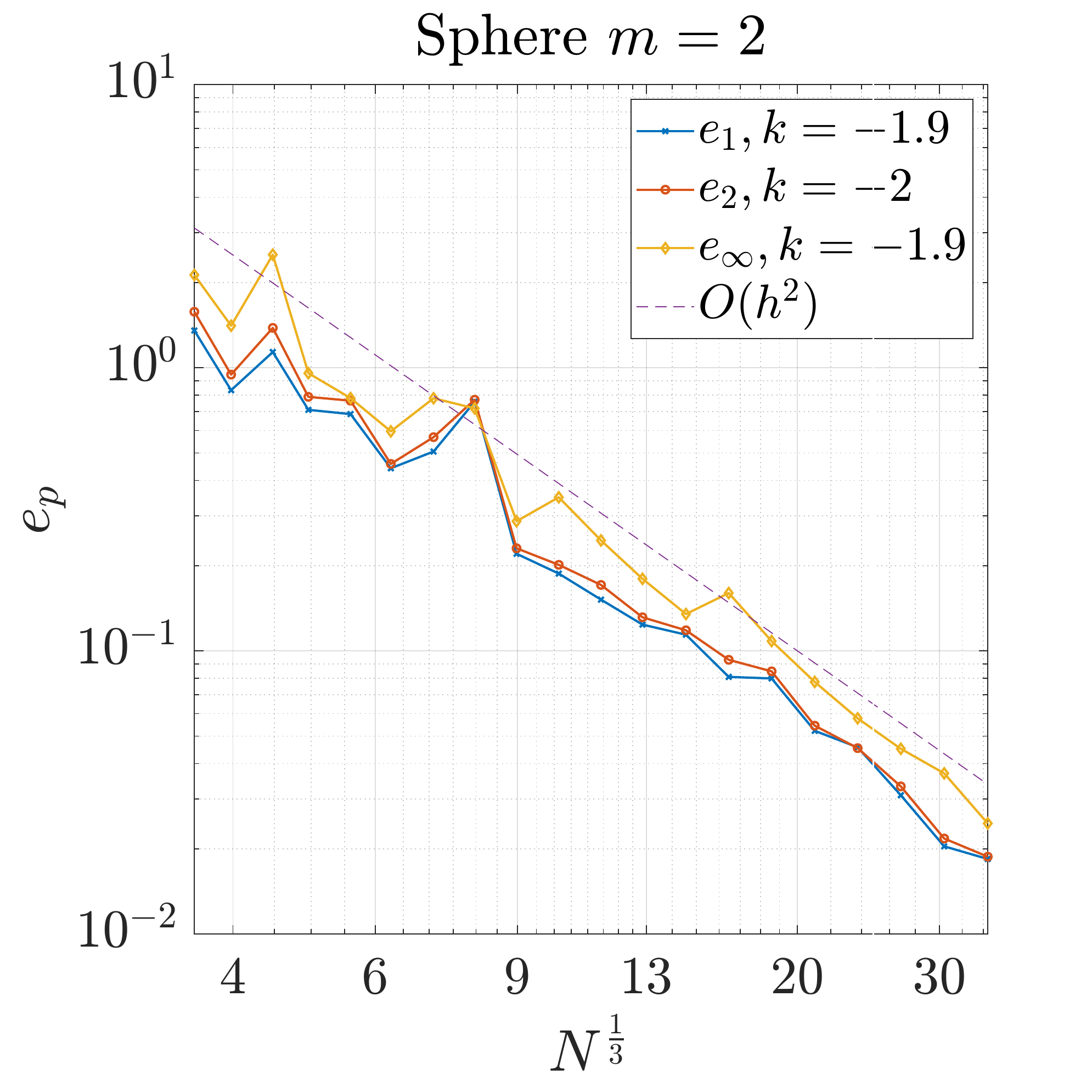}
	\includegraphics[width=0.32\linewidth]{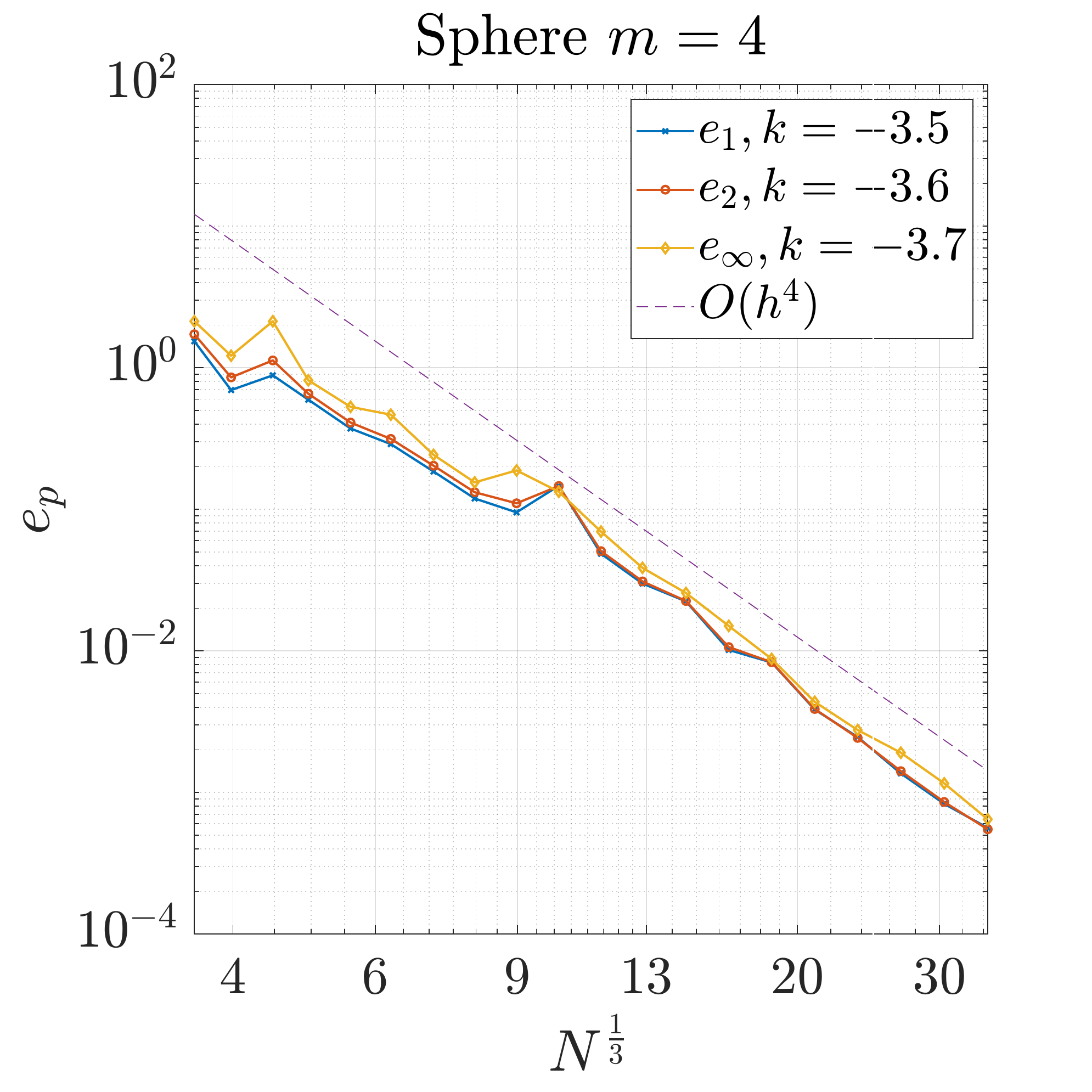}
	\includegraphics[width=0.32\linewidth]{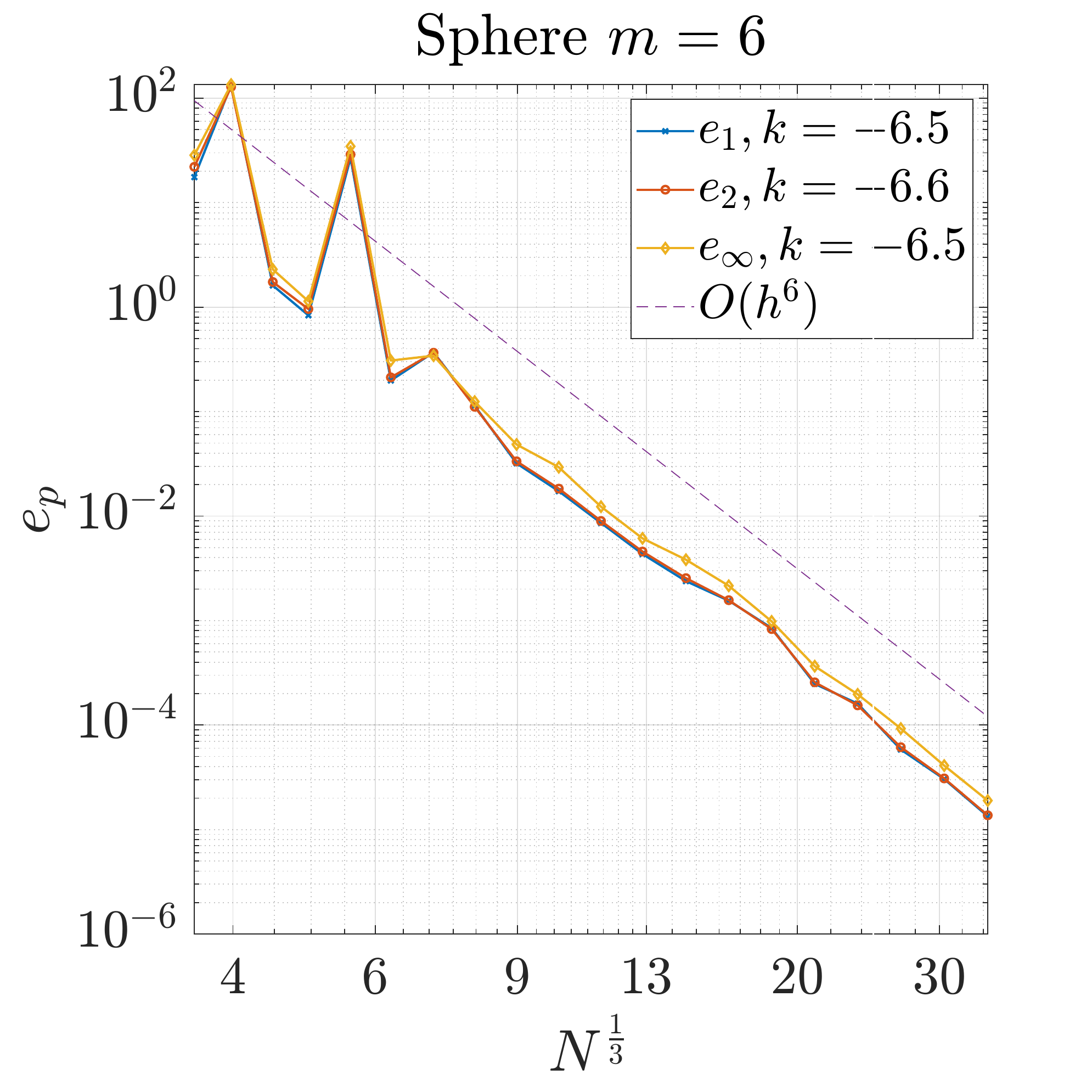}
	
	\includegraphics[width=0.32\linewidth]{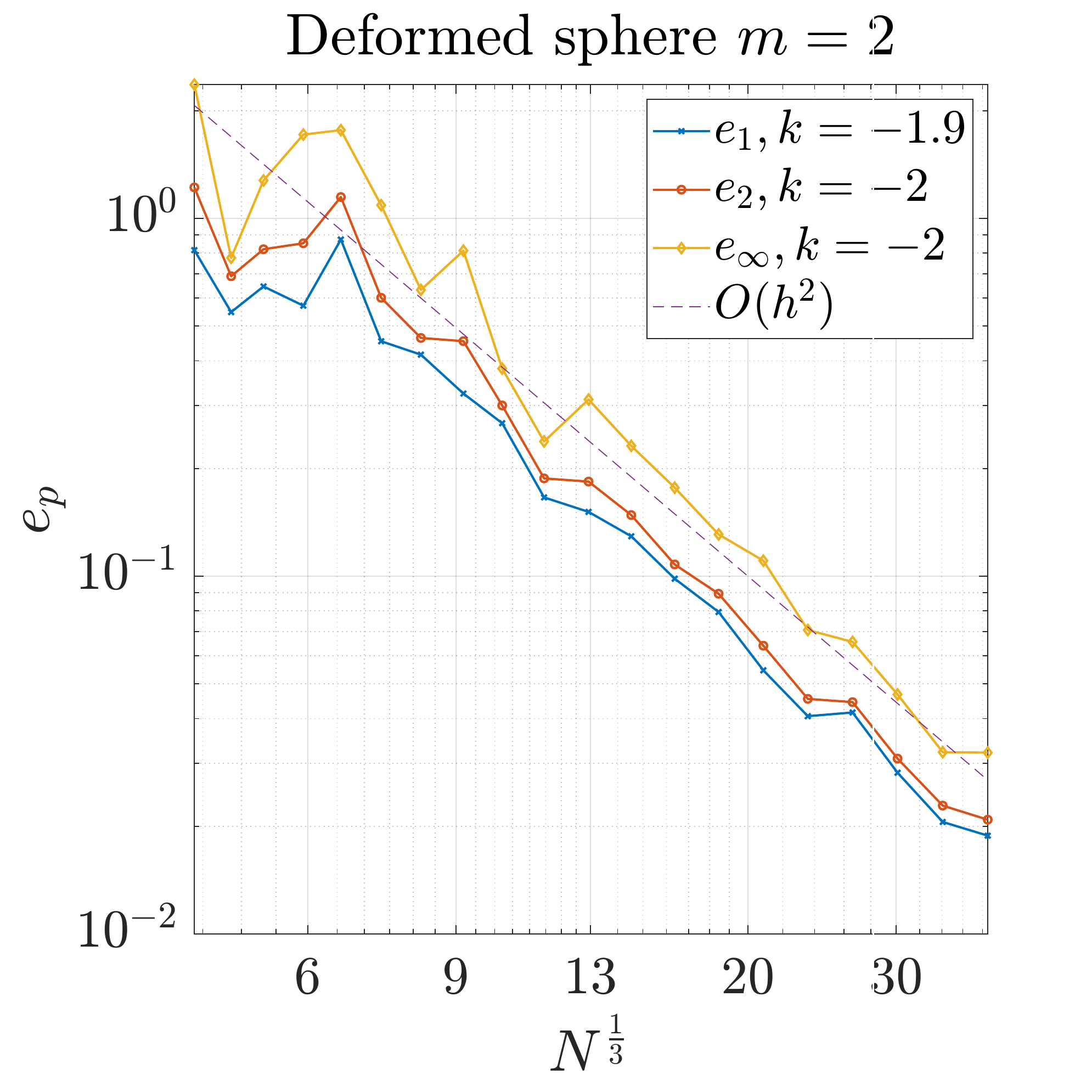}
	\includegraphics[width=0.32\linewidth]{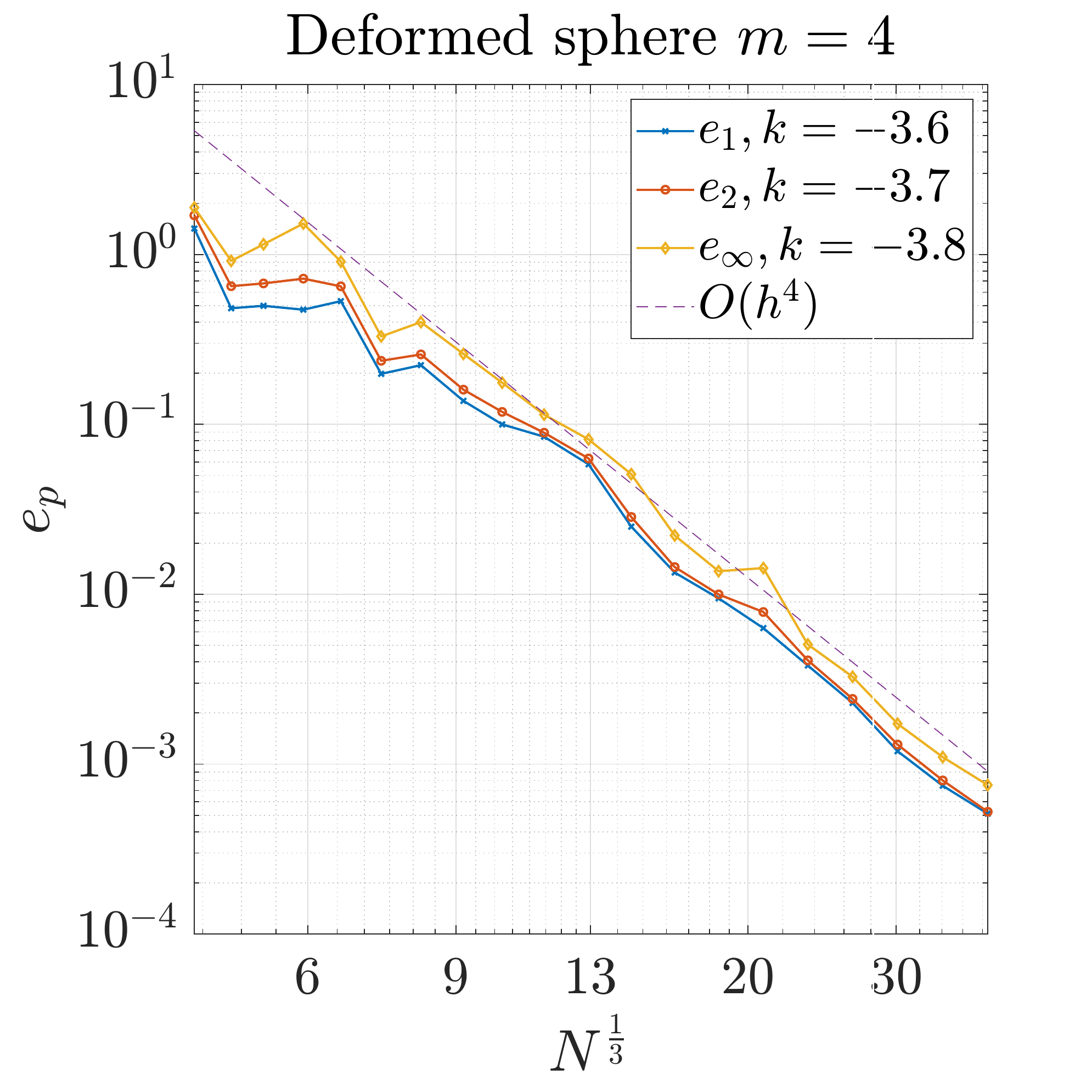}
	\includegraphics[width=0.32\linewidth]{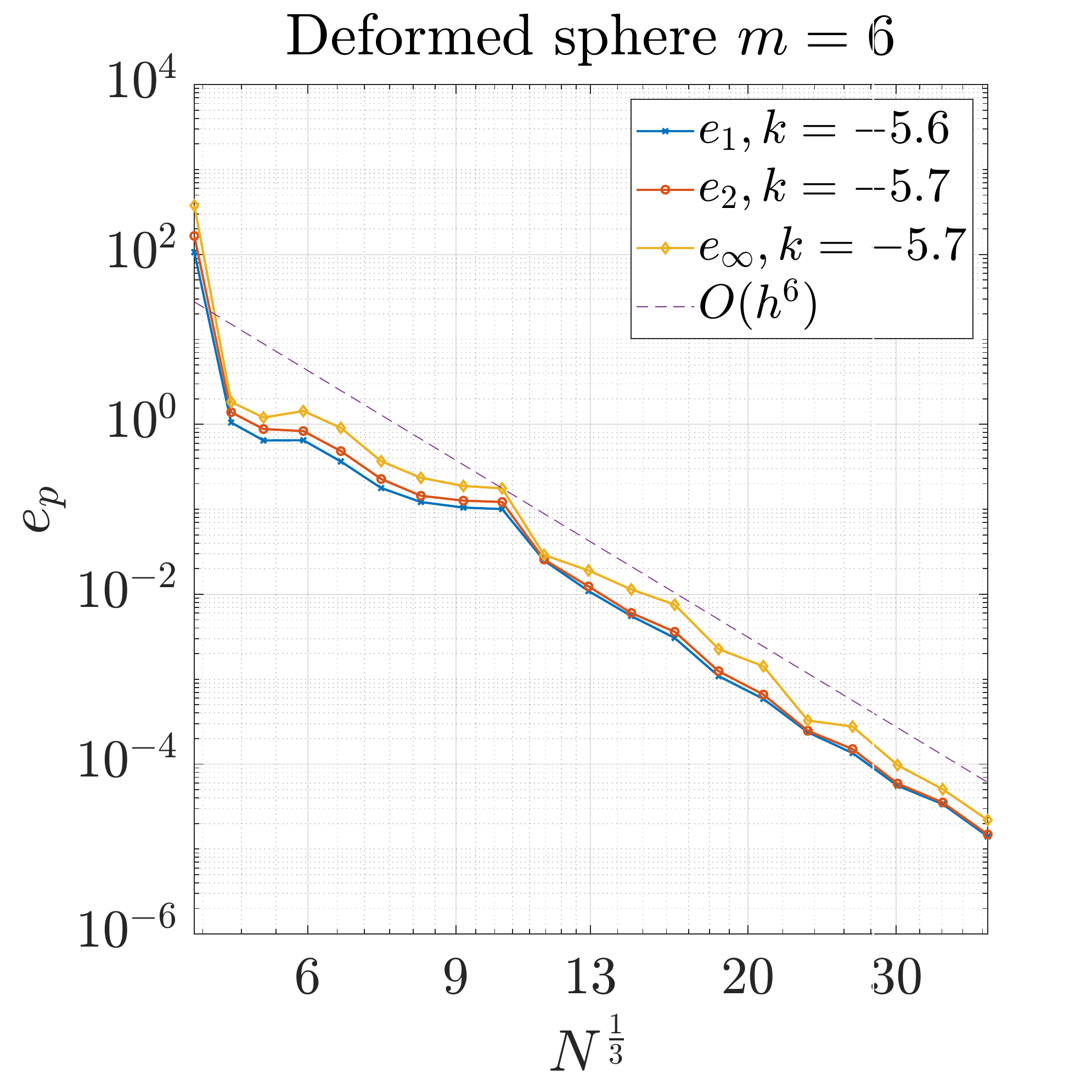}
	\caption{Error in the solution to Poisson's equation with respect to the number of nodes.}
	\label{fig:poiss_conv}
\end{figure}

Once the numerical solution $\hat{u}$ is obtained, we observe the convergence behaviour of the solution through error norms defined as 
\begin{align}
	e_1   & = \frac{\|\hat{u} - u_a\|_1}{\|u_a\|_1}, \quad \|u_a\|_1 = \frac{1}{N} 
	\sum_{i=1}^N |u_a^i|,
	\label{eq:e1}                                                               \\
	e_2   & = \frac{\|\hat{u} - u_a\|_2}{\|u_a\|_2}, \quad \|u_a\|_2 = \sqrt{\frac{1}{N} \sum_{i=1}^N |u_a^i|^2},
	\label{eq:e2}                                                               \\
	e_{\infty} & = \frac{\|\hat{u} - u_a\|_\infty}{\|u_a\|_\infty}, \quad \|u_a\|_\infty = \max_{i=1, \ldots, N}|u_a^i|.
	\label{eq:einf}
\end{align}
In Figure~\ref{fig:poiss_conv} we can see that for all three geometries the
solution converges with the expected order of accuracy according to the order of augmenting monomials. 

\begin{figure}[h]
	\centering
	\includegraphics[width=0.32\linewidth]{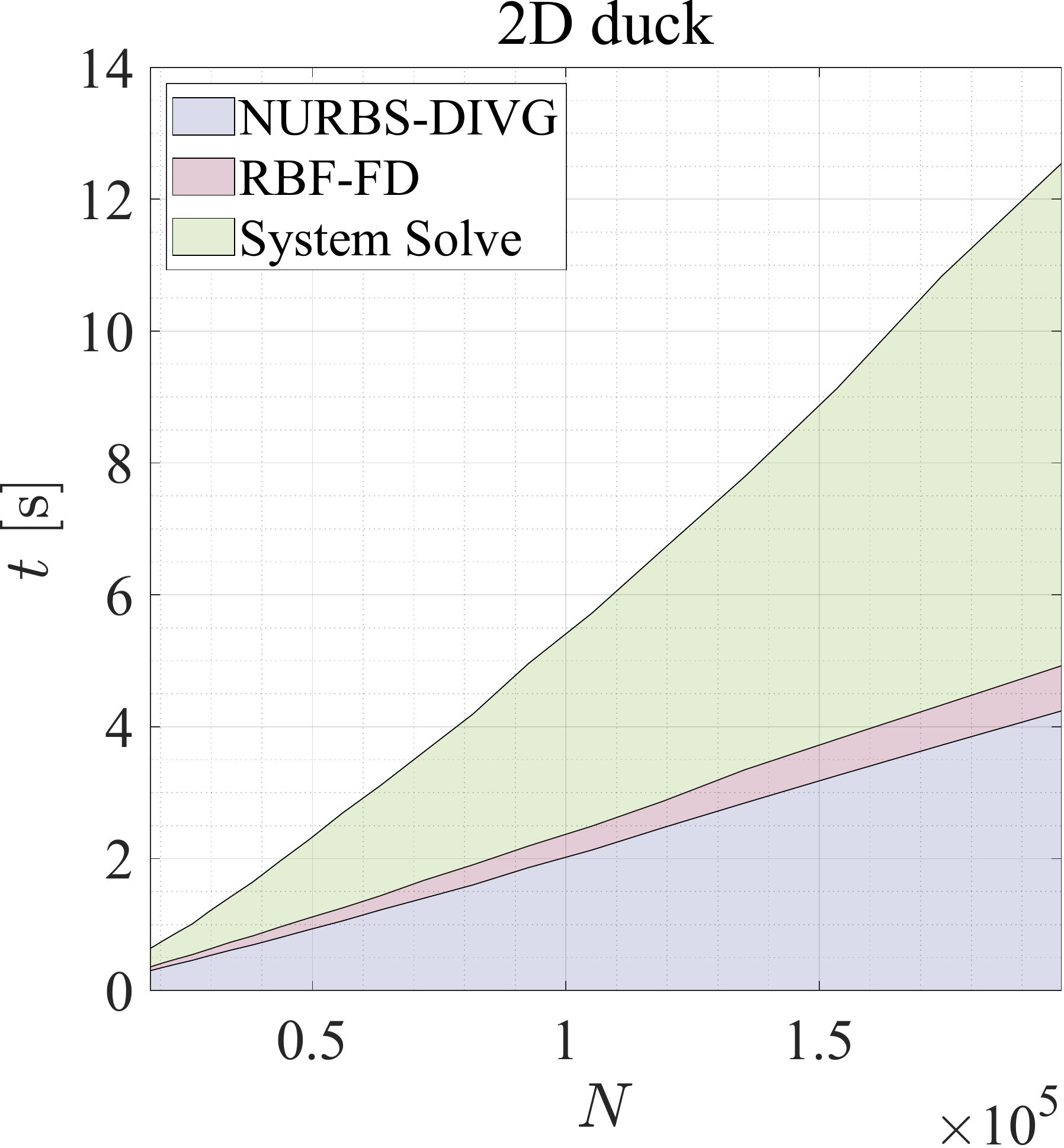}
	\includegraphics[width=0.32\linewidth]{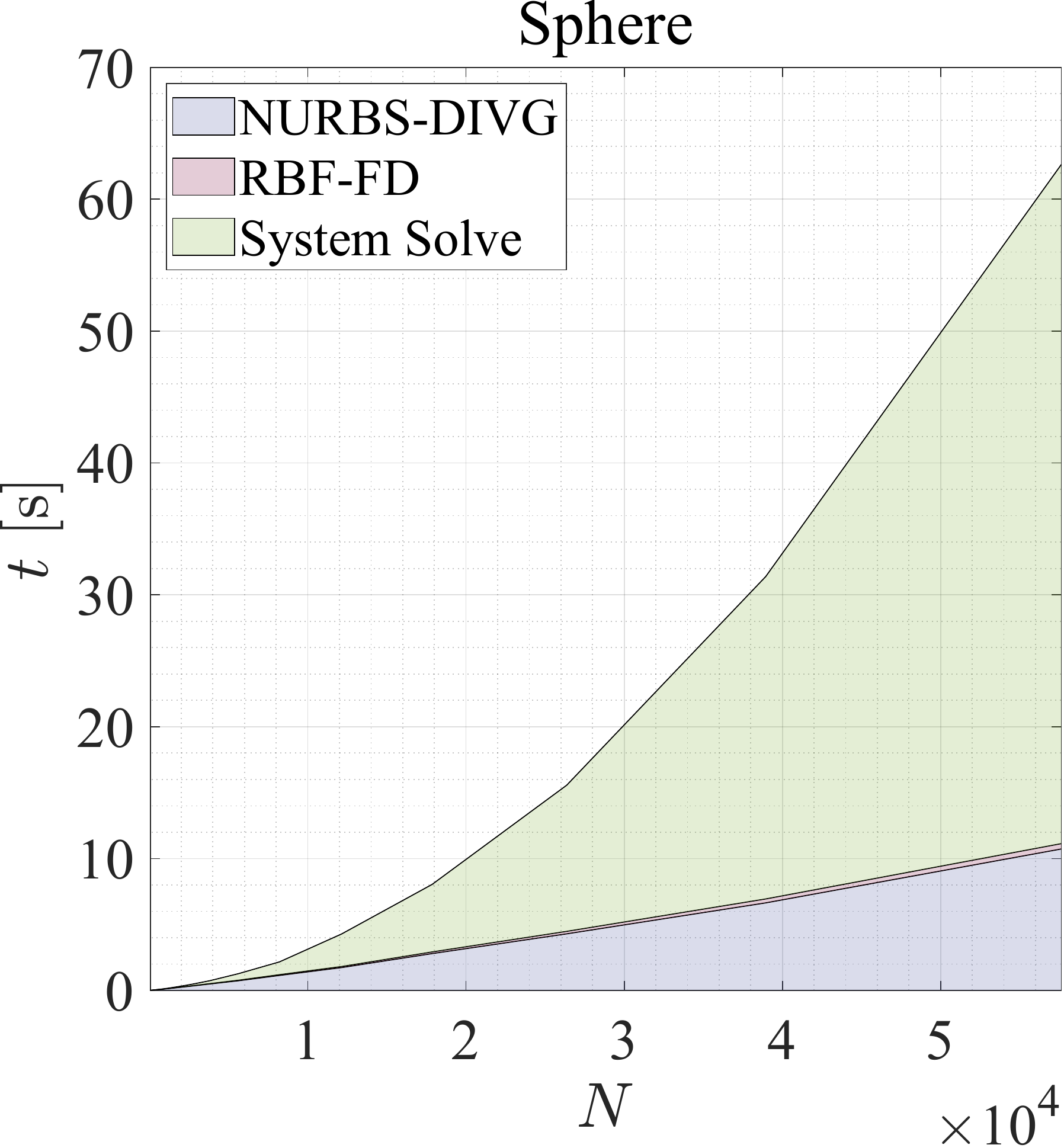}
	\includegraphics[width=0.32\linewidth]{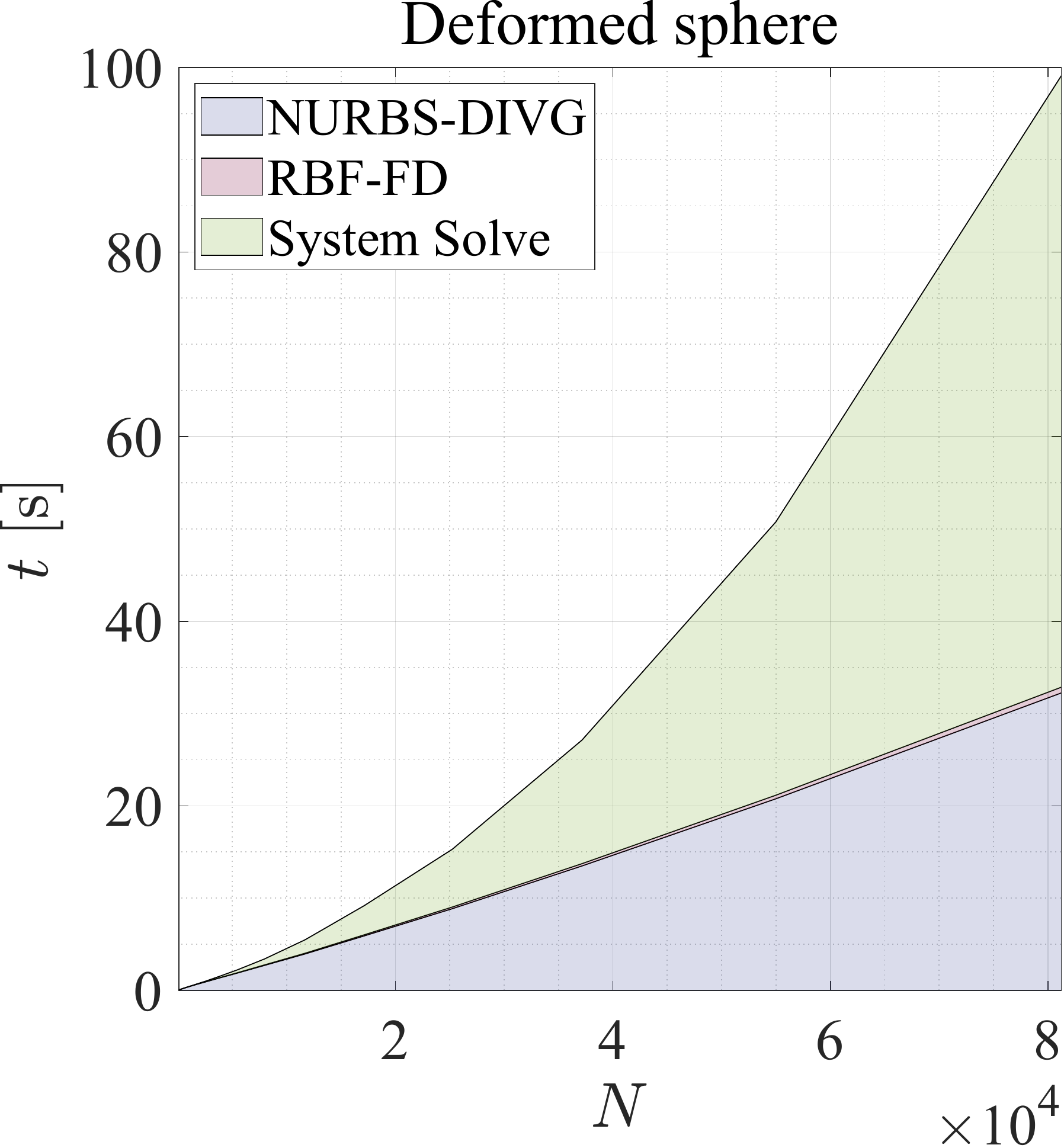}
	\caption{Execution times broken down to separate solution procedure modules.}
	\label{fig:exec_time}
\end{figure}

Next, we assess the execution time of solving Poisson's equation with second
order monomial augmentation. In Figure~\ref{fig:exec_time} the execution times
for all three geometries are broken down to core modules of the solution
procedure. We measure execution time for generation of nodes using proposed NURBS-DIVG algorithm\footnote{For a more in-depth
analysis of computational complexity, see the DIVG~\cite{slak2019generation} and
sDIVG~\cite{duh2021fast} papers.}. 
\diff{The} generation of stencils and \diff{the computation} of stencil weights are measured together as \diff{the} RBF-FD part
of the solution procedure. Separately, we also measure \diff{the cost of} sparse matrix assembly
(which is negligible~\cite{depolli2022parallel}) and \diff{the solution of the corresponding linear system}; with an
increasing number of nodes, \diff{this solve} ultimately dominates the execution
time~\cite{depolli2022parallel}. In \diff{the} 2D duck case, the computational times of
solving the system and filling the domain are of the same order
as the number of nodes is still relatively small. However, in both 3D cases we
can clearly see that \diff{the cost of} solving \diff{the linear} system scales super-linearly and soon dominates the
overall computational cost. The RBF-FD part (as expected) scales almost linearly
(neglecting the $\mathcal{O}(N\log{N})$ resulting from $k$-d tree in stencil
selection). 

\subsection{Linear elasticity - Navier-Cauchy equation}
\begin{figure}[h]
	\centering
	\raisebox{+0.1\height}{
	\includegraphics[height=0.33\linewidth]{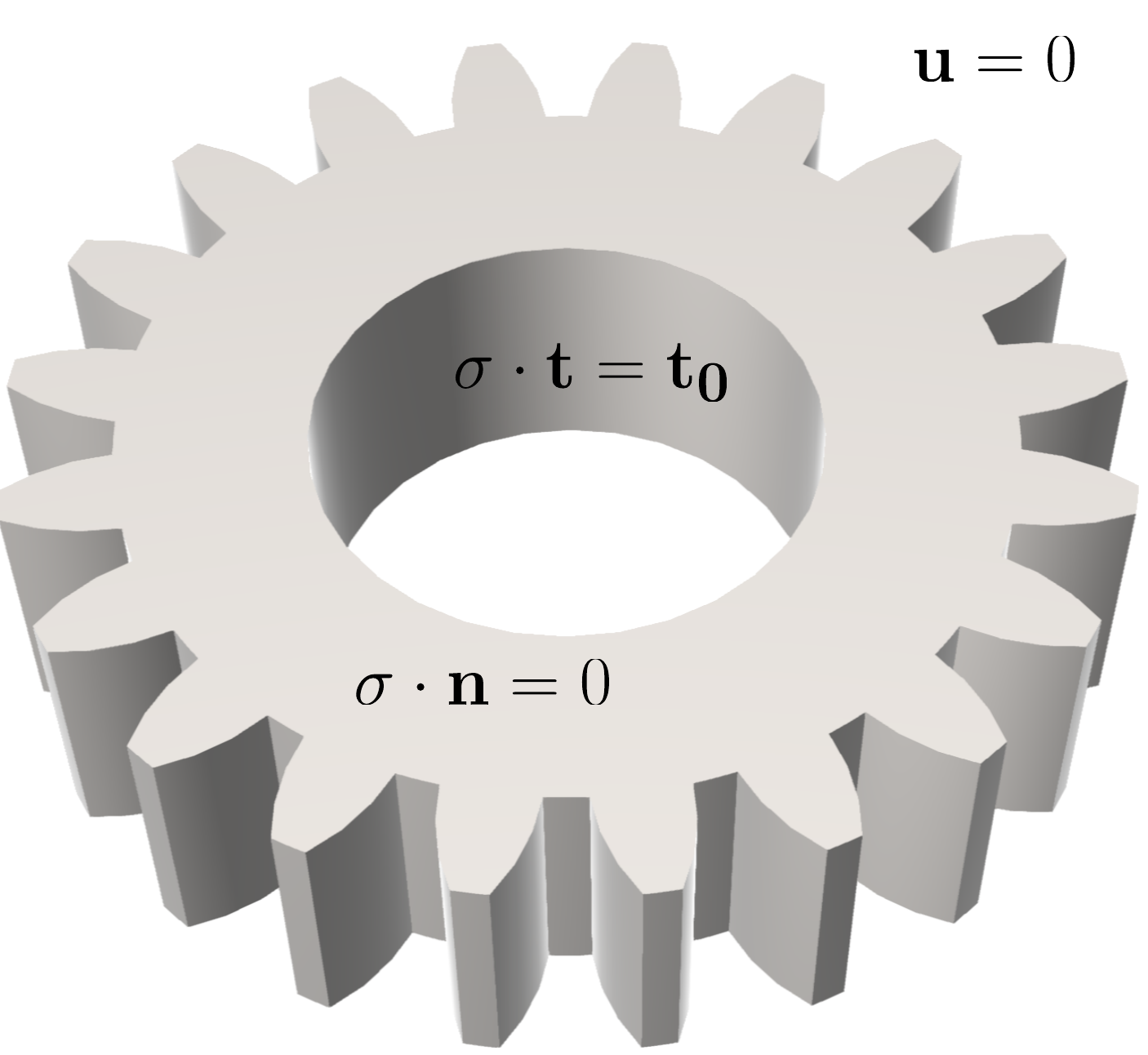}}
	\centering
	\includegraphics[trim={3cm 1cm 0 1cm},clip,  height=0.40\linewidth]{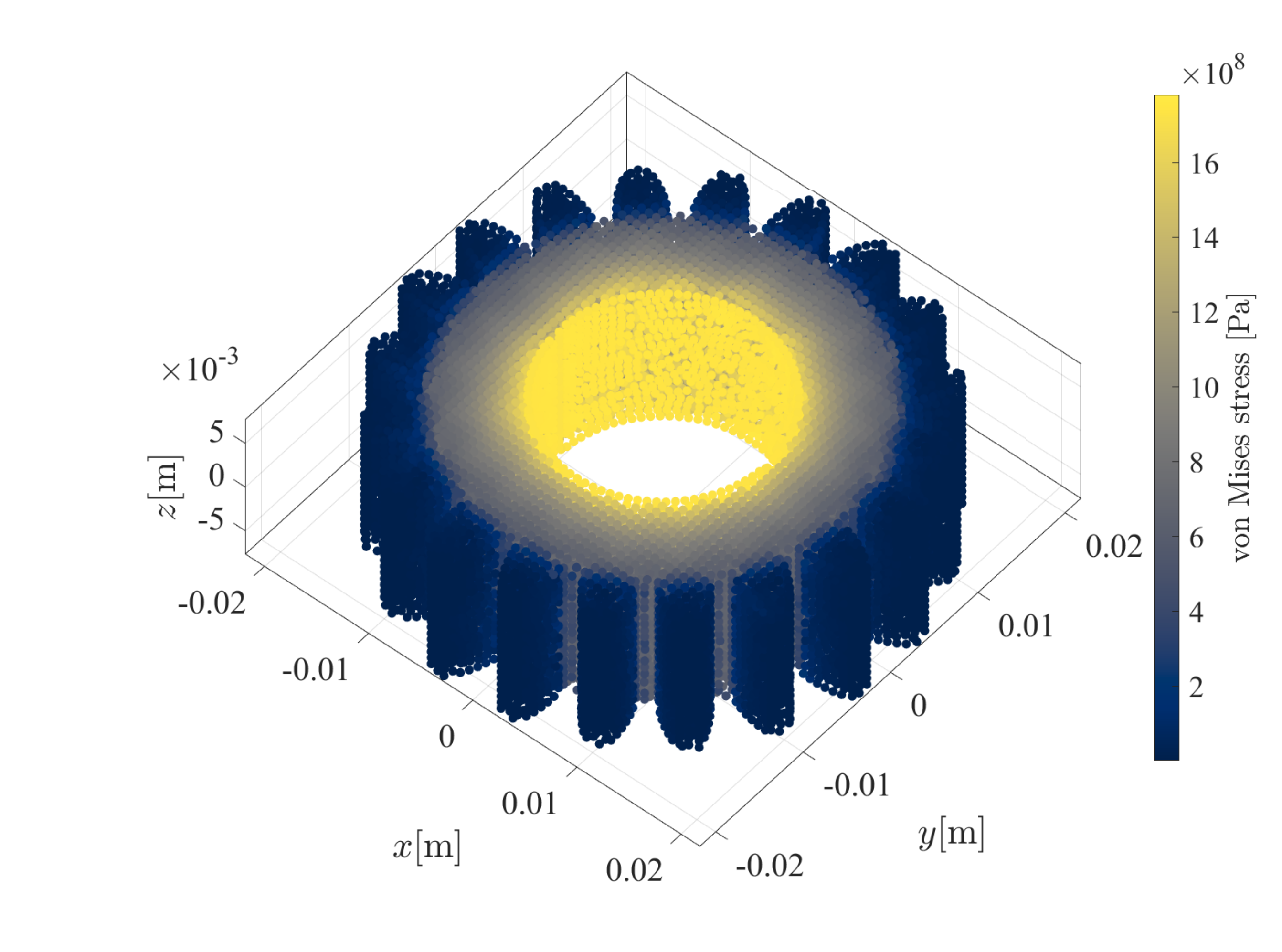}
	
	\caption{Scheme of the linear elasticity example (left) accompanied with the RBF-FD solution in terms of von Misses stress (right).
  \diff{The gear model is made of 84 patches.}}
	\label{fig:gear_scheme}
\end{figure}

\begin{figure}[H]
	\centering
	\includegraphics[trim={0cm 1cm 0cm 1cm},clip,  height=0.40\linewidth]{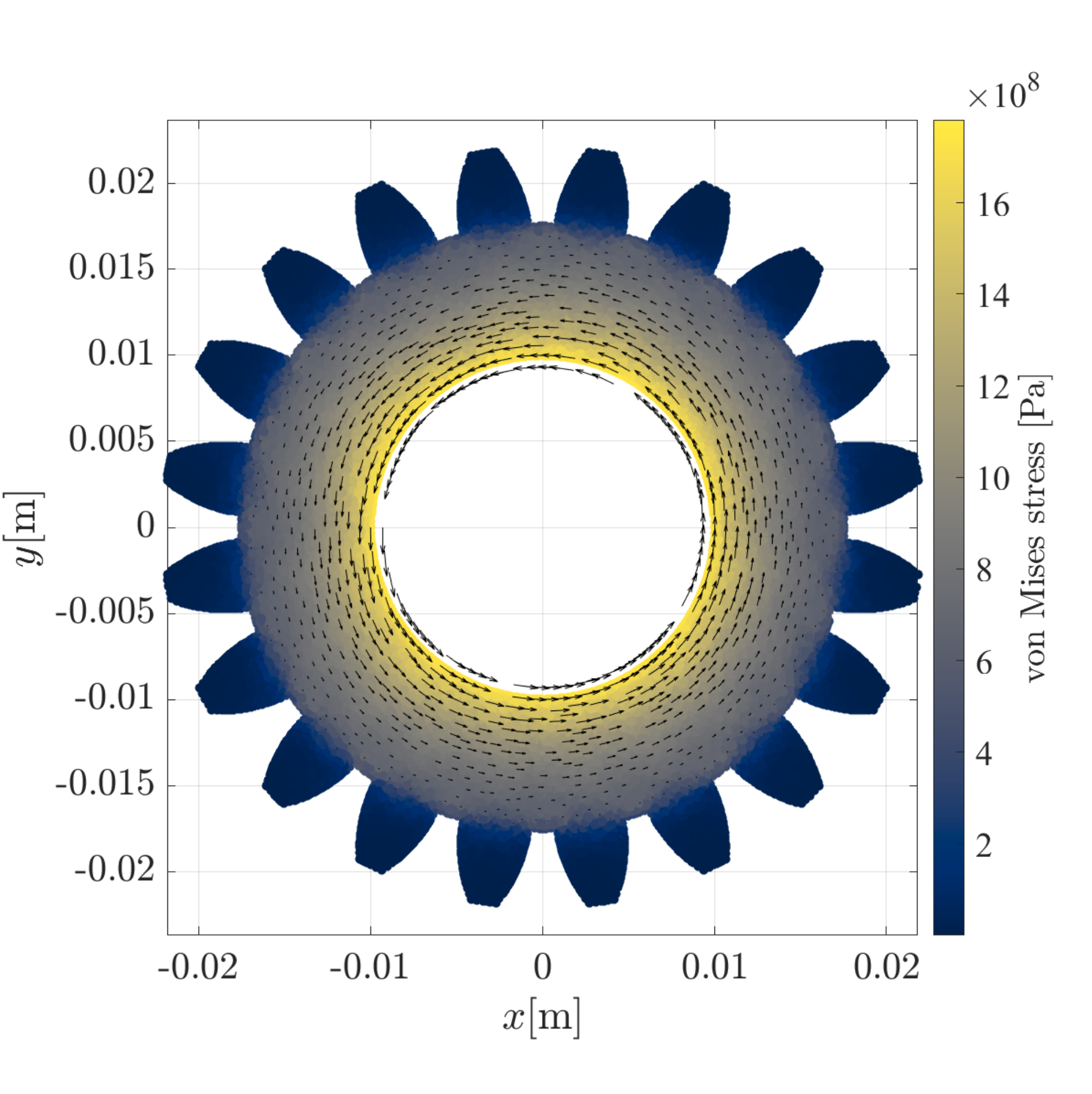}
	\includegraphics[trim={0cm 1cm 0cm 1cm},clip,  height=0.40\linewidth]{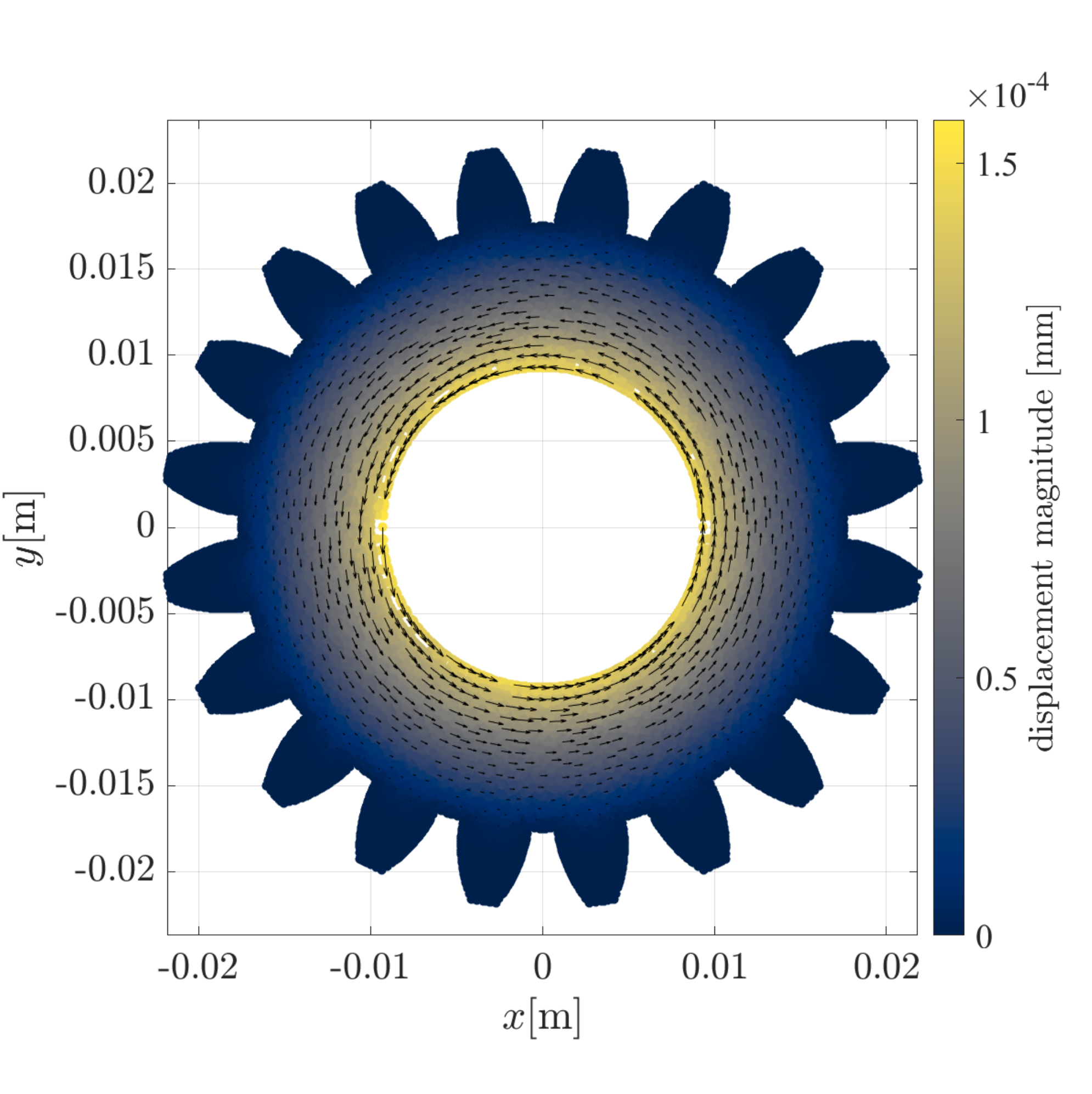}
	
	\caption{The von Mises stress (left) and the displacement magnitude (right) at $z=0$ cross section accompanied with a quiver plot of the displacement field.}

	\label{fig:gear_cross}
\end{figure}

In \diff{the} previous section, we established confidence in the presented solution procedure by obtaining expected convergence rates in solving Poisson's equation on three different geometries in 2D and 3D. In this section we apply NURBS-DIVG to a more realistic case from linear elasticity, governed by the Navier-Cauchy equation

\begin{equation}
	\frac{E}{2\left(\nu  + 1 \right)}\left( \nabla ^2 \mathbf{u} + \frac{1}{1 - 2 \nu }\nabla \left( {\nabla  \cdot \mathbf{u}} \right) \right) = 0
\end{equation}   
where $\b u$ stands for the displacement vector, and Young's modulus $E =  \SI{72.1e9}{\pascal}$ and Poisson's ratio $\nu = 0.33$ define material properties. The displacement and the stress tensor ($\sigma$) are related via Hooke's law  
\begin{equation}
	\label{eq:elasticity}
	\sigma = \frac{E}{\nu  + 1} \left(    
	\frac{1}{1 - 2 \nu } \mathrm{tr}(\varepsilon) I + \varepsilon \right), \quad
	\varepsilon = \frac{\nabla \bf{u} + (\nabla \mathbf{u})^\mathsf{T}}{2},
\end{equation}
with $\varepsilon$ and $I$ standing for strain and identity tensors. We observe
a 3D gear object that is subjected to an external torque resulting in a
tangential traction $t_0 = \SI{1e3}{\pascal}$ on axis, while the gear teeth are
blocked, \emph{i.e.} the displacement is zero $\mathbf u = \SI{0}{\metre}$.
\diff{The} top and bottom surfaces are free, \emph{i.e.} traction free boundary
conditions apply. In summary
\begin{align}
	\bf{u} &= \SI{0}{\metre}, \quad \text{on } \Gamma_{\text{teeth}},\\
	\sigma\cdot \mathbf{n} &= \SI{0}{\pascal}, \quad \text{on } \Gamma_{\text{free}}, \\
	\sigma\cdot \mathbf{t} &= t_0 \mathbf{t}, \quad \text{on } \Gamma_{\text{axis}}.
\end{align}
The case is schematically presented in Figure~\ref{fig:gear_scheme} together with von Mises stress scatter plot. The stress is highest near the axis where the force is applied, and gradually fades towards blocked gear teeth. 
Displacement and von Mises stress are further demonstrated in Figure~\ref{fig:gear_cross} at $z = \SI{0}{\metre}$ cross section, where we see how the gear is deformed due to the applied force. \diff{All results} were computed using 41210 scattered nodes generated by NURBS-DIVG.   

\subsection{Transient heat transport}
\diff{The} last example is focused on \diff{the} transient heat equation
\begin{equation}
	\label{eq:heat}
	\frac{\partial T}{\partial t} = \lambda \nabla^2 T + q,
\end{equation}
\diff{where $T$ stands for temperature, $\lambda$ for thermal conductivity, and $q$ for the heat source. The goal is to solve heat transport within the duck model subject to the Robin boundary condition}
\begin{equation}
	\frac{\partial T}{\partial t} + T = 0 
\end{equation} 
and a heat source within the domain
\begin{equation}
	q = 5 e^{10||\mathbf x - \mathbf x_0||}
\end{equation} 
with $\mathbf x_0 = (0,0,0.2)$, \diff{the}
initial temperature set to $0$ throughout the domain, \diff{ and $\lambda = 2$.} \diff{Time marching is performed via implicit stepping}
\begin{equation}
	\label{eq:disc}
	\frac{T_2 - T_1}{\Delta t} = \lambda \nabla^2 T_2 + q,
\end{equation}
\diff{where $T_1$ and $T_2$ stand for the temperature in the current and the
  next time step respectively and $\Delta t$ represents the time step. The
  spatial discretization of the Laplace operator is done using RBF-FD with
  $m=2$. We used a time step of $\Delta t = 3 \cdot 10^{-4}$ and $3000$ iterations
to reach the steady state using the criterion $T_2 - T_1 < 3\cdot 10^{-6}$) at
$t=0.9$.}

\begin{figure}[h]
	\centering
	\includegraphics[width=0.45\linewidth]{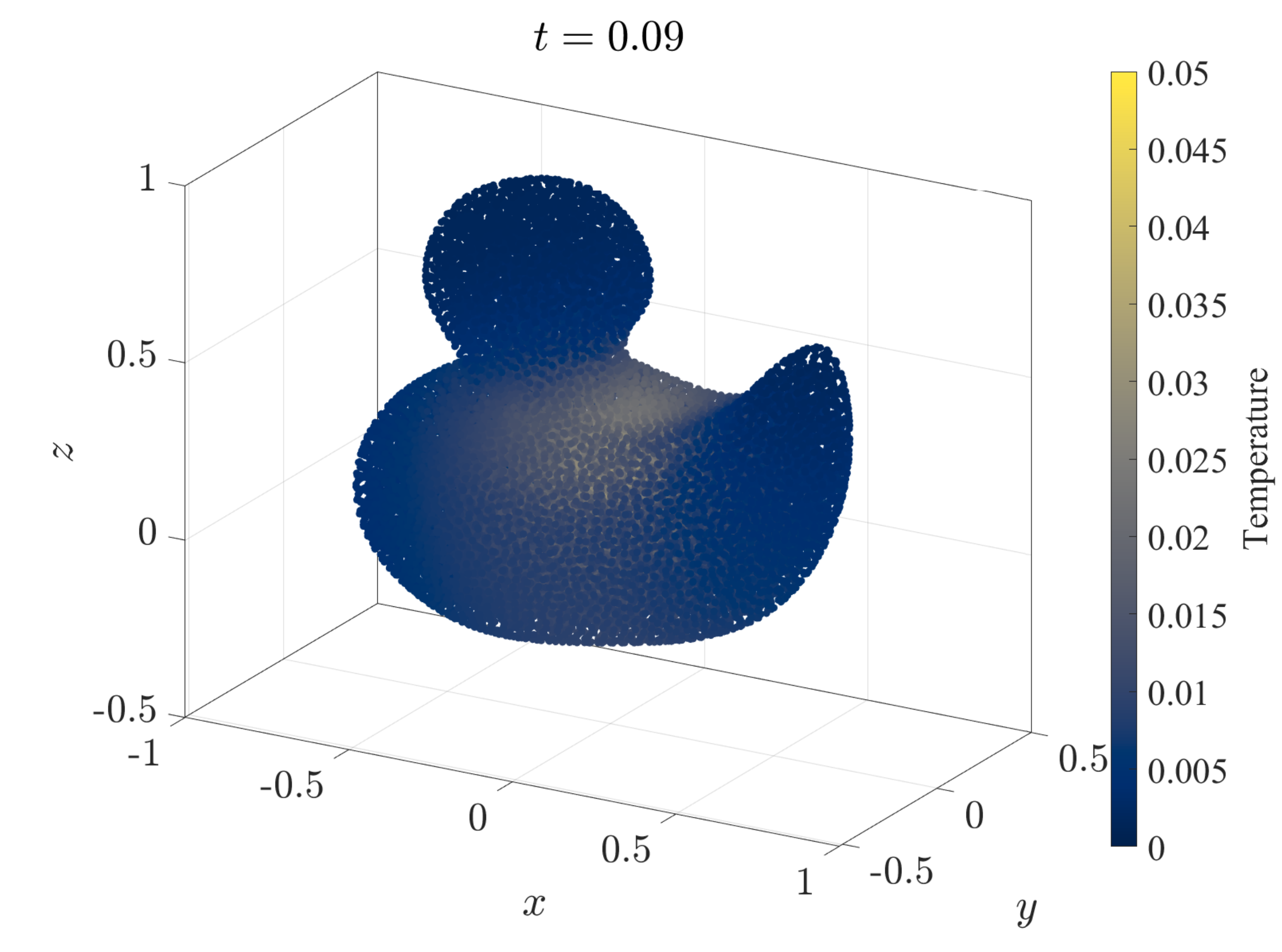}
	\includegraphics[width=0.45\linewidth]{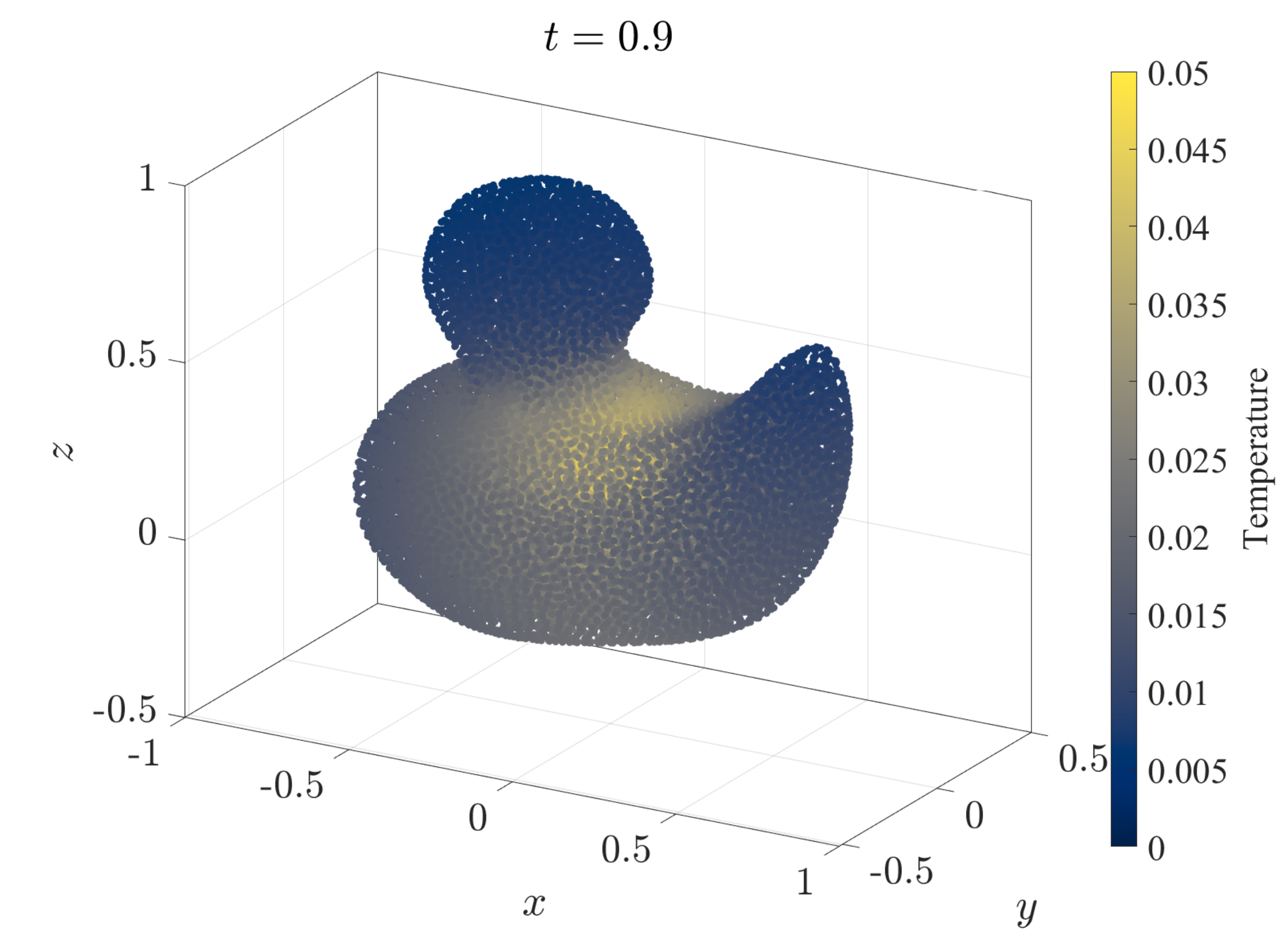}
  \caption{Heat transport within a 3D duck. \diff{The model is based on}~\cite{tdduckmodel} and consists of only 1 patch.}
	\label{fig:heat}
\end{figure}

\begin{figure}[h]
	\centering
	\includegraphics[trim={3cm 14cm 3cm 5cm},clip,  width=1\linewidth]{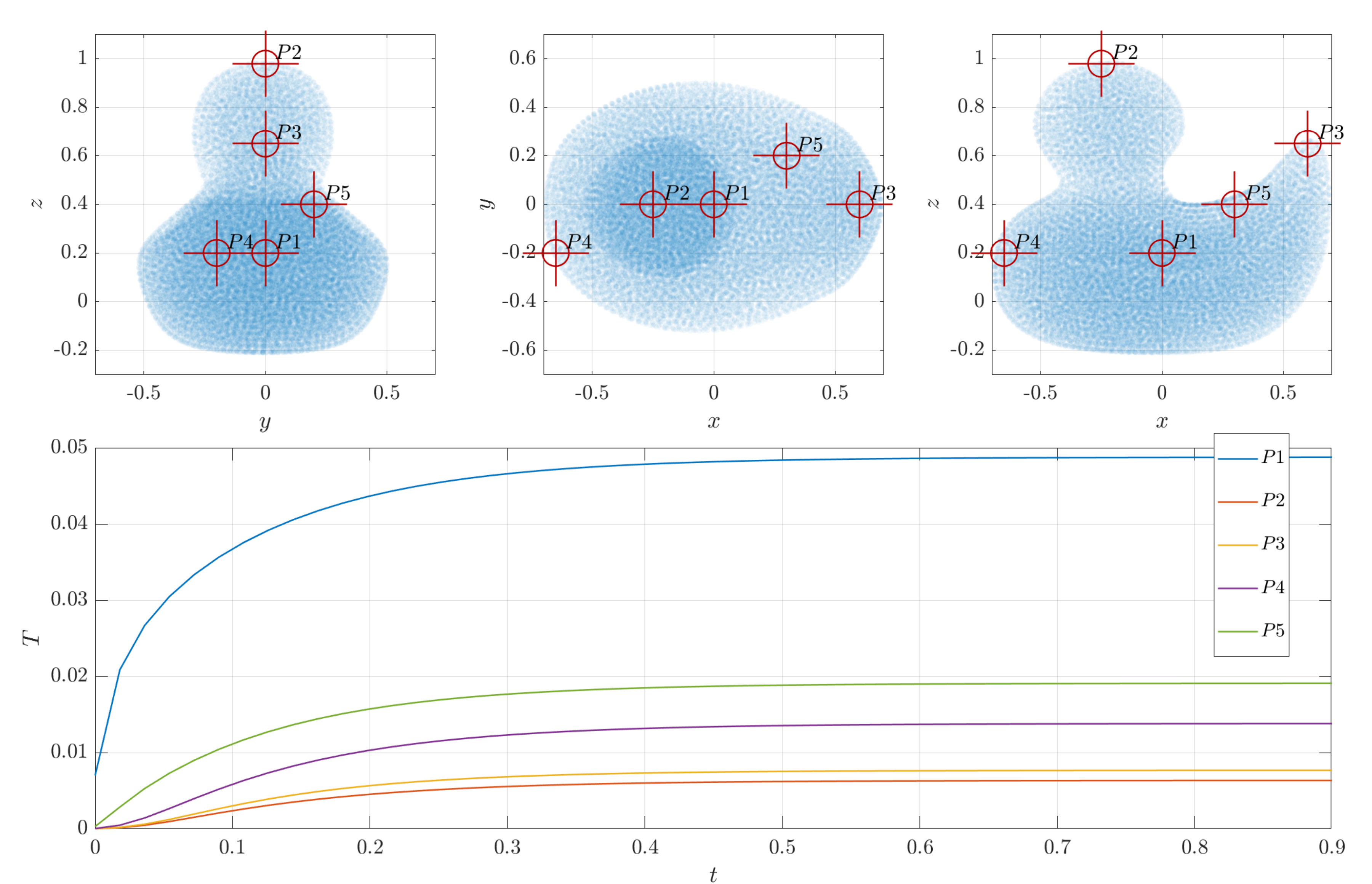}

	\caption{Time evolution of the temperature at five control points.}
	\label{fig:heat_time}
\end{figure}

Figure~\ref{fig:heat} shows the temperature scatter plot computed with RBF-FD on 21956 nodes generated with the proposed NURBS-DIVG at two different times (first at the beginning of the simulation and second at the steady state). 
In Figure~\ref{fig:heat_time}, the time evolution of the temperature at five control points $P1-5$ is shown. Control point $P1$ is located at the heat source, $P2-4$ at \diff{the} most distant points from the source, and $P5$ asymmetric with respect to the y-axis. As one would expect, at the source the temperature rises immediately after the beginning of the simulation and also reaches the highest value, while the rise is a bit delayed and lower at the distant points that are closer to the boundary where the heat exchange with surroundings takes place. Once the heat exchange with \diff{the surroundings} matches the heat generation at the source, the system reaches steady-state. 

\section{Conclusions}
\label{sec:conc}
In this paper, we presented a meshless algorithm, NURBS-DIVG, for generating
quasi-uniform nodes on domains whose boundaries are defined by CAD models
consisting of multiple NURBS patches. The NURBS-DIVG algorithm is able to deal with complex geometries with sharp edges and concavities, supports refinement, and
can be generalized to higher dimensions. We also demonstrated that node layouts
generated with NURBS-DIVG are of sufficiently high quality for meshless
discretizations, first by directly assessing the quality of these node sets,
then by using RBF-FD to solve the Poisson equation with mixed Dirichlet-Neumann
boundary conditions on different domains to high-order accuracy. Finally, we
demonstrated NURBS-DIVG in conjunction with RBF-FD in tackling two more
challenging test cases: first, the stress analysis of a gear subjected to an
external force governed by the Navier-Cauchy equation; and second, a
time-dependent heat transport problem inside a duck. This work advances the
state of the art in fully-autonomous, meshless, isogeometric analysis. All
algorithms presented in this work are implemented in C++ and included in our
in-house open-source meshfree library \emph{Medusa}~\cite{medusa,slak2021medusa}, see the \emph{Medusa wiki}~\cite{medusawiki} for usage
examples. The interface to all CAD files was implemented via Open
Cascade~\cite{opencascade}.

\begin{acknowledgements}
The first and third authors acknowledge the financial support from the Slovenian Research Agency research core funding No. P2-0095, research project J2-3048, and research project N2-0275. The second author was partially supported by the United States National Science Foundation (NSF) grant CISE CCF 1714844.

Funded by National Science Centre, Poland under the OPUS call in the Weave programme 2021/43/I/ST3/00228. This research was funded in whole or in part by National Science Centre (2021/43/I/ST3/00228). 

For the purpose of Open Access, the author has applied a CC-BY public copyright licence to any Author Accepted Manuscript (AAM) version arising from this submission.
\end{acknowledgements}

\section*{Declarations}
The authors declare that they have no conflict of interest.

\section*{Data availability}
The datasets generated during and/or analysed during the current study are available from the corresponding author on reasonable request. For some practical examples see~\cite{medusawiki}.

\bibliography{references}
\bibliographystyle{spmpsci}
\end{document}